\documentclass[12pt]{article}
\usepackage{amsmath}
\usepackage{graphicx}
\usepackage{enumerate}
\usepackage{natbib}
\usepackage{url} 

\newcommand{\blind}{0}

\RequirePackage{amsthm,amsmath,amsfonts,amssymb,subcaption,booktabs,multirow}
\usepackage{pdflscape}
\usepackage{floatrow}
\usepackage{makecell}
\usepackage{threeparttable}
\allowdisplaybreaks
\RequirePackage[colorlinks,citecolor=blue,urlcolor=blue]{hyperref}
\RequirePackage{graphicx}
\theoremstyle{plain}

\newtheorem{theorem}{Theorem}[section]
\newtheorem{corollary}{Corollary}[section]

\newtheorem{remark}{Remark}[section]
\newtheorem{definition}{Definition}[section]
\newtheorem{assumption}{Assumption}
\theoremstyle{remark}

\newcommand{\fo}{f^{\circ}}
\newcommand{\go}{g^{\circ}}
\newcommand{\xx}{x}
\newcommand{\x}{X}
\newcommand{\bs}{\boldsymbol{s}}
\newcommand{\bk}{\boldsymbol{k}}
\newcommand{\bl}{\boldsymbol{l}}
\newcommand{\bu}{\boldsymbol{u}}
\newcommand{\bo}{\boldsymbol{0}}
\newcommand{\bv}{\boldsymbol{v}}
\newcommand{\ba}{\boldsymbol{a}}
\newcommand{\bka}{\boldsymbol{\kappa}}
\newcommand{\cf}{\int_{\Omega}(f-\fo)^2d\mathcal{P}_{\x}}
\newcommand{\chf}{\int_{\Omega}(\hat f-\fo)^2d\mathcal{P}_{\x}}
\newcommand{\ctf}{\int_{\Omega}(\tilde f-\fo)^2d\mathcal{P}_{\x}}

\newcommand{\cehf}{\sum_{i=1}^{n}\sum_{j=1}^{m}\left( \hat f(\x_{ij})-\fo(\x_{ij})\right)^2}
\newcommand{\cetf}{\sum_{i=1}^{n}\sum_{j=1}^{m}\left( \tilde  f(\x_{ij})-\fo(\x_{ij})\right)^2}
\newcommand{\var}{\operatorname{Var}}
\newcommand{\cov}{\operatorname{Cov}}
\newcommand{\nn}{\mathcal{FNN}}
\newcommand{\sm}{s_{\text{min} }}
\newcommand{\sM}{s_{\text{max} }}
\newcommand{\relu}{\operatorname{ReLU}}
\newcommand{\ham}{\operatorname{Ham}}
\newcommand{\kl}{\operatorname{KL}}
\newcommand{\diag}{\operatorname{diag}}
\usepackage{tikz}
\newcommand{\ts}{\tilde{s}}
\newcommand{\sso}{s_{0}}
\newcommand{\ddo}{d_{0}}
\newcommand{\ds}{\varrho}
\newcommand{\pef}{\mathcal{E}_{\fo}(\hat{f}_{nm})}

\newcommand{\pefh}{\mathcal{E}_{\fo}(\hat{f})}
\newcommand{\delf}{\Delta_{nm}(\hat{f}_{nm})}
\newcommand{\vc}{\operatorname{VCdim}}
\newcommand{\hmfhf}{\hat{\mathcal{E}}_{\fo}(\hat{f})}

\newcounter{cnstcnt}
\newcommand{\xc}{%
   \refstepcounter{cnstcnt}%
   \ensuremath{c_{\thecnstcnt}}
}
\newcommand{\jc}[1]{\ensuremath{c_{\ref{#1}}}}

\newcounter{bnstcnt}
\newcommand{\xb}{%
   \refstepcounter{bnstcnt}%
   \ensuremath{b_{\thebnstcnt}}
}
\newcommand{\jb}[1]{\ensuremath{b_{\ref{#1}}}}

\newcounter{dbnstcnt}
\newcommand{\xdb}{%
   \refstepcounter{dbnstcnt}%
   \ensuremath{B_{\thedbnstcnt}}
}
\newcommand{\jdb}[1]{\ensuremath{B_{\ref{#1}}}}

\usepackage{enumerate}

\begin{document}

\def\spacingset#1{\renewcommand{\baselinestretch}%
{#1}\small\normalsize} \spacingset{1}


	\if0\blind
	{
		\title{\bf Deep Regression for Repeated Measurements}
		\author{
			Shunxing Yan 
                \thanks{Shunxing Yan is the first author,  E-mail: sxyan@stu.pku.edu.cn}
                \\    
			School of Mathematical Sciences, Center for Statistical Science, \\
			Peking University, Beijing, China\\ 
			Fang Yao
			\thanks{
				Fang Yao is the corresponding author, E-mail: fyao@math.pku.edu.cn.
			}
            \hspace{.2cm} 
			\\
			School of Mathematical Sciences, Center for Statistical Science, \\
			Peking University, Beijing, China\\ 
                Hang Zhou 
                \thanks{hgzhou@ucdavis.edu}
                \\Department of Statistics, University of California Davis, USA
		}  
		\date{}
		\maketitle
	} \fi

\if1\blind
{
  \bigskip
  \bigskip
  \bigskip
  \begin{center}
    {\LARGE \bf Deep Regression for Repeated Measurements}
\end{center}
  \medskip
} \fi

\bigskip
\begin{abstract}
Nonparametric mean function regression with repeated measurements serves as a cornerstone for many statistical branches, such as longitudinal/panel/functional data analysis.  
In this work, we investigate this problem using fully connected deep neural network (DNN) estimators with flexible shapes.
A novel theoretical framework allowing arbitrary sampling frequency is established by adopting empirical process techniques to tackle clustered dependence. 
We then consider the DNN estimators for H\"older target function and illustrate a key phenomenon, the phase transition in the convergence rate, inherent to repeated measurements and its connection to the curse of dimensionality. 
Furthermore, we study several examples with low intrinsic dimensions, including the hierarchical composition model, low-dimensional support set and anisotropic H\"older smoothness. We also obtain new approximation results and matching lower bounds to demonstrate the adaptivity of the DNN estimators for circumventing the curse of dimensionality. 
Simulations and real data examples are provided to support our theoretical findings and practical implications. 
\end{abstract}

\noindent%
\vfill

\newpage
\spacingset{1.9} 

\section{Introduction}
\subsection{Literature Review}
Repeated measurements collected from a sample of subjects have been extensively investigated in statistics, finding applications across fields like biomedicine, epidemiology, economics, and engineering. 
These observations typically possess dependence within subjects while being independent across subjects, also referred to as ``clustered dependence''.
Consequently, several statistical branches have emerged, including longitudinal data analysis \citep{diggle2002analysis,weiss2005modeling,hedeker2006longitudinal},  panel data analysis \citep{chamberlain1984panel, baltagi2008econometric, hsiao2022analysis}, and functional data analysis \citep{ramsay2005,ferraty2006nonparametric,hsing2015theoretical}.
Among these, estimating the mean function is a fundamental problem and serves as a crucial step for subsequent analyses. 

Classical nonparametric regression methods \citep{rice1983smoothing, fan1992variable} have addressed this for longitudinal data that typically have finite measurements per subject \citep[i.e., finite sampling frequencies,][]{brumback1998smoothing, lin2000nonparametric}, while cases that allow the sampling frequency to infinity were established in the context of functional data by \citet{cai2011optimal} and \citet{zhang2016sparse}. 
This line of research on traditional methods spanning over 20 years reflects the challenges and significance of establishing a unified method with theory for various sampling schemes, especially exploring the effect of sampling frequencies and indicating how densely such measurements are needed to achieve the optimal rate of convergence. 
Moreover, traditional nonparametric methods are usually for univariate or low-dimensional covariates and suffer from the curse of dimensionality. 
Hence, estimators based on modern tools like neural networks are desired, given the success achieved by deep learning in recent years.

Deep neural networks (DNNs) with ReLU activation function have rapidly received increasing attention with a vast literature emerging, due to its computational power and approximation ability \citep{telgarsky2016benefits,fan2021selective}. 
\citet{yarotsky2017error,schmidtaos,yarotsky2018optimal, kohler2021rate} have established non-asymptotic approximation bounds for sparsely connected or more implementable ReLU DNNs.
Recently, \citet{shen2020}, \citet{lu2021deep} and  \citet{shen2022optimal} established the (nearly) optimal rate of approximation of fully connected neural networks in terms of both width and depth. 
A crucial preponderance of DNNs is the adaptability to various low intrinsic-dimensional structures that may circumvent curse of dimensionality. 
Without specific structure knowledge, DNNs automatically approximate the composition models \citep{bauer2019deep, schmidtaos, kohler2021rate} that contain the common semiparametric models, and various smoothness settings \citep{barron1993universal,suzuki2018adaptivity,suzuki2021deep}. 
Another perspective is that the predictor lies on a low-dimensional subset.
\citet{shaham2018provable} used the 4-layer networks to approximate functions with suitable smoothness on low-dimensional manifolds.
\citet{schmidt2019deep,chen2019efficient,nakada2020adaptive} considered DNN estimators for more general smooth functions with manifold or Minkowski dimension sets input, and \citet{cloninger2021deeparxiv} showed the adaptivity of DNN estimator to intrinsic dimensionality even beyond the domain for some cases.
A comprehensive theory for fully connected DNN estimators of low-dimensional inputs was then derived in \citet{huangjian}. 
Recently, \citet{zhang2023effective} introduced a novel effective Minkowski dimension determining the regression rate.

Benefiting from the aforementioned advantages, DNNs are gaining attention as an emerging nonparametric approach in various statistical problems. 
There has been a series of work that successfully established theories of nonparametric DNN regression for independent observations \citep{bauer2019deep,schmidtaos,suzuki2018adaptivity,farrell2021deep,kohler2021rate}, which here is referred to as {\it cross-sectional} model/setting that measures a single observation for each subject and is distinguished from repeated measurements setting. 
\citet{covnet} proposed a novel CovNet to estimate the covariance function for fully and grid-observed functional data, which approximates well and facilitates a direct eigendecomposition.
\citet{fan2022noise} studied the Huber DNN estimator for heavy-tail and \citet{deephazard,deepcox} considered survival models. 
Recently, some noteworthy works \citep{fan2022factor, bhattacharya2023deep} considered the high dimensional problem, including factor augmented sparse throughput model and nonparametric interaction models with rigorous theory and matching lower bounds.

\subsection{Challenges of DNN Regression for Repeated Measurements}
Compared with the cross-sectional setting, the repeated measurements model is composed of $n$ subjects (often referred to as sample size), and the $i$th individual is measured $m_i$ times. 
Here we assume $m_i\asymp m$ for convenience, where $a_n\lesssim(\gtrsim) b_n$ indicates $a_n\leq(\geq) C b_n$ for some constant $C>0$, and $a_n\asymp b_n$ means $a_n\lesssim b_n$ and $b_n\lesssim a_n$.  
To estimate the mean function, we adopt the pooling strategy that combines observations from all subjects. It works well for both sparse and dense designs and has been widely 
used in many contexts, such as random effects models, generalized estimating equation approaches  \citep{zeger1988models,lin2000nonparametric}, and functional data models   \citep{cai2011optimal,zhang2016sparse}, among others.

Along with broad applicability and excellent properties of pooling repeated measurements, the clustered dependence structure, especially $m$ can be finite or growing with the sample size $n$ at any rate, brings theoretical difficulties,  which renders the techniques and results of the DNN regression based on independent observations \citep{gyorfi2002distribution,koltchinskii2006local}.
Different from the kernel methods with closed-form estimates and spline/wavelet based methods with specific-structure basis, the DNN suffers greater challenges. 
Recently, \citet{wang2021estimation} studied this problem without generalization error, while a key assumption on the vanishing maximal eigenvalue of covariance matrix (as $m \rightarrow \infty$) contradicts its positive lower bound. 
Therefore, a new theoretical framework is desired for clustered dependent data to which existing empirical process results are not applicable.

Similar to the cross-sectional setting, 
the repeated measurements model also suffers from the curse of dimensionality. 
However, the latter behaves differently due to an interesting phenomenon, namely the {\it phase transition} in the sense that the convergence rate changes when the number of measurements grows from sparse to dense sampling schemes. 
Here, the curse of dimensionality is accompanied by larger errors under sparse design and the requirement for higher sampling frequency to achieve a faster and even a parametric rate.
We shall study this phenomenon for DNN estimators in multi-dimensional and complex clustered structures, which would guide the practical sampling design for improving estimation efficiency in statistical scenarios.

One of our objectives is to circumvent the curse of dimensionality for repeated measurements as in cross-sectional settings. 
For this purpose, we need accurate characterizations of the DNN approximation ability in various function spaces that are not fully understood in the existing literature. 
For instance, the hierarchical composition model used in \citet{kohler2021rate} requires that the smoothness of each function is no less than one and the neural network has a specific shape, which limits its applicability. Moreover, there are no results available on the approximation of fully connected ReLU networks to anisotropic smoothness functions. These problems urge us to develop some new approximation results that are meaningful for both cross-sectional and repeated measurements settings.

\subsection{Our Contributions}
Motivated by the aforementioned challenges and issues, we establish a comprehensive framework for DNN regression in estimating mean functions for repeated measurements. In particular,
we exploit a conditioning argument to decouple the randomness between and within subjects, which is the first attempt to adopt and thus allow the use of empirical process tools for general clustered dependence. 
This can also be utilized for further theoretical analysis in relevant models and problems. 
The main contributions of this paper are summarized as follows.
\begin{enumerate}[(1)]

    \item We develop a series of oracle bounds on prediction error for estimators based on minimizing the empirical least squares loss function in the repeated measurements setting. 
    To the best of our knowledge, it is the first work in this setting and may serve as a cornerstone for developing theories of various statistical models coupled with DNNs.
    
    \item We derive the upper bound on the prediction error using fully connected DNNs (and also regression splines as a by-product) when the true function belongs to the H\"older space, as well as the matching minimax lower bound.
    It follows that a phase transition phenomenon occurs for DNN estimators, which reveals the impact of dimensionality on the convergence rate in a discrepant way with cross-sectional models.    
    \item To circumvent the curse of dimensionality, we study several representative low intrinsic-dimensional structures, including the hierarchical composition model, low dimensional input and anisotropic H\"older smoothness, and analyze the minimax optimal rates and phase transition phenomena, respectively. 
    We also develop some new results having their own merits even for cross-sectional settings, for example, representing the hierarchical composition model to allow the smoothness smaller than one and deriving a new approximation result for anisotropic H\"older smooth functions.
\end{enumerate}

{\spacingset{1.2} 
\begin{table}[htbp]
  \caption{A summary of convergence rates for ReLU DNN estimators for various function classes and sampling frequencies. All rates are minimax optimal.} \label{tab:addlabel}
  \begin{threeparttable}   
    \centering \scriptsize
    \begin{tabular}{ccccc}
    \toprule
    \multicolumn{1}{c}{\multirow{2}[4]{*}{Function Class}} & \multicolumn{1}{c}{\multirow{2}[4]{*}{ 
    \makecell[c]{ Cross-Sectional\\Regression $(m=1)$}
    }} & \multicolumn{3}{c}{Repeated Measures Regression $(m \geq 2)$} \\
\cmidrule{3-5}   &   & \multicolumn{1}{c}{Sparse}    & \multicolumn{1}{c}{Phase transition}  & \multicolumn{1}{c}{Dense} \\
    \midrule
    H\"older &  \makecell[c]{ $n^{-2s/(2s+d)}$ \\  \citep{kohler2021rate} }  & $(nm)^{-2s/(2s+d)}$ & $m \asymp n^{d/2s}$  & $n^{-1}$ \\
    \midrule
    Hierarchical & \makecell[c]{ $n^{- 2\gamma/(2\gamma+1)}$ \\ \citep{kohler2021rate} }  & $(nm)^{- 2\gamma/(2\gamma+1)}$ & $m \asymp n^{1/2\gamma}$ &$n^{-1}$ \\
    \midrule
    Manifold Supp & \makecell[c]{  $n^{-2s/(2s+d_{\mathcal{M}})}$ \\ \citep{huangjian}} & $(nm)^{-2s/(2s+d_{\mathcal{M}})}$ & $m \asymp n^{d_{\mathcal{M}}/2s}$  & $n^{-1}$ \\ 
    \bottomrule
    Anisotropic &  \makecell[c]{  $n^{-2\ts/(2\ts+d)}$ \\ \citep[sparse net][]{suzuki2021deep}}  
    & $(nm)^{-2\ts/(2\ts+d)}$ & $m \asymp n^{d/2\ts}$  & $n^{-1}$ \\
    \midrule
    \end{tabular}%
\begin{tablenotes}    
\scriptsize
\item[a] Rates for $m\geq 2$ are established in this work. Here, ``Phase transition" indicates the order when phase transition in the convergence rate occurs. ``Sparse" refers to sampling frequencies $m$ sparser than the order where phase transition occurs, while ``Dense" pertains to those denser. 
\end{tablenotes}           
\end{threeparttable} 
\end{table}%
}
\vspace{-1.2cm}  
\subsection{Organization}
The rest of the article is organized as follows. 
The model settings, estimators and key oracle inequalities are presented in Section \ref{sec:pre}, based on which 
Section \ref{sec:holderph} presents the convergence analysis for H\"older space. 
We then provide in Section \ref{sec:cir} two useful examples to circumvent the curse of dimensionality. 
Section \ref{sec:sim} illustrates the theoretical findings via simulation experiments and Section \ref{sec:realdata} presents two real data examples, while further discussion is given in  Section \ref{sec:con}. 
The additional theoretical and numerical results, all proofs and auxiliary lemmas are deferred to the Supplementary Material for space economy.
\section{Model Estimation and Oracle Inequalities}\label{sec:pre}
\subsection{Model and Estimation}

Given a random function $Z(\xx)$ defined on $\left[0,1\right]^{d}$. 
Let $Z_{1}(\xx), Z_{2}(\xx), ..., Z_{n}(\xx)$ be independent and identically distributed (i.i.d.) realizations. 
In practice, these realizations cannot be fully observed and are measured for the $i$th subject at $m_{i}$ points as ${\x_{i1}, \x_{i2}, ...,\x_{im_{i}}}$, inevitably contaminated with additive noise.
Under random design assumption, we assume that $\x_{ij}, 1 \leq i \leq n, 1\leq j \leq m_i$ are i.i.d. from a distribution $\mathcal{P}_{\x}$ which is supported on some $\Omega \subseteq \left[0,1\right]^{d}$ and independent of $Z_{i}, 1 \leq i \leq n$.  
Here, $\Omega$ is assumed to have a well-defined volume measure and $\mathcal{P}_{\x}$ is absolutely continuous with respect to the 
volume measure. 
Let $\fo(\xx)=\mathbb{E}\left[ Z(\x)|\x=\xx \right]$ be the regression function (i.e., the mean structure) of interest. 
Then the repeated measurement responses follow  
\begin{equation}\label{equ:model}
    Y_{ij}=Z_{i}(\x_{ij})+\epsilon_{ij}=\fo(\x_{ij})+U_{i}(\x_{ij})+\epsilon_{ij}, \hspace{0.2in}  1\leq i \leq n, 1\leq j \leq m_i,
\end{equation}
where $U_{i}(\xx)=Z_{i}(\xx)-\fo(\xx)$ are individual stochastic parts and the $\epsilon_{ij}$ are independent noise variables. 
A key feature of this repeated measurements model is that observations from the same subject are correlated/dependent, while those from different subjects are independent, which is distinct from the cross-sectional setting based on independent observations.
In the sequel, to simplify exposition, we assume that each subject is measured $m$ times, i.e.,  $m_1=...=m_{n}=m$ (also referred to as sampling frequency). 

With the available data $\left\{ \left(Y_{ij},\x_{ij} \right): 1\leq i \leq n, 1\leq j \leq m \right\}$, our target is to estimate the mean function $\fo(\xx)$. 
We pool observations from all subjects for estimation, and define the nonparametric least square estimator within a suitably chosen function class  $\mathcal{F}_{nm}$, minimizing the empirical risk as follows, 
\begin{equation}\label{equ:estimator}
    \hat{f}_{nm} \in \arg\min_{f\in \mathcal{F}_{nm}}\frac{1}{nm}\sum_{i=1}^{n}\sum_{j=1}^{m}\left(Y_{ij}-f(X_{ij})\right)^{2}. 
\end{equation}
To evaluate the accuracy for an estimator $\hat{f}_{nm}$ quantitatively, we introduce the mean squared prediction error
$$
\pef:=\mathbb{E}\left[\int_{\Omega}(\hat{f}_{nm}(\xx)-\fo(\xx))^2d\mathcal{P}_{\x}\right], 
$$
which is also equal to the excess risk 
$\mathbb{E}[|Y-\hat{f}_{nm}(\x)|^{2}]-\mathbb{E}[\left|Y-\fo(\x)\right|^{2}]$. 

In \eqref{equ:estimator}, the candidate class $\mathcal{F}_{nm}$ plays a central role in estimation procedures for the trade-off of approximation error and estimation error.  
In this paper, we mainly consider the neural network class. 
A ReLU feedforward neural network with architecture $(L,(W_{L}, W_{L-1},..., W_{0}))$ can be written as a composition of a series of functions
\begin{equation}\label{equ:fnn0}
f(\xx)=\mathcal{L}_{L} \circ  \relu  \circ \mathcal{L}_{L-1} \circ  \relu  \circ \mathcal{L}_{L-2} \circ  \relu  \circ \cdots \circ \mathcal{L}_{1} \circ  \relu  \circ \mathcal{L}_{0}(\xx),    
\end{equation}
where $\mathcal{L}_{i}(\xx)=L_{i} x+b_{i}$ is a multivariate linear function with $L_{i} \in \mathbb{R}^{W_{i+1} \times W_{i}}, b_{i} \in \mathbb{R}^{W_{i}}$, $x\in \mathbb{R}^{W_{0}}$, and the operator $\relu:\mathbb{R}^{W_{i}}\rightarrow \mathbb{R}^{W_{i}}$ 
applies rectified linear unit (ReLU) function $a\mapsto \max \{a, 0\}$ to each component of its input. 
Then we define 
$$
\nn (d, L, W)= \left\{f: \mathbb{R}^{d} \rightarrow \mathbb{R} \text { is of the form \eqref{equ:fnn0} with } W_{0}=d \text{ and } W_{i}\leq W \text{, } 1\leq i\leq L \right\},
$$
 and 
$$
\nn (d, L, W, \beta)=\mathcal{T}_{\beta} \nn (d, L, W) =\left\{\mathcal{T}_{M}f: f \in \mathcal{T}_{\beta} \nn (d, L, W)\right\},
$$ 
where 
$\mathcal{T}_{\beta}f(\xx)=f(\xx)I(|f(\xx)|\leq\beta)+\operatorname{sign}(f(\xx))\beta I(|f(\xx)|>\beta)$ is the truncation operator. 
Specifically, we call $L$ and $W$ the depth and width of the neural networks, respectively.

\subsection{Oracle Inequalities}
We first establish some universal bounds on the mean squared prediction error of the proposed estimators without imposing specific $\mathcal{F}_{nm}$. 
Denote the infinity norm $\|f\|_\infty=\sup_{x\in \Omega} |f(x)|$ for a function $f$ defined on $\Omega$. 
We make the following assumptions. 

\begin{assumption}[Regression function and candidate class]\label{ass:fobound} There exist some positive number $\xdb\label{b1}$, such that the target function $\|\fo\|_{\infty}\leq \jdb{b1}$ and 
$\|f\|_{\infty}\leq \jdb{b1}$ for each $f\in \mathcal{F}_{nm}$.
\end{assumption}

\begin{assumption}[Random process and noise] \label{ass:subexp}
The random processes $U_{i}(x)$ are continuous, and there exist nonnegative numbers $\xdb\label{b2}$ and $\xdb\label{b3}$, such that for $1\le i \le n, 1\le j \le m$,
\begin{equation}\label{equ:proexp}
    \mathbb{E}\left[\exp\left\{\left|\frac{ U_i(\x_{ij})  }{\jdb{b2}} \right|\right\}\right]\leq 1, \hspace{0.2in} 
    \mathbb{E}\left[\exp\left\{\left|\frac{ \epsilon_{ij}}{\jdb{b3}}\right| \right\}\right]\leq 1.
\end{equation}
\end{assumption}

Assumption \ref{ass:fobound} is a standard condition in empirical process and nonparametric regression. 
For example, it guarantees a probability bound of $\pef$
by its empirical version. 
To satisfy this assumption, the candidate functions are often truncated. 
Remarkably, our later results are non-asymptotic with respect to $\jdb{b1}$, so it can be relaxed as an increasing sequence with respect to $mn$ when the exact bound of $\|\fo\|_{\infty}$ is unknown. 
Assumption \ref{ass:subexp} that requires sub-exponential random parts is also standard.

Now we specify a condition to characterize the complexity of $\mathcal{F}_{nm}$. Define  
$\mathcal{F}_{nm}(r)=\left\{ f\in \mathcal{F}_{nm}: \int_{\Omega}\left( f(\xx)-\fo(\xx) \right)^2 d\mathcal{P}_{\x} \leq r  \right\}. $
We say a function $\phi_{nm}(r)$ is sub-root if and only if  $\phi_{nm}(r)/\sqrt{r}$ is nonnegative and nonincreasing. Let the i.i.d. Rademacher variables $\left\{\sigma_{ij}, 1\leq i \leq n, 1\leq j\leq m\right\}$ be uniformly chosen from $\left\{-1, +1\right\}$. 
For mathematical rigor, we introduce the definition $\mathbb{E}\sup_{f\in \mathcal{F}}V(f)=\sup \left\{ \mathbb{E}\sup_{f \in \mathcal{F}^\ast} V(f): \mathcal{F}^\ast \subset \mathcal{F}, \mathcal{F}^\ast \ \mbox{is finite}\right\}$ for a random process $V$
when the index set $\mathcal{F}$ is uncountable.

\begin{definition}[Rademacher fixed point]\label{ass:fixpoint}
Let $r^{*}_{nm}$ be a positive number and $\phi_{nm}(r)$ be  a sub-root function for $r\geq r^{*}_{nm}$.  Assume that $\phi_{nm}( r^{*}_{nm})\leq r^{*}_{nm}$ and 
$$
\phi_{nm}(r) \geq \mathbb{E} \left[\sup_{f\in \mathcal{F}_{nm}(r)}\left| \frac{1}{nm}\sum_{i=1}^{n}\sum_{j=1}^{m} \sigma_{ij} \left( f(\x_{ij}) - \fo(\x_{ij})\right) \right| \right]
$$
for  $r \geq r^{*}_{nm}$. Then $r^{*}_{nm}$ is called the Rademacher fixed point of $\mathcal{F}_{nm}$. 
\end{definition}

The positive number $r^{*}_{nm}$ is usually referred to as the fixed point of the continuity modulus of the Rademacher average that characterizes the neighborhood complexity of $\fo$ in  $\mathcal{F}_{nm}$. 
In the empirical process with independent random variables  \citep{bartlett2005local,koltchinskii2006local}, the fixed point is found to be of the same order as the convergence rate of the estimators that minimize an empirical loss criterion over a given function class. 
Our first main result in the next theorem establishes a non-asymptotic bound on the prediction error in the repeated measures model using $r^{*}_{nm}$.

\begin{theorem}\label{thm:thm1} 
Consider the repeated measures model \eqref{equ:model} and the estimator $\hat{f}_{nm}$ obtained by \eqref{equ:estimator}. 
Under Assumptions \ref{ass:fobound} and \ref{ass:subexp}, there is a  
universal constant $\xc\label{c001}$ such that 
\begin{equation}
\begin{aligned}
\pef \leq &
\jc{c001}\left(\inf_{f\in\mathcal{F}_{nm}} \int_{\Omega}\left(f-\fo\right)^2d\mathcal{P}_{\x} 
+\frac{\jdb{b1}^2+\jdb{b2}^2}{n}\right.
\\ & \quad \quad + \left.(r^{*}_{nm} + \frac{1}{nm})(\jdb{b1}^2+\jdb{b2}^2(\log n)^2+\jdb{b3}^2(\log nm)^2)\right).\\
\end{aligned} \label{equ:thm1}
\end{equation}
\end{theorem}

The result above can also take optimization error into account using standard techniques, and Remark S.5.1 in the Supplementary Material discusses the lower and upper bounds when the estimator is not the global minimizer. 
By Theorem \ref{thm:thm1}, the prediction error of an empirical risk minimizer $\hat{f}_{nm}$ can be bounded by the approximation error $\inf_{f\in\mathcal{F}_{nm}}  \int_{\Omega}\left(f-\fo\right)^2d\mathcal{P}_{\x}$ and the estimation error $(\jdb{b1}^2+\jdb{b2}^2) n^{-1}+ (r^{*}_{nm} + (nm)^{-1})(\jdb{b1}^2+\jdb{b2}^2(\log n)^2+\jdb{b3}^2(\log nm)^2)$.  
Compared with the cross-sectional model
\begin{equation} \label{equ:modelcs}
Y_{i}=\fo(\x_{i})+\epsilon_{i}, \hspace{0.2in} 1\leq i \leq n, 
\end{equation}
which can be viewed as a special case of \eqref{equ:model} that measures a single observation for each subject.  Accordingly, letting $m=1$ and $\jdb{b2}=0$ yields a parallel result for model (\ref{equ:modelcs}). We remark that the term in (\ref{equ:thm1}) involving $n^{-1}$  plays an interesting role in the repeated measurements model, coupling with   
the terms $\inf_{f\in\mathcal{F}_{nm}}  \int_{\Omega}\left(f-\fo\right)^2d\mathcal{P}_{\x}+(r^{*}_{nm} + (nm)^{-1})(\log nm)^2$ depending on the total number of observations $n m$. 
To minimize the right-hand side (r.h.s.) of (\ref{equ:thm1}), we choose the function class $\mathcal{F}_{nm}$ that usually expands with the sample size to balance the trade-off between two errors. Then the upper bound on the r.h.s. can be shown of the order $n^{-1}+(nm)^{-\alpha}$ up to logarithmic factors for specific problems. 
This reveals some fundamentally different behavior of nonparametric regression estimators in the repeated measurements model. No matter how large $m$ is, even when the $n$ functions are fully observed, one cannot make the prediction error faster than the order of $n^{-1}$ that is customarily called the parametric rate in terms of sample size. 

The logarithmic terms in (\ref{equ:thm1}) are caused by technical issues due to the sub-exponential tails in Assumption \ref{ass:subexp}, and will vanish if ones assume bounded variables as in \citet{gyorfi2002distribution}. 
A similar result could be also given with sub-Gaussian assumptions and is omitted. 
It is worth noting that the bound $\jdb{b1}$ for the target function $\fo$ and $f\in \mathcal{F}_{nm}$ can be relaxed as a sequence that slowly grows with $nm$, especially when one has no prior information about such functions.
For example, we can set $\jdb{b1}=\log n$ for the chosen function class $\mathcal{F}_{nm}$ as in (S.3) in the Supplementary Material, and the convergence rate in the result remains nearly unchanged except for a logarithmic factor. In the sequel, we do not distinguish whether $\jdb{b1}$ is finite or increases with $mn$ when no confusion arises. For convenience, we denote $\iota_{nm}=\jdb{b1}^2+\jdb{b2}^2(\log n)^2+\jdb{b3}^2(\log nm)^2$.

Next, we extend the result in Theorem \ref{thm:thm1} for several situations of interest. The reason is that the bound (\ref{equ:thm1}) utilizes the Rademacher fixed points for $\mathcal{F}_{nm}$ which is an excellent distribution-dependent tool for measuring the complexity of the candidate function class. However, this can only be computed when the distribution is known, which is often not the case in practice. 
One proposal is to estimate $r^{*}_{nm}$ by an empirical version like Theorem 4.2 in \citet{bartlett2005local} that is data-dependent. Alternatively, a feasible resolution is to derive an upper bound using additional knowledge about the complexity of the class.

\begin{corollary}\label{cor:cor1} Under the same model and assumptions as Theorem \ref{thm:thm1}, we have 

(i) random/uniform covering numbers:
$$
\begin{aligned}
  \pef 
\leq 
\xc\label{c:cover} &\left(\inf_{f\in\mathcal{F}_{nm}} \int_{\Omega}\left(f-\fo\right)^2d\mathcal{P}_{\x} 
+\frac{\jdb{b1}^2+\jdb{b2}^2}{n}\right.\\& + \left.
\frac{\iota_{nm}}{nm}
\left(  \mathbb{E}\left[ \log \mathcal{N}\left(\frac{\jdb{b1}}{nm}, \left.\mathcal{F}_{nm}\right|_{\mathcal{X}}, L^2 \right) \right] + \log nm  \right) 
\right),\\
\end{aligned}
$$
and the bound in uniform covering number is of the same form, but with the expectation of $ \log \mathcal{N}\left((nm)^{-1}\jdb{b1}, \left.\mathcal{F}_{nm}\right|_{\mathcal{X}}, L^2 \right) $ replaced by its supremum $\log \mathcal{N}_{nm}\left((nm)^{-1}\jdb{b1}, \mathcal{F}_{nm}, L^2 \right)$.

(ii) VC dimension:
$$
    \begin{aligned}
  \pef &
\leq 
\xc\label{c:vc}\left(\inf_{f\in\mathcal{F}_{nm}} \int_{\Omega}\left(f-\fo\right)^2d\mathcal{P}_{\x} 
+\frac{\jdb{b1}^2+\jdb{b2}^2}{n}
+\frac{  \vc(\mathcal{F}_{nm}) \iota_{nm} \log nm}{nm} \right).\\
\end{aligned}
$$
where $\jc{c:cover}$ and $\jc{c:vc}$ are universal constants. 
\end{corollary}

Lastly, we present the case that $\mathcal{F}_{nm}$ is a set of truncated fully connected ReLU feedforward neural networks $\nn (d, L, W, \jdb{b1})$, which demonstrates how the width and depth of the network determine the estimation error explicitly.

\begin{corollary}\label{cor:networksize}
Under the same model and assumptions as Theorem \ref{thm:thm1} and set the function class $\mathcal{F}_{nm}=\nn (d, L, W, \jdb{b1})$ with $\log LW \geq 1$, then
$$
\begin{aligned}
  \pef &
\leq 
\xc\label{c:nn}\left(\inf_{f\in\mathcal{F}_{nm}} \int_{\Omega}\left(f-\fo\right)^2d\mathcal{P}_{\x} 
+\frac{\jdb{b1}^2+\jdb{b2}^2}{n} + 
\frac{L^2W^2 \log LW \iota_{nm} \log nm}{nm}  
\right), \\
\end{aligned}
$$
where $\jc{c:nn}$ is a positive constant. 
\end{corollary}

Since this result is the first in the repeated measurement setting, we may only compare it with existing results for the cross-sectional setting (\ref{equ:modelcs}) by letting $m=1$.
First, our error bounds work for sub-exponential noise, while  
\citet{gyorfi2002distribution} assumes the response to be bounded and \citet{schmidtaos} requires the noises to be not heavier than Gaussian. 
Second, the estimation error in our results is bounded in various ways, including Rademacher fixed point, random/uniform covering numbers and VC dimension. This allows easy extensions for future works, even if $\mathcal{F}_{nm}$ is not a feedforward neural network class.\label{ju1}
In particular, these results are usually sharper than the functional covering number which leads to the further boundedness requirement of network weights and bias.

\section{H\"older Space and Phase Transition}\label{sec:holderph} 
In this section, we apply these oracle inequalities to the situation that the target function belongs to the H\"older space that has general smoothness assumptions and plays an important role in nonparametric regression problems. 
We first establish a theoretical bound of the DNN estimator considered above, then demonstrate a matching lower bound and analyze its optimal convergence rate, while the phase transition owing to repeated measurements is revealed. 
Assume $\Omega = [0,1]^d$ in this section for brevity.
Recall the definition of H\"older space. 

\begin{definition}[H\"older space  $\mathcal{C}^{s}\left(\Sigma\right)$] \label{def:holder}
Let $s$ be a positive number and $\lfloor s \rfloor$ be the largest integer strictly smaller than $s$. The $s$-H\"older norm of a function $f:\Sigma\rightarrow\mathbb{R}$ is defined as 
$$
\|f\|_{\mathcal{C}^{s}}=
\sum_{\boldsymbol{\alpha}:|\boldsymbol{\alpha}|<s}\left\|\partial^{\boldsymbol{\alpha}} f\right\|_{\infty}+\sum_{\boldsymbol{\alpha}:|\boldsymbol{\alpha}|=\lfloor s \rfloor} \sup _{\substack{\xx, \xx' \in \Sigma \\ \xx \neq \xx'}} \frac{\left|\partial^{\boldsymbol{\alpha}} f(\xx)-\partial^{\boldsymbol{\alpha}} f(\xx')\right|}{\|\xx-\xx'\|^{s-\lfloor s \rfloor}},
$$
where  $\partial^{\alpha}:=\partial^{\alpha_{1}}_{1}\partial^{\alpha_{2}}_{2} \ldots \partial^{\alpha_{d}}_{d}$  and $|\boldsymbol{\alpha}|:=\alpha_{1}+\alpha_{2}+\dots+\alpha_{d}$ with $\boldsymbol{\alpha}=\left(\alpha_{1}, \ldots, \alpha_{d}\right) \in \mathbb{N}^{d}$.
H\"older space $\mathcal{C}^{s}\left(\Sigma\right)$ consists of all functions for which the $s$-H\"older norm is finite. Typical examples include Lipschitz functions and continuously differentiable functions on a compact set. 
\end{definition}

\begin{assumption}[H\"older smoothness] \label{ass:holder}
The target function $\fo$ belongs to the H\"older space $\mathcal{C}^{s}\left([0,1]^d\right)$ 
for a given positive number $s$, and $\|\fo \|_{\mathcal{C}^{s}}\leq \xdb\label{holdernorm}$  for some constant $\jdb{holdernorm}>0$.
\end{assumption}

The non-asymptotic error bound of the ReLU feedforward neural network estimator is determined by a trade-off between the estimation error and the approximation error in Corollary \ref{cor:networksize}, associated with the depth and width. 
Combining Corollary \ref{cor:networksize} and the approximation error in \citet{huangjian}, we arrive at the following result.

\begin{theorem}\label{thm:holdernn}
Consider the model \eqref{equ:model} and the estimator $\hat{f}_{nm}$ obtained by \eqref{equ:estimator}. Suppose that Assumptions \ref{ass:fobound}, 
\ref{ass:subexp} and \ref{ass:holder} hold, and set
$\mathcal{F}_{nm} = \nn\left(d, \xc\label{c500} L \log L, \xc\label{c501} W \log W, \jdb{b1}\right)
$ with $\jc{c500}, \jc{c501}$ not depending on $L, W$.    
If we specify the neural network in a flexible way satisfying 
$
LW=\lfloor \xc\label{c55} (nm)^{ d/(4s+2d) } (\log nm)^{  - 4d/(2s+d) } \rfloor
$
for some constant $\jc{c55}$ not depending on $n$ and $m$, then  
$$
\pef \leq \xc\label{c58} ( n^{-1}+ (nm)^{- 2s/(2s+d)}(\log nm)^{16s/(2s+d)} ), 
$$
where $\jc{c58}$ is a constant free of $n$ and $m$. 
\end{theorem}

Let $\mathcal{P}_{Z,\epsilon}$ be the probability measure for the random function $Z$ and noise $\epsilon$, and $\mathcal{Q}_{Z,\epsilon}$ be a collection of $\mathcal{P}_{Z,\epsilon}$ satisfying some certain conditions. 
The next theorem provides a lower bound on the minimax convergence rate for repeated measurement regression estimators. 
\begin{theorem}\label{thm:low1}
Consider the repeated measures model \eqref{equ:model} and assume that $\mathcal{P}_{\x}$ has a uniform distribution on $\Omega$.
Under the Assumptions \ref{ass:subexp} and \ref{ass:holder}, we have 
$$
    \inf_{\tilde{f}}\sup_{ \mathcal{P}_{Z,\epsilon} \in \mathcal{Q}_{Z,\epsilon} 
    }  \mathcal{E}_{\fo}(\tilde{f})
\geq \xc\label{c60} (n^{-1} + (nm)^{-2s/(2s+d)}).
$$
where the infimum is taken over all possible estimators $\tilde{f}$ based on the available data set $\left\{ \left(Y_{ij},\x_{ij} \right): 1\leq i \leq n, 1\leq j \leq m \right\}$, the supremum is taken over all distribution $\mathcal{P}_{Z,\epsilon}\in \mathcal{Q}_{Z,\epsilon}$ satisfying Assumptions \ref{ass:subexp} and \ref{ass:holder}, and $\jc{c60}$ is a positive number not depending on $n$ and $m$. 
\end{theorem}

Combining the two theorems yields that the optimal convergence rate is of the order $O(n^{-1} + (nm)^{-2s/(2s+d)})$ (up to logarithmic factors) for H\"older smoothness index $s\in \mathbb{R}_{+}$. 
It matches the existing minimax rate $n^{-2s/(2s+d)}$ in cross-sectional case \citep{stone1980,schmidtaos,kohler2021rate} by setting $m=1$, 
and $n^{-1} + (nm)^{-2s/(2s+1)}$ in repeated measures regression for RKHS with $O(k^{-2s})$ eigendecay \citep{cai2011optimal} and $n^{-1} + (nm)^{-4/5}$ for local linear regression by assuming $s=2,d=1$ \citep{zhang2016sparse}. 
Intriguingly, we observe a phase transition phenomenon in the minimax convergence rate. To see this, when the sampling frequency $m \gtrsim n^{d / 2s}$, the rate reaches $n^{-1}$ which is usually called the parametric rate of convergence in terms of the sample size $n$. After that, no matter how many more repeated measurements are collected, the order $n^{-1}$ cannot be improved, even when the $n$ functions are fully observed which is considered the ideal situation in functional data analysis. On the other hand, 
When $m \lesssim n^{d / 2s}$ that is viewed as the sparse sampling scheme, the convergence is slowed to an of order $(nm)^{-2s/(2s+d)}$ that depends on the total number of observations $nm$ and in accordance to the so-called nonparametric rate.
For $d=1$, the phenomenon was studied in integer Sobolev smoothness \citep{cai2011optimal} and continuously differentiable functions \citep{zhang2016sparse}, as well as in Riemannian manifold space  \citep{shao2022intrinsic} and online estimation setting \citep{yang2022online}. 

This phase transition reveals the impact of dimensionality on the convergence rate in two folds. 
First, when the sample size is $n$, the optimal rate $n^{-1}$ is not improvable. To achieve this rate, the sampling frequency should be of order $m \asymp n^{d / 2s}$. 
A higher dimension $d$ requires a much larger number of measurements that has to grow exponentially in terms of the dimension $d$. 
Second, when the subjects are sparsely sampled, the rate $(nm)^{-2s/(2s+d)}$  would be fairly slow when the dimension is relatively large compared to smoothness, hence the curse of dimensionality behaves similarly to the cross-sectional setting with $nm$ independent observations.  
We close the section by mentioning that the same technique can be applied to other commonly used $\mathcal{F}_{nm}$ to establish the associated optimal convergence analysis, which is also illustrated for regression spline estimators in Section S.2 in the Supplementary Material. 

\section{Circumventing Curse of Dimensionality}\label{sec:cir}
In modern statistical modeling and machine learning, the prevalence of multi-dimensional data has posed the formidable challenge of the curse of dimensionality, particularly salient for nonparametric regression. 
Our previous findings indicate that while the convergence of DNN-based estimators in repeated measurements differs from that in cross-sectional context, the curse of dimensionality persists in smoothness class like H\"older space, which requires a much larger sampling frequency to achieve the optimal parametric rate. 
In this work, we examine several low intrinsic-dimensional function spaces, encompassing the hierarchical composition model, manifold input and anisotropic H\"older smoothness. For brevity, this section only shows the results for the first two cases, while additional results are given in Section S.1.2 in Supplementary Material. Demonstrably, DNN estimators exhibit adaptability to these structures without prior knowledge and specialized design, thereby mitigating the curse of dimensionality to a certain extent.

\subsection{Hierarchical Composition Structure}\label{sec:hier}

The classical approach in statistics to circumvent the curse of dimensionality is to make semiparametric model assumptions, such as single index models, generalized additive models and partially linear models, so that faster convergence can be achieved.
In fact, these model assumptions can be subsumed within a broader framework of the hierarchical composition structure \citep{bauer2019deep,schmidtaos,kohler2021rate} that has optimal cross-sectional regression rate $n^{-2\gamma/(2\gamma+1)}$ with some index $\gamma$.
Now, we introduce the definition of the Hierarchical composition model $\mathcal{H}^{l, \mathcal{G}}(\Sigma)$ in this paper.

\begin{definition}[Hierarchical composition model $\mathcal{H}^{l, \mathcal{G}}(\Sigma)$] Suppose $\Sigma \subseteq \mathbb{R}^d$, $l\in \mathbb{N}$ and $\mathcal{G}$ is a $l$ height tree with nodes belonging to $(0, \infty)\times \mathbb{N}$, then the hierarchical composition model $\mathcal{H}^{l, \mathcal{G}}(\Sigma)$ is recursively defined as follows. For $l=1$ and  $\mathcal{G}=((s,K))$, 
$$
\begin{aligned}
\mathcal{H}^{l, \mathcal{G} }(\Sigma)=&\Big\{f: \Sigma \rightarrow \mathbb{R}:  f(\xx)=g(x_{\pi(1)}, x_{\pi(2)}, ...,x_{\pi(K)}) \text{, where } \\
&\ \ \quad  g: \mathbb{R}^{K} \rightarrow \mathbb{R} \text{ satisfies } \|g\|_{\mathcal{C}^{s}}\leq M 
\text{ and } \pi:[K] \rightarrow[d]
\Big\}.
\end{aligned}
$$
For $l\geq 2$, given $\mathcal{G}=((s,K);\mathcal{G}_{1},...,\mathcal{G}_{K} )$ with $(s,K)\in (0,\infty)\times \mathbb{N}$, 
$$
\begin{aligned}
\mathcal{H}^{l, \mathcal{G} }(\Sigma)=\Big\{&f: \Sigma \rightarrow \mathbb{R}:  f(\xx)=g\left(h_{1}(\xx),h_{2}(\xx), ..., h_{K}(\xx)\right) \text{, where } \\
& g: \mathbb{R}^{K} \rightarrow \mathbb{R} \text{ satisfies } \|g\|_{\mathcal{C}^{s}}\leq M  
\text{ and } h_{k} \in \mathcal{H}^{l_{k}, \mathcal{G}_{k}}(\Sigma) \text{ with } l_{k}\leq l-1 \text{ for } 1\leq k\leq K
\Big\}.
\end{aligned}
$$
Here, $M$ is a positive constant and $x_{\pi(j)}$ is the $\pi(j)$-th component of $\xx$.
\end{definition}

Note that our definition of the hierarchical composition model contains more information than that in \citet{kohler2021rate}. 
Due to the difficulty in deriving the convergence upper bound, Theorem 1 in \citet{kohler2021rate} demands that the smoothness of functions in each layer should be no less than 1, which limits its applicability. 
Therefore, in the above definition, we introduced $\mathcal{G}$ with a tree data structure to store the smoothness and dimension of functions at each level to deal with the transport of sub-smoothness in the composition.  
Further, we define the intrinsic smoothness-dimension ratio $\gamma$ by recursion. 
Initially, let $\gamma=s/K$ for $l=1$ and $\mathcal{G}=((s,K))$. 
For $l\geq 2$, assume the root node of $\mathcal{G}$ is $(s,K)$ and $\mathcal{G}_{1}, \mathcal{G}_{2}, ...,\mathcal{G}_{K}$ are sub-trees whose roots are children of $(s,K)$. 
Then we define 
$$\gamma= \min\left(s/K, \min_{1\leq k\leq K}\gamma_{k}\cdot \min(1,s)\right),$$ 
where $\gamma_{k}$ are the intrinsic smoothness-dimension ratios of $\mathcal{H}^{l-1, \mathcal{G}_{k} }(\Sigma)$. 
Alternatively, the intrinsic smoothness-dimension ratio could be calculated directly. 
For each node $G=(s_{G},K_{G})$ of $\mathcal{G}$, define the effective smoothness index as 
$
s_{G}^{*}=s \cdot\prod_{(s',K')\in\mathcal{A}(G)}\min(1,s'),
$
where $\mathcal{A}(G)$ denotes the ancestor set of $G$ in the tree $\mathcal{G}$. One can verify that 
$
\gamma=\min_{G\in \operatorname{node}(\mathcal{G})} s_{G}^{*}/K_{G}, 
$
where the minimum value is taken over all nodes $G$ of tree $\mathcal{G}$. 

We derive a new result that characterizes the approximation ability of deep ReLU feedforward neural network to the hierarchical composition function, which sets the stage for studying the error bound and optimal rate of convergence.  

\begin{theorem}\label{lem:hie}
Suppose that $f \in \mathcal{H}^{l, \mathcal{G}}(\Omega)$, $\|f\|_{\infty}\leq \jdb{b1}$ and $\gamma$ is defined as above. 
Then, for any $L,W\geq 3$, there exists some neural network 
$f_{\nn}\in \nn\left(d, \xc\label{c506} L \log L, \xc\label{c507} W \log W, \jdb{b1}\right)$ satisfying
$$
|f_{\nn}(x)-f(x)|\leq \xc\label{c508}\left(LW\right)^{-2 \gamma},
$$
for every $x\in \Omega$, where
$\jc{c506}$ and $\jc{c507}$ are constants not depending on $L$, $W$ and $f$, and 
$\jc{c508}$ is a constant free of $L$ and $W$.
\end{theorem}

Now we assume the target function $\fo$ has such a hierarchical structure, and arrive at the following non-asymptotic error bound for repeated measurements model by coupling with Corollary \ref{cor:networksize}, which shows how this alleviates 
the curse of dimensionality.
\begin{assumption}[Hierarchical structure] \label{ass:f0Hiera}
The target function $\fo(\xx)$ belongs to the hierarchical composition model $\mathcal{H}^{l, \mathcal{G}}(\Omega)$ for given $l$, $\mathcal{G}$ and $M$. 
\end{assumption}

\begin{theorem}\label{thm:hieup}
Consider model \eqref{equ:model} and the estimator $\hat{f}_{nm}$ obtained by \eqref{equ:estimator}. Suppose that Assumptions \ref{ass:fobound}, 
\ref{ass:subexp} and \ref{ass:f0Hiera} hold, and set $\mathcal{F}_{nm} = \nn\left(d, \jc{c500} L \log L, \jc{c501} W \log W, \jdb{b1}\right).
$
If we specify the neural network in a flexible way satisfying 
$
LW=\lfloor \xc\label{c62} (nm)^{1/(4 \gamma +2)} (\log nm)^{ -4/(2\gamma +1) } \rfloor
$
for some constant $\jc{c62}$ not depending on $n$ and $m$, then 
$$
\pef \leq \xc\label{c61} ( n^{-1}+ (nm)^{-2\gamma/(2\gamma+1)   }(\log  nm)^{ 16\gamma/(2\gamma+1)} ), 
$$
where $\jc{c61}$ is a constant free of $n$ and $m$. 
\end{theorem}

The following lower bound is constructed for the repeated measurement setting and the hierarchical composition model. 
\begin{theorem}\label{thm:hielow}
Consider the repeated measures model \eqref{equ:model} and assume that $\mathcal{P}_{\x}$ has a uniform distribution on $\Omega$. 
Assume that there exists a node $G=(s_{G},K_{G})$ such that $s_{G}^{*}/K_{G}=\gamma$ and $K_{G}\leq d$. 
Under the Assumption \ref{ass:subexp} and \ref{ass:f0Hiera}, we have
$$
    \inf_{\tilde{f}}\sup_{\mathcal{P}_{Z,\epsilon} \in \mathcal{Q}_{Z,\epsilon}}  \mathcal{E}_{\fo}(\tilde{f})
\geq \xc\label{c63} (n^{-1} + (nm)^{-2\gamma/(2\gamma+1)}),
$$
where the infimum is taken over all possible estimators $\tilde{f}$ based on the available data set $\left\{ \left(Y_{ij},\x_{ij} \right): 1\leq i \leq n, 1\leq j \leq m \right\}$, 
the supremum is taken over all distributions $\mathcal{P}_{Z,\epsilon} \in \mathcal{Q}_{Z,\epsilon}$  satisfying Assumptions \ref{ass:subexp} and \ref{ass:f0Hiera}, 
and $\jc{c63}$ is a positive number not depending on $n$ and $m$. 
\end{theorem}

The previous two theorems demonstrate that the minimax optimal rate in the hierarchical composition model becomes $O(n^{-1} + (nm)^{-2\gamma/(2\gamma+1)})$ up to some logarithmic factors. 
That means the intrinsic smoothness-dimension ratio $\gamma$, rather than the extrinsic smoothness-dimension ratio $s/d$, determines the convergence rates under the assumption of the representational function set $\mathcal{H}^{l, \mathcal{G}}(\Omega)$, while one typically has $\gamma<s/d$.
Accordingly, the phase transition phenomenon inherent to repeated measurements takes place at the order $m \asymp O(n^{1/2\gamma})$ instead of $O(n^{d/2s})$.  
Hence, the curse of dimensionality is lessened, and exponentially less sample frequency is required to achieve the parametric rate.

In traditional works, hierarchical structure often appears as model assumptions in one of its specific forms, such as assuming that the true function satisfies a single-index model, additive model, varying coefficient model, or their combinations. Given a particular structure, tailored estimation methods are designed. A closely related work is \citet{horowitz2007} which used the penalized least squares method to estimate the known composition of the additive model.
For the adaptivity of traditional tools to unknown hierarchical structures, \citet{schmidtaos} demonstrated that wavelet series do not adapt to unknown hierarchical structures. \citet{suzuki2021deep} extended this conclusion on non-adaptivity to more general linear estimators, including local polynomials, RKHS, regression splines, and wavelets.
For completeness, we also show that the linear estimators cannot be adaptive to hierarchical models without specific knowledge and designs in the repeated measurements model, even if the true function is of a single index form; see Section S.1.1 in Supplementary Material for more details.

\subsection{Low-dimensional Support}\label{sec:lowipz}
Considering predictors in a low-dimensional structure is sensible in many applications, such as image processing, natural language processing and bio-medical studies. 
There is rich literature to study regression on an unknown manifold $\mathcal{M}$, especially for local polynomial methods \citep{bickel2007local,aswani2011regression,cheng2013local}. 
\citet{bickel2007local} revealed that local polynomial regression can adapt to the unknown low dimensional structure and achieve the optimal convergence rate $O(n^{-2s/(2s+d_{\mathcal{M}})})$, \citet{aswani2011regression} and \citet{cheng2013local} learned the manifold structure first and consider further estimation.  
Recent works in the cross-sectional regression have shown that DNNs are also adaptive to intrinsic structures  \citep[e.g.][]{chen2019efficient,schmidt2019deep,nakada2020adaptive,huangjian}, classification \citep[e.g.][]{liu2021besov,zhang2024nonparametric}, and more advanced scenarios \citep[e.g.][]{dahal2022deep,xu2023sample}. 
The previous regression works are all cross-sectional, and now we consider the repeated measures model when $\Omega = \mathcal{M}$.

\begin{assumption}[Manifold support] \label{ass:manifold}
The predictor $\x$ is supported on a compact Riemannian manifold $\mathcal{M}$ that is $d_{\mathcal{M}}$-dimensional and isometrically embedded in $[0,1]^d$ with conditional number $1/\tau$ and a finite area. 
The target function $\fo(\xx)$ belongs to the H\"older space $\mathcal{C}^{s}\left( \mathcal{M} \right)$ for a given positive number $s$.
\end{assumption}
In this assumption, the condition number is defined to be $1/\tau$, where $\tau$ is the largest number such that, for any $r<\tau$, the open normal bundle about $\mathcal{M}$ of radius $r$ is embedded in $\mathbb{R}^d$. Intuitively, a smaller condition number means the manifold is flatter, and it is easier to approximate by DNNs. 
Definition \ref{def:holder} cannot be directly generalized to manifolds, a detailed definition of $\mathcal{C}^{s}\left( \mathcal{M} \right)$ is provided in Section S.1.3 in the Supplementary Material.

\begin{theorem}\label{thm:maniup}
Consider model \eqref{equ:model} and the estimator $\hat{f}_{nm}$ obtained by \eqref{equ:estimator}.
Suppose that Assumptions \ref{ass:fobound}, 
\ref{ass:subexp}, \ref{ass:holder} and \ref{ass:manifold} hold, set $\mathcal{F}_{nm} = \nn\left(d, \xc\label{c503p} L \log L, \xc\label{c504p} W \log W, \jdb{b1}\right).
$
If we specify the neural network in a flexible way satisfying 
$
LW=\lfloor \xc\label{c76} (nm)^{d_{\mathcal{M}} / (4s+2d_{\mathcal{M}})} (\log nm)^{-4d_{\mathcal{M}}/(2s+d_{\mathcal{M}}) } \rfloor
$ 
for some constant $\jc{c76}$ not depending on $n$ and $m$, then 
$$
\pef \leq \xc\label{c77} ( n^{-1}+ (nm)^{-2s/(2s+d_{\mathcal{M}}) }(\log nm)^{16s/(2s+d_{\mathcal{M}})}), 
$$
where $\jc{c77}$ is a constant independent of $n$ and $m$. 
\end{theorem}

Combining with the minimax lower bound in Theorem \ref{thm:low1}, it follows that the ReLU DNN estimator attains the optimal minimax rate $O(n^{-1}+(nm)^{-2s/(2s+d_{\mathcal{M}})})$ up to logarithmic factors under the manifold support set assumption. 
The convergence rate is usually faster than $O(n^{-1}+(nm)^{-2s/(2s+d)})$ since the intrinsic dimension $d_{\mathcal{M}}\le d$. 
Accordingly, the phase transition takes place at $O(n^{d_{\mathcal{M}}/2s})$ rather than $O(n^{d/2s})$, i.e., which hints that substantially less repeated measurements are required to obtain the parametric rate.

We close this section by mentioning that the results developed for the repeated measurements model are readily applicable to other low-dimensional structures like the approximate manifold and Minkowski dimension set. 
Some error bounds use recent results\citep{huangjian} have been shown in Section S.1.3 in the Supplementary Material.

\section{ Simulation Studies} \label{sec:sim}
{We conduct simulation studies to demonstrate the adaptivity of the DNN estimators and the phase transition phenomenon. Under various unknown structures, the comparisons with RKHS method in \citet{cai2011optimal} and local linear method in \citet{zhang2016sparse} are also provided since they are two available common methods of allowing arbitrary sampling frequency $m$. The numerical behaviors of the regression spline method in Section S.2 of the Supplementary Material for cases with $d\leq 5$ that our computational capacity allows are also shown. For more simulated examples on hierarchical structure with common closed forms and on anisotropic smoothness, see Section S.3 in the Supplementary Material. 

The Adam optimization algorithm is used with a maximum of 200 epochs and early stopping with a patience of $10$ for regularization.
In all simulations, the $d$-dimensional covariate $X$ is generated from the uniform distribution on $\Omega \subseteq [0,1]^d$.
Recalling the model setting $Y_{ij}=\fo(\x_{ij})+U_{i}(\x_{ij})+\epsilon_{ij}$, we take the individual stochastic part  $U_{i}(\xx)=\sum_{k=1}^{\infty}\sum_{l=1}^{d}\frac{\sqrt{3}\xi_{ikl}}{\sqrt{d}k}\cos(k\pi\xx_{l})$ with independent scores $\xi_{ikl} \sim N(0,0.1^2)$, and measurement errors $\epsilon_{ij} \sim N(0,0.01^2)$. Let $a=0.3$ and $b=1.6$, we consider the following cases:
\begin{align*}
& \text { Case } 1 : \fo(x) =
\sum_{ k_1,...,k_5\geq 1} a^{ \max\{k_1,k_2,k_3,k_4,k_5\}} \prod_{l=1}^{5}\cos(2\pi l^{-1} b^{k_l} \xx_l), \Omega = [0,1]^5,  \\
& \text { Case } 2 : \fo(x) = \sum_{l=1}^5 \sum_{ k\geq 1} a^{k} \cos(2\pi l^{-1} b^{k_l} \xx_l), \Omega = [0,1]^5, \\
& \text { Case } 3 : \fo(x) =
\sum_{ k_1,...,k_5\geq 1} a^{ \max\{k_1,k_2,k_3,k_4,k_5\}} \prod_{l=1}^{5}\cos(2\pi l^{-1} b^{k_l} \xx_l),\\
&\quad\quad\quad\quad \Omega = \{(\cos(2\pi t_1), \sin(2\pi t_1),\sin(2\pi t_2),t_1,t_2):t_1,t_2\in[0,1]\},\\
& \text { Case } 4 : \fo(x) = \sin\left( \frac{12 \pi }{d(d+1)} \sum_{l=1}^{d} l x_l \right) , \Omega = [0,1]^d, d=10,20,30
\end{align*}
Cases 1-3 take the form of series expansions, ensuring their appropriate smoothness inversely related with the values of $a$ and $b$.
Case 2 represents an additive function, thus serving as an example of the hierarchical model. 
In Case 3, the support $\Omega$ is an intrinsically 2-dimensional manifold. 
Case 4 contains large-dimensional examples with low intrinsic dimensions due to the hierarchical structure.

For each case, we conduct 50 Monte Carlo runs and the datasets of smaller $n$ and $m$ are subsets of larger ones. 
A separate i.i.d. set, comprising $25\%$ of the training sample size, is used for hyperparameter selection of each method. 
For DNN estimators, the networks are setting to be $(L,W) = (1,50),(2,100),(3,200),(4,400),(5,600),(5,800)$ with equal width hidden layers.
Each example trains all networks with different depths and widths, then we select the model with the minimal validation error. 
The estimate $\hat{f}$ is then evaluated using empirical mean squared error, computed over an independent set of $10^4$ independently and uniformly distributed points. 
Results of DNN estimators are shown in Figure \ref{fig:simualtion}, and the comparisons are in Table \ref{tab:bijiao}.

\begin{figure}[!ht]
\spacingset{1.36}
\scriptsize
    \centering
    \includegraphics[width=\textwidth]{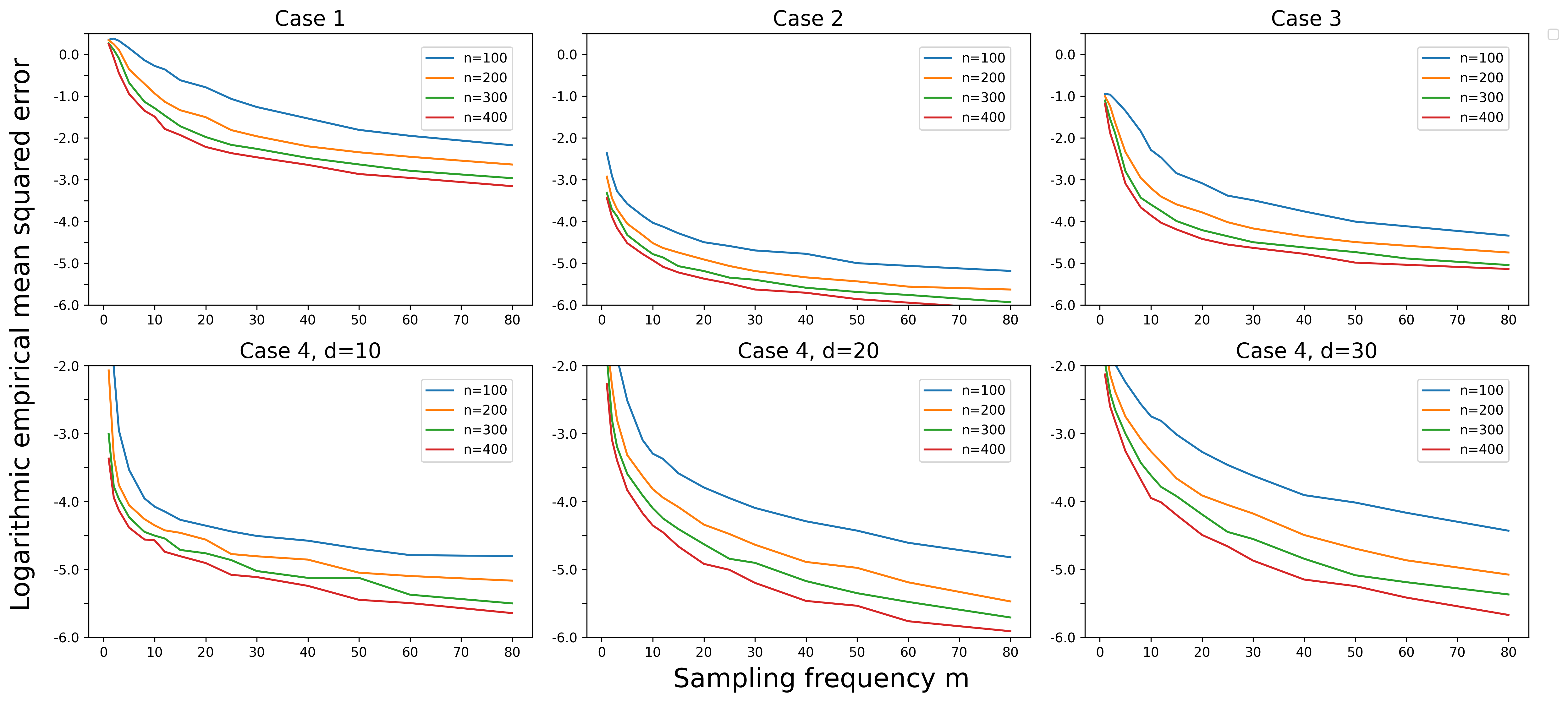}
	\caption{Finite sample mean squared error of the DNN estimators for $n=100,200,300,400$ and sampling frequencies $m=1,2,3,5,8,10,12,15,20,25,30,40,50,60,80$. Each subplot portrays sampling frequencies ($m$) on the horizontal axis and logarithmic empirical mean squared error on the vertical axis. }\label{fig:simualtion}
\end{figure}
\vspace{-0.5cm}

\begin{table}[ht]
\spacingset{1.36}
\centering
\scriptsize
    \caption{Finite sample mean square error of the DNN estimator, regression spline estimator (for cases with $d\leq 5$), RKHS and local linear regression estimators. }
    \setlength{\tabcolsep}{0.03in}
	\label{tab:bijiao}
 \begin{tabular}{ccccccccccccccccccc}
		\hline
		  \multicolumn{2}{c}{}& \multicolumn{5}{c}{Case 1}&&\multicolumn{5}{c}{Case 2}&&\multicolumn{5}{c}{Case 3}\\
    \cline{3-7}\cline{9-13} \cline{15-19}
    \multicolumn{2}{c}{Sample size $n$} & \multicolumn{2}{c}{$100$}& & \multicolumn{2}{c}{$300$} & &\multicolumn{2}{c}{$100$}& & \multicolumn{2}{c}{$300$} & &\multicolumn{2}{c}{$100$}& & \multicolumn{2}{c}{$300$} \\
    \cline{3-4} \cline{6-7}\cline{9-10} \cline{12-13} \cline{15-16}\cline{18-19}
    \multicolumn{2}{c}{Sampling frequency  $m$} 
    & $10$ &$30$&& $10$ &$30$&& $10$ &$30$&& $10$ &$30$&& $10$ &$30$&& $10$ &$30$\\
    \hline
    &DNN  & 0.7598 & 0.2839 && 0.2762 & 0.1043& & 0.0178 & 0.0092 && 0.0084 & 0.0045 & & 0.1016 & 0.0305 && 0.0275 & 0.0112 \\
    & Spline 
    & 1.5913 & 0.8543 && 0.7995 & 0.3080&
    & 0.2105 & 0.0570 && 0.0570 & 0.0333 & 
    & 0.2170 & 0.0518 && 0.0464 & 0.0100 \\
    & RKHS  & 0.8289 & 0.4630 && 0.4652 & 0.2205& & 0.0264 & 0.0130 && 0.0129 & 0.0069 & & 0.1090 & 0.0401 && 0.0390 & 0.0149 \\
    & Local Linear  
    & 1.3112 & 0.9221 && 0.8854 & 0.6147 &
    & 0.6633 & 0.4535 && 0.5883 & 0.2353 & 
    & 0.2371 & 0.1395 && 0.1472 & 0.1260 \\
    \hline
    \hline
    \multicolumn{2}{c}{}& \multicolumn{5}{c}{Case 4: $d=10$ }&&\multicolumn{5}{c}{Case 4: $d=20$}&&\multicolumn{5}{c}{Case 4: $d=30$}\\
    \cline{3-7}\cline{9-13} \cline{15-19}
    \multicolumn{2}{c}{Sample size $n$} & \multicolumn{2}{c}{$100$}& & \multicolumn{2}{c}{$300$} & &\multicolumn{2}{c}{$100$}& & \multicolumn{2}{c}{$300$} & &\multicolumn{2}{c}{$100$}& & \multicolumn{2}{c}{$300$} \\
    \cline{3-4} \cline{6-7}\cline{9-10} \cline{12-13} \cline{15-16}\cline{18-19}
    \multicolumn{2}{c}{Sampling frequency  $m$} 
    & $10$ &$30$&& $10$ &$30$&& $10$ &$30$&& $10$ &$30$&& $10$ &$30$&& $10$ &$30$\\
    \hline
    &DNN  & 0.0170 & 0.0110 && 0.0111 & 0.0066 & 
    & 0.0370 & 0.0167 && 0.0166 & 0.0074 & 
    & 0.0642 & 0.0268 && 0.0269 & 0.0105 \\
    & RKHS  & 0.1996 & 0.1100 && 0.1104 & 0.0573 & 
    & 0.1820 & 0.1364 && 0.1363 & 0.0754 & 
    & 0.1119 & 0.0967 && 0.0963 & 0.0752 \\
    & Local Linear  & 0.5246 & 0.3767 && 0.3988 & 0.2447 & 
    & 0.4644 & 0.3049 && 0.3305 & 0.2074 & 
    & 0.2644 & 0.2064 && 0.2012 & 0.1790 \\
    \hline
	\end{tabular}
\end{table}
Analyzing the results in Figure \ref{fig:simualtion}, we observe that for a fixed value of sample size $n$, as $m$ increases initially, the error significantly decreases. However, if $m$ continues to increase beyond some points, the gain becomes less significant, which supports our theoretical findings on the phase transition phenomenon. 
Compared to Case 1 with similar series structures, Case 2 has a more explicit additive structure and thus converges faster. 
Case 3 and Case 1 share the same $\fo$, differing only in their input support, and the results exhibit substantially improved convergence rates and a more pronounced phase transition phenomenon, aligning well with our theoretical guidance.
Case 4 benefits from the composition structure that allows these large dimensional functions to converge with relatively small samples. 
Nevertheless, the larger $d$ is, the larger $m$ is needed for the same $n$ to obtain the same convergence error, because the coefficients of the $nm$ term depending on $d$.
Table \ref{tab:bijiao} presents a comparison of the performance of the four methods, highlighting the strengths of the DNN estimator in all cases. As expected by theoretical developments, the DNN estimator is preferred over linear estimators when dealing with unknown hierarchical structures in Case 2. 
With low-dimensional input, the DNN estimator surpasses the local linear estimator, which is also adaptive to unknown manifold structures \citep{bickel2007local}. 
The superiority of the DNN estimator is also evident in the high-dimensional cases (Case 4).
}

\vspace{-0.5cm}
\section{ Real Data Applications}\label{sec:realdata}
\vspace{-0.3cm}
\subsection{Airline Delay Example}
\begin{figure}[!ht]
	\label{fig:example}\spacingset{1.36}
	\centering
	\begin{subfigure}{.61\textwidth}
	   \centering				   \includegraphics[width=\textwidth,height=2.6in]{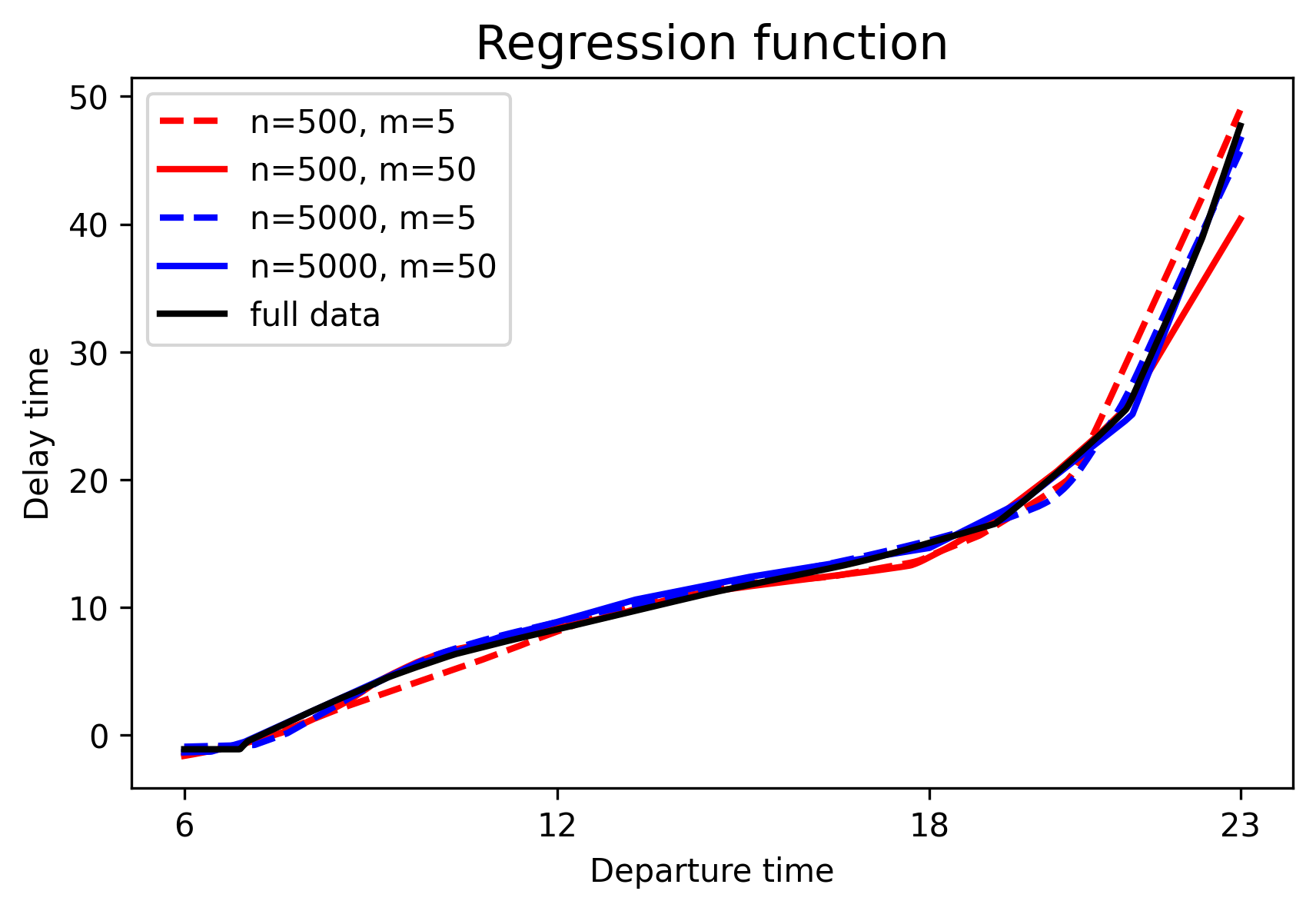} 
	   \caption{Estimated functions}\label{fig:air1}
	\end{subfigure}
	\begin{subfigure}{.38\textwidth}
		\centering
		\includegraphics[width=\textwidth,height=2.6in]{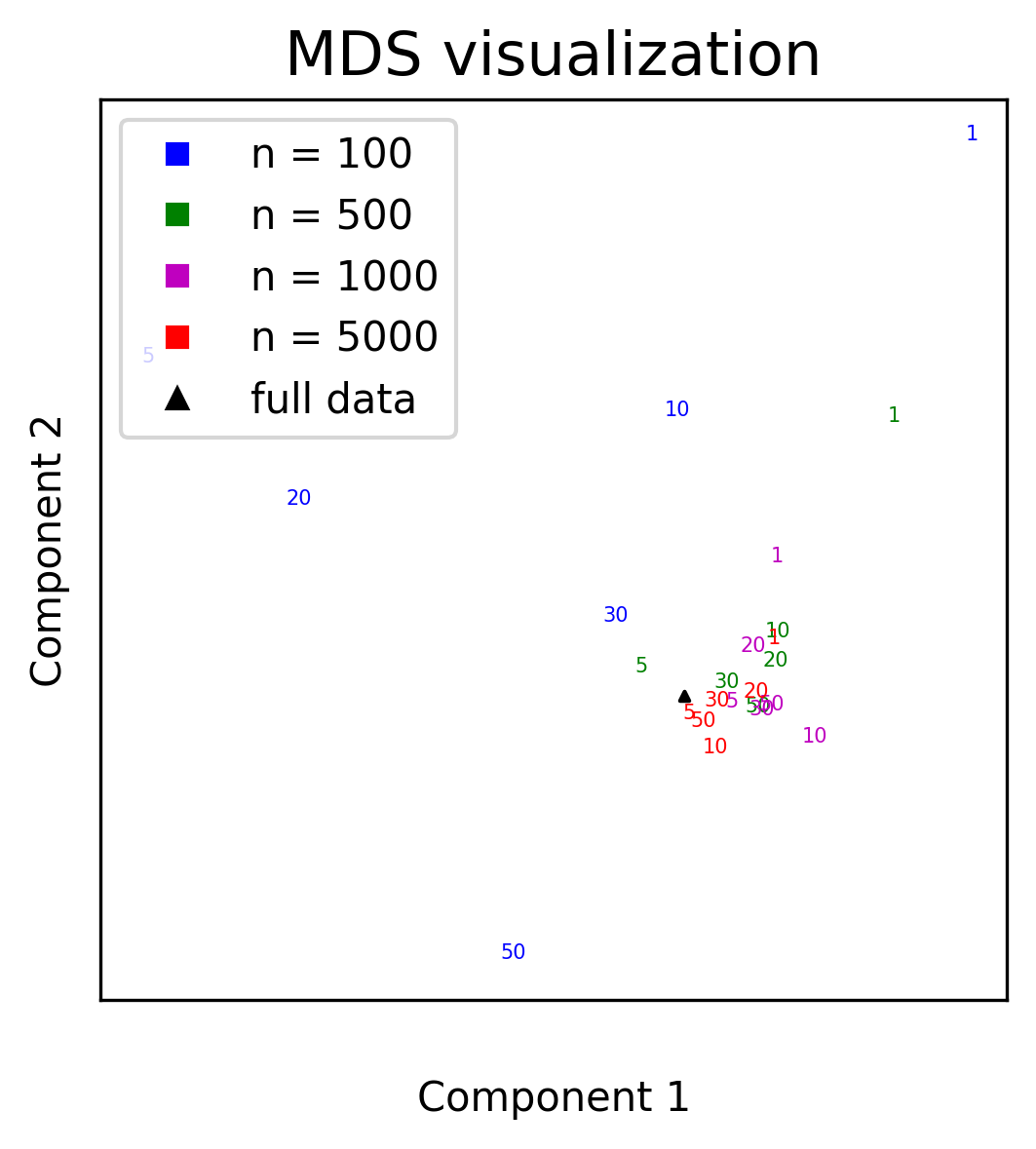}
		\caption{Multiple-dimensional scaling}\label{fig:mds1}
	\end{subfigure}

    \caption{ (a) Graphs depicting estimated function profiles for distinct sample sizes and sampling frequencies. 
    (b) Multiple-dimensional scaling visualization. The numerical values denote $m$, the colors represent $n$, and the black triangles symbolize the estimation with full data.
    }\label{fig:air}
\end{figure}

The airline dataset contains information about commercial flight arrivals and departures in the USA from October 1987 to April 2008, which is publicly available at 
\href{https://community.amstat.org/jointscsg-section/dataexpo/dataexpo2009}{https://community.}
\href{https://community.amstat.org/jointscsg-section/dataexpo/dataexpo2009}{amstat.org/jointscsg-section/dataexpo/dataexpo2009}. 
With nearly 120 million records involving 304 airports, daily records per airport range from 1 to 1269. Our study focuses on delay durations during peak hours (6:00 to 23:00) in 2008, considering airports with at least 50 daily flights in this period. 
Assuming data independence across airports and days, we treat daily data for each airport as individual subjects. Time acts as a covariate, while flight delay times constitute measurements. This aligns with our clustered dependent model, as delay times within the same airport on a given day exhibit correlation. 
There are a total of $8234$ subjects in 2008, we divide the training and validation sets in a $4:1$ ratio. 
The network construction is the same as in Simulation.
For illustration, we randomly sample 50 measurements from each subject as full sample based on which the estimation is taken as the baseline.
To illustrate the proposed method and theory, in a sequentially growing manner, we also train subsamples $n = 100, 500, 1000, 5000$ and sampling frequencies $m = 1,5,10,25,50$, while the rest data with the sample size $0.25 n$ and sampling frequencies $m$ are used as the validation sets. 
The results are shown in Figure \ref{fig:air}.

As depicted in Figure \ref{fig:air1}, estimated functions progressively approach the baseline with larger sample sizes and increased sampling frequencies. 
Figure \ref{fig:mds1} provides a more intuitive confirmation of our points. 
When fixing the number of subjects $n$, increasing $m$ results in shifting towards some sub-center. As $n$ increases, these sub-centers become closer to the baseline obtained from the full data. This provides evidence for the phase transition and certain practical guidance on different roles played by the sample size and sampling frequency.

{\color{black}
\subsection{PM2.5 Example}\label{sec:pm2p5}
The second example concerns the PM2.5 concentration problem with the dataset extracted from \citet{zheng2024dynamic} available at \href{https://github.com/FlyHighest/Dynamic-Synthetic-Control}{https://github.com/FlyHighest/Dynamic-Synthetic-Control}. 
The dataset includes air pollution and meteorological variables of 94 monitoring stations in Beijing and nearby provinces (Tianjin, Hebei, Shandong, Shanxi). 
The variables include pollution concentrations of PM$2.5$, PM$10$, $\text{SO}_2$, $\text{NO}_2$, $\text{O}_3$ and CO, wind speed, humidity, dew point temperature, and air pressure. 
Each station records 72 hourly measurements, with the last 24 taken during Beijing's orange air pollution alert, impacting all stations.  
In this study, each station is taken as a subject, and we focus on the first $48$ records of each station, i.e., $m=48$. 
The concentrations of $\text{SO}_2$, $\text{NO}_2$, $\text{O}_3$ and CO, meteorological conditions including wind speed, humidity, dew point temperature, and air pressure are taken as covariates, while the logarithmic concentration of PM2.5 is the response. Given the complicated relationship between these variables and PM2.5 \citep{chen2017detecting,zou2022estimation}, we frame this as an 8-dimensional regression problem.

For comparison, we consider three commonly used statistical models: linear model, single index model \citep[R package `PLSiMCpp',][]{PLSiMCpp} and additive model \citep[R package `gam',][]{gampackage}, using working independence correlation during estimation.  
Meanwhile, the RKHS estimation \citep{cai2011optimal} and the local linear regression method \citep{zhang2016sparse} are also compared. 
We randomly select $20\%$ of the subjects as test data to compare the performance of the proposed DNN estimator with alternative methods.
Then, the remaining subjects are randomly split into training and validation sets in $4:1$ ratio. The validation set is used for hyperparameter selections, including early stops, network size, and tuning parameters of other methods. The DNN structures are the same as in the simulation. 
Now, with $n_{train} = 62$ and $m =48$, the testing errors of the all methods are shown in Table \ref{tab:pm2.5}.

\begin{table}[!ht] 
\setlength{\tabcolsep}{0.1in}
\spacingset{1.36}
\caption{Mean square testing error of the DNN estimator and the other five estimators.} \label{tab:pm2.5}
\centering
\begin{tabular}{ccccccc}
\hline
& DNN & Linear & Single index & Additive  & RKHS & Local linear\\ 
\hline
Testing error &0.1142 & 0.4449 & 0.1541 & 0.1310 & 0.1520 & 0.1691
\\ 
\hline
\end{tabular}
\end{table}

From Table \ref{tab:pm2.5}, the DNN estimator exhibits the lowest testing error. The linear regression model performs poorly due to bias introduced by its overly strong linear assumptions, and the single index model shows improvement. The additive model performs relatively well, indicating a significant proportion of the effects of covariates on PM2.5 can be attributed separately. However, it still underperforms the DNN estimator, suggesting that interactions between predictors are also important and cannot be ignored.
The RKHS estimator and local linear estimator outperform both the linear regression and single index models due to no model bias. Nevertheless, they still fall short of the performance of the DNN estimator. 

\begin{figure}[!ht]
\spacingset{1.36}
\scriptsize
    \centering
    \includegraphics[width=\textwidth]{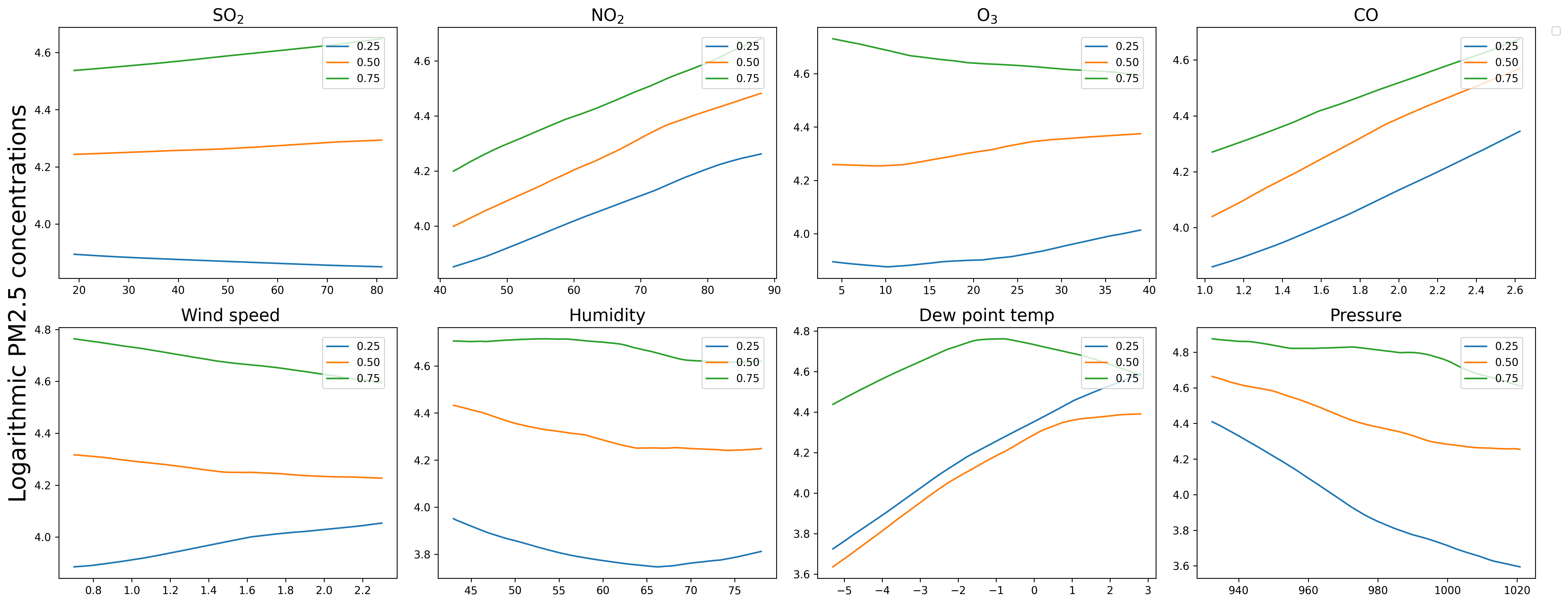}
	\caption{ Plots of the estimated function with respect to each covariate. In each subplot, other covariates are held at their 0.25/0.50/0.75 quantile.}\label{fig:pm}
\end{figure}

Figure \ref{fig:pm} displays the marginal relation between the log PM2.5 concentration and each covariate, with the remaining covariates set to their respective q-quantiles.
For most variables, the response exhibits similar trends across different quantiles, suggesting that the function is closely approximated by an additive model. 
Among the four pollutants, NO$_2$ and CO are most significantly correlated with PM2.5, mainly due to their similar sources in the considered region, including traffic exhaust and industrial emissions. Higher wind speeds tend to homogenize PM2.5 concentrations. 
Humidity and dew point temperature have complex effects on PM2.5 levels: they elevate PM2.5 levels by promoting particle growth and secondary aerosol formation, while also bringing precipitation to wash out particulates from the air. 
Higher air pressure typically means cold and clean air from Siberia, hence leading to lower PM2.5 levels.
}

\section{Conclusion and Discussion}\label{sec:con}
Clustered dependent regression with arbitrary sampling frequency $m$ is a cornerstone for related statistical fields such as longitudinal and functional data analysis. 
We made the attempt to adopt and (thus) allow the use of empirical process tools for general clustered dependence problems, obtaining a series of oracle inequalities. 
These inequalities are new and their proof inspires further studies.
Then a phenomenon, called  ``phase transition'', i.e. repeated measures regression can achieve a parametric rate when increasing the sampling frequency $m$, imparts a unique influence on the asymptotic performance and the curse of dimensionality. 
By considering three further examples, we show the proposed estimator can circumvent the curse of dimensionality, such as achieving a parametric rate with lower sampling frequency and getting faster convergence with sparse design.

The neural networks employed in this article are fully connected feedforward neural networks, while other popular networks are certainly of interest, like convolutional neural networks and deep residual networks. 
More of their properties need to be further understood.
Our oracle inequalities distinguish the approximation error and quantify estimation errors in multiple means (e.g., Rademacher complexity and covering numbers) since we hope that our results will still work with new results about network construction and complexity measure in the future. 
For example, the oracle inequalities could be applied to overparameterized convolutional residual networks with norm constraints using recent advanced results developed in \citet{zhang2024nonparametric}.

\newpage
\setcounter{section}{0}  
\setcounter{equation}{0}
\renewcommand\thesection{S.\arabic{section}}
\renewcommand{\theequation}{S.\arabic{equation}}
\theoremstyle{plain}
\newtheorem{stheorem}{Theorem}
\newtheorem{slemma}{Lemma}
\renewcommand{\theslemma}{S.\arabic{slemma}}
\renewcommand{\thestheorem}{S.\arabic{stheorem}}
\newtheorem{sdefinition}{Definition}
\renewcommand{\thesdefinition}{S.\arabic{sdefinition}}
\newtheorem{sassumption}{Assumption}
\renewcommand{\thesassumption}{S.\arabic{sassumption}}

\theoremstyle{remark}
\renewcommand{\thetable}{S.\arabic{table}}
\renewcommand{\thefigure}{S.\arabic{figure}}

\begin{center}
\LARGE
\textbf{Supplement to ``Deep Regression for
repeated measurements''}
\end{center}

In this supplementary material, we provide the additional results, auxiliary lemmas, and all proofs for the main manuscript. 
Section \ref{sec:cirs} further elucidates more results on circumventing the curse of dimensionality from Section \ref{sec:cir}. 
In Section \ref{app:spline}, we extend the discussion on regression splines, a supplementing result of Section \ref{sec:holderph}.
Additional numerical results are shown in Section \ref{sec: anr}. 
Then Section \ref{app:lem} presents a collection of technical lemmas.  The detailed proofs of theorems and corollaries can be found in Section \ref{app:pftc} and proofs of the the lemmas are given in Section \ref{app:prflem}.

\section{Additional Results on Circumventing Curse of Dimensionality}\label{sec:cirs}

\subsection{Suboptimality of Linear Estimators for Hierarchical Composition Structure}\label{sec:linearnonada}
\citet{schmidtaos} demonstrated that wavelet series estimators are not adaptive to the unknown composition structure in the regression function, even for the single index model. \citet{suzuki2021deep} then extended to the multi-index model and a broader class of linear estimators. 
Here, the linear estimators for cross-sectional regression are any estimators with the form 
$$
\tilde{f}_{L}(\xx)=\sum_{i=1}^n Y_i \varphi_i\left(\xx, X_1,X_2, \ldots, X_n\right), 
$$
where $\varphi_i, i=1,2...,n$ are measurable functions that only depend on $x, X_1, X_2, \ldots, X_n$. 
It is easy to verify that the linear estimator class includes local polynomials, RKHS, regression splines, wavelet series and so on.  
Similarly, one can show, in the repeated measurements setting, the linear estimators
$$
\tilde{f}_{L}(\xx)=\sum_{i=1}^n\sum_{j=1}^m Y_{ij} \varphi_{ij}\left(\xx, X_{11}, \ldots, X_{nm}\right),
$$
are also not adaptive to unknown composition structures. 
For simplicity, we only discuss the single index models.

\begin{stheorem}\label{thm:lowlinear}
Consider the repeated measurement model \eqref{equ:model} and assume that $\mathcal{P}_{\x}$ has a uniform distribution on $\Omega = [0,1]^d$, $d\geq 3$. 
Then we have 
$$
\inf_{\tilde{f}_{L}}\sup_{\mathcal{P}_{Z,\epsilon} \in \mathcal{Q}_{Z,\epsilon}} \mathcal{E}_{\fo}(\tilde{f}_{L})  
\geq \xc\label{c691} \left( n^{-1} + (nm)^{-(2s+d)/(2s+2d)} \right), 
$$
where the infimum is taken over all linear estimators $\tilde{f}_{L}$ based on the available data set $\left\{ \left(Y_{ij},\x_{ij} \right): 1\leq i \leq n, 1\leq j \leq m \right\}$ with form $
\tilde{f}_{L}(\xx)=\sum_{i=1}^n\sum_{j=1}^m Y_{ij} \varphi_{ij}\left(\xx, X_{11}, \ldots, X_{nm}\right)$ and the supremum is taken over all distributions $\mathcal{P}_{Z,\epsilon} \in \mathcal{Q}_{Z,\epsilon}$ satisfying Assumption \ref{ass:subexp} and $\fo(\xx) = g(A^T\xx + b)$ for some $\|g\|_{ \mathcal{C}^{s}([0,1])} \leq 1$, $A\in \mathbb{R}^{d}, b\in \mathbb{R}$ s.t. $A\xx+b \in [0,1]$, $\forall \xx \in [0,1]^d$, and $\jc{c691}$ is a positive number not depending on $n$ and $m$. 
\end{stheorem}

Theorem \ref{thm:lowlinear} provides a lower bound that the convergence rate of linear estimators cannot be faster than $ n^{-1} + (nm)^{-(2s+d)/(2s+2d)}$, which is slower than the optimal rate $ n^{-1} + (nm)^{-2s/(2s+1)}$ achieved by the DNN estimator given in Theorem \ref{thm:hieup} when $2s+d<2sd$. 
This illustrates that the linear estimators are not adaptive to the hierarchical composition structure, i.e., the minimax rate cannot be achieved when the structure is unknown. 
Therefore, traditional methods often assume a specific structure of the real model and then design the tailored estimation method for that known structure.
In contrast, 
the DNNs can adapt to unknown low intrinsic dimension structures through optimization because the minimizer can learn the unknown hierarchical structures without specific knowledge and designs.

\subsection{Anisotropic H\"older Smoothness} \label{sec:anisoh}
In the previous subsection, we have studied that the function of interest obeys a general presumed structural assumption and has shown success in theoretical terms. 
However, assuming a hierarchical composition structure for each regression problem is impractical. To confront this challenge, another effective approach is to consider the functions that exhibit varying smoothness along different coordinate directions, as improving smoothness isotropically can be a daunting task  \citep{Anirban2014anisot}.
This is also practically sensible, given that the coordinates may have distinct meanings, such as spatial and temporal features, thus different smoothness is often preferred for different coordinates.
Based on these considerations, we introduce the definition of anisotropic H\"older space as follows.

\begin{sdefinition} 
Let $\boldsymbol{s} =\left( s _{1}, s _{2}, \ldots,  s _{d}\right),  s _{i}>0$ be the smoothness index. 
Assume $d\geq 2$, the anisotropic $\boldsymbol{s}$-H\"older norm of a function $f:\Sigma \in \mathbb{R}^d  \rightarrow\mathbb{R}$ is defined as 
$$
\|f\|_{\mathcal{C}^{\boldsymbol{s}}}=
\sum_{i=1}^{d}\sum_{\alpha< s_{i}}\left\|\partial^{\alpha}_{i} f\right\|_{\infty}
+ \sum_{i=1}^{d}\sum_{\alpha = \lfloor s_{i} \rfloor}
\sup _{\substack{\xx, \xx' \in \Sigma \\  \xx' = \xx + te_{i}, t\neq 0}} \frac{\left|\partial^{\alpha}_{i} f(\xx)-\partial^{\alpha}_{i} f(\xx')\right|}{\|\xx-\xx'\|^{s_{i}-\lfloor s_{i} \rfloor}},
$$
where $\partial^{\alpha}_{i}$ is the $\alpha$-th order partial derivative of $f$ with respect to the $i$-th component $e_{i}$ with $\alpha \in \mathbb{N}$. Furthermore, the  anisotropic H\"older space $\mathcal{C}^{\boldsymbol{s}}\left(\Sigma\right)$ consists of all functions for which the anisotropic $\boldsymbol{s}$-H\"older norm is finite.
\end{sdefinition}

Some existing works consider similar anisotropic smoothness functions for cross-sectional regression.
\citet{HANG2021337} investigated the nonparametric regression problem using SVMs with anisotropic Gaussian RBF kernels and established the optimal convergence rates. 
They assume knowing function smoothness in each direction to establish the theory, and their method needs to tune $(d+1)$ parameters.
Recently, \citet{suzuki2021deep} developed nonparametric regression theory via sparsely connected DNNs. 
From the previous literature, the optimal rate for cross-sectional regression is $O(n^{-2\ts/(2\ts + d)})$, where $\ts$ is the harmonic mean of the smoothness in each direction
$$
\ts := \left(\frac{1}{d}\sum_{i=1}^{d}\frac{1}{s_i} \right)^{-1}. 
$$
Now, we first build the approximation theorem of fully connected DNNs for anisotropic H\"older function.

\begin{stheorem}[Approximation for anisotropic H\"older function]\label{thm:aniso}
Suppose that $f \in \mathcal{C}^{\bs}\left(\Omega\right)$ with the anisotropic $\boldsymbol{s}$-H\"older norm $\| f  \|_{\mathcal{C}^{\bs}}\leq \xdb\label{anisoholdernorm}$ for some $\jdb{anisoholdernorm}>0$ and  $\|f\|_{\infty}\leq\jdb{b1}$. 
Then, for any $L,W\geq 3$, there is some neural network 
$f_{\nn}\in \nn\left(d, \xc\label{c66} L \log L, \xc\label{c67} W \log W, \jdb{b1}\right)$
satisfying 
$$
|f_{\nn}(\xx)-f(\xx)|\leq \xc\label{c68}\left(LW\right)^{-2 \tilde{s} / d},
$$
for each $\xx\in \Omega$, where
$\jc{c66}$ and $\jc{c67}$ are constants not depending on $L$, $W$ and $f$, and 
$\jc{c68}$ is a constant free of $L$ and $W$.

\end{stheorem}

By assuming the anisotropic H\"older smoothness, we obtain the non-asymptotic upper bound and the matching minimax lower bound that lead to the optimal convergence rate for the repeated measurements model.

\begin{sassumption}[Anisotropic H\"older smoothness] \label{ass:anisoholder}
The target function $\fo(\xx)$ belongs to the anisotropic H\"older space $\mathcal{C}^{\bs}\left(\Omega\right)$ for a given positive  $\bs$ and $\|\fo \|_{\mathcal{C}^{\bs}}\leq \jdb{anisoholdernorm}$. 
\end{sassumption}

\begin{stheorem}\label{thm:anisoholdernn}
Consider model \eqref{equ:model} and the estimator $\hat{f}_{nm}$ obtained by \eqref{equ:estimator}. Suppose that Assumptions \ref{ass:fobound},
\ref{ass:subexp} and \ref{ass:anisoholder} hold, and set 
$\mathcal{F}_{nm} = \nn\left(d, \jc{c66} L \log L, \jc{c67} W \log W, \jdb{b1}\right).$
 If we specify the neural network in a flexible way satisfying 
$
LW=\lfloor \xc\label{c69} (nm)^{ d/(4\tilde{s}+2d)  } (\log nm)^{ -4d/(2\tilde{s}+d)} \rfloor
$ 
for some certain constant $\jc{c69}$ not depending on $n$ and $m$, then 
$$
\pef \leq \xc\label{c70} ( n^{-1}+ (nm)^{-2\tilde{s}/(2\tilde{s}+d)}(\log nm)^{ 16\tilde{s}/(2\tilde{s}+d)} ), 
$$
where $\jc{c70}$ is a constant free of $n$ and $m$. 
\end{stheorem}

\begin{stheorem}\label{thm:anisolow}
Consider the repeated measures model \eqref{equ:model} and assume that $\mathcal{P}_{\x}$ has a uniform distribution on $\Omega$.
Under the Assumptions \ref{ass:subexp} and \ref{ass:anisoholder}, we have
$$
    \inf_{\tilde{f}}\sup_{  \mathcal{P}_{Z,\epsilon} \in \mathcal{Q}_{Z,\epsilon} }  \mathcal{E}_{\fo}(\tilde{f})
\geq \xc\label{c71} (n^{-1} + (nm)^{-2\tilde{s}/(2\tilde{s}+d)}), 
$$
where the infimum is taken over all possible estimators $\tilde{f}$ based on the available data set $\left\{ \left(Y_{ij},\x_{ij} \right): 1\leq i \leq n, 1\leq j \leq m \right\}$, the supremum is taken over all distributions $\mathcal{P}_{Z,\epsilon} \in \mathcal{Q}_{Z,\epsilon}$ satisfying Assumption \ref{ass:subexp} and \ref{ass:anisoholder}, 
and $\jc{c71}$ is a positive number not depending on $n$ and $m$. 
\end{stheorem}

Similar conclusions can be drawn from these theorems that our ReLU DNN estimator attains the minimax optimal rate $O(n^{-1} + (nm)^{-2\tilde{s}/(2\tilde{s}+d)})$ up to logarithmic factors for the repeated measurement model \eqref{equ:model}, with the phase transition occurs at $m\asymp O(n^{-d/2\tilde{s}})$. It is worth noting that in the anisotropic framework $\fo \in \mathcal{C}^{\bs} $ in contrast to $\fo \in \mathcal{C}^{\sm}$ in the isotropic case with the convergence rate $O(n^{-1} + (nm)^{-2\sm/(2\sm+d)})$ by Theorem \ref{thm:holdernn}. 
Knowing $\min_{1\leq i\leq d} s_{i}\leq \tilde{s}\leq \max_{1\leq i\leq d} s_{i}$, the harmonic mean $\tilde{s}$ can be significantly larger than $\sm$ when the smoothness is highly anisotropic.

\subsection{Low-dimensional Support} \label{sec:lowip}

For completeness of the paper, we first briefly introduce some of the basic concepts about manifolds used. 

\begin{sdefinition}[Chart]
A chart for $\mathcal{M}$ is a pair $(U, \phi)$ such that $U \subset \mathcal{M}$ is open and $\phi: U \to \mathbb{R}^d$, where $\phi$ is a homeomorphism, i.e. $\phi$ is a bijective, both $\phi$ and the inverse of $\phi$ are continuous. 
\end{sdefinition}
\begin{sdefinition}[$C^k$ compatible]
We say two charts $(U, \phi)$ and $(V, \psi)$ on $\mathcal{M}$ are $C^k$ compatible if and only if 
$$
\phi \circ \psi^{-1}: \psi(U \cap V) \mapsto \phi(U \cap V) \quad \text { and } \quad \psi \circ \phi^{-1}: \phi(U \cap V) \mapsto \psi(U \cap V)
$$
are both $C^k$. 
\end{sdefinition}
\begin{sdefinition}[Atlas]
An atlas for $\mathcal{M}$ is a collection $\left\{\left(U_i, \phi_i\right)\right\}_{i \in I}$ of charts such that $\bigcup_i U_i=\mathcal{M}$. 
Then a $C^k$ atlas for $\mathcal{M}$ is an atlas consisting of pairwise $C^k$ compatible charts. 
\end{sdefinition}
\begin{sdefinition}[Smooth manifold and Riemannian manifold]
 A smooth manifold is a manifold together with a $C^{\infty}$ atlas. A Riemannian manifold is a smooth manifold $\mathcal{M}$ endowed with Riemannian metric. 
\end{sdefinition}
\begin{sdefinition}[H\"older space $\mathcal{C}^{s}\left( \mathcal{M} \right)$]
Let $\mathcal{M}$ be a $d_{\mathcal{M}}$-dimensional Riemannian manifold isometrically embedded in $\mathbb{R}^d$. Let $\left\{\left(U_i, \phi_i\right)\right\}_{i \in \mathcal{I}}$ be a $C^{\infty}$ atlas of $\mathcal{M}$ where the $\phi_i^{\prime}s$ are orthogonal projections onto tangent space. 
Let $s$ be a positive number.
We say $f\in \mathcal{C}^{s}\left( \mathcal{M} \right)$ if and only if, for any chart $\left(U_i, \phi_i\right)$, the composition $f \circ \phi^{-1}: \phi(U) \mapsto \mathbb{R}$ satisfies $f \circ \phi^{-1} \in \mathcal{C}^{s}\left( \phi(U) \right)$. 
\end{sdefinition}

It is important to note that the above definition does not depend on the choice of the atlas $\left\{\left(U_i, \phi_i\right)\right\}_{i \in \mathcal{I}}$. Given chart $(U, \phi)$, consider an alternate chart $(V, \psi)$ such that the intersection $V \cap U$ is not empty.
In this context, the function $f \circ \psi^{-1}$ can be decomposed as $\left(f \circ \phi^{-1}\right) \circ \left(\phi \circ \psi^{-1}\right)$. 
Because $f \circ \phi^{-1}$ is s-H\"older smooth and $\phi \circ \psi^{-1}$ is $C^{\infty}$ smooth, their composition is still $\mathcal{C}^{s}$ smooth.

We also extend to the approximate manifold and Minkowski dimension set, coupling with the latest approximation theory of DNN estimators based on fully connected ReLU neural networks \citep{huangjian}.  
We keep Assumption \ref{ass:holder} that $\fo \in \mathcal{C}^{s}\left( [0,1]^d\right)$. Based on it, we assume that the predictor $X$ is supported on some lower dimensional sets.

Let $\mathcal{M}$ be a compact $d_{\mathcal{M}}$-dimensional Riemannian submanifold and define the approximate manifold as 
$$
\mathcal{M}_\rho=\left\{x \in[0,1]^d: \inf \left\{\|x-y\|_2: y \in \mathcal{M}\right\} \leq \rho\right\}, \rho \in(0,1). 
$$
Since, in practice, observations may not fall strictly on a manifold but still tend to have a low-dimensional structure, as the empirical studies in \citet{carlsson2009topology} reveal for image data. 
Combining Corollary \ref{cor:networksize} and the approximation theory in the proof of Theorem 6.1 in \citet{huangjian}, we have the following result.

\begin{sassumption}[Approximate manifold] \label{ass:nearmanifold}
Assume the predictor $X$ is supported on an approximate manifold $\mathcal{M}_\rho$. 
\end{sassumption}

\begin{stheorem}\label{thm:nearmanifold}
Consider the repeated measures model \eqref{equ:model} and the estimator $\hat{f}_{nm}$ obtained by \eqref{equ:estimator}. Suppose that Assumptions \ref{ass:fobound}, 
\ref{ass:subexp}, \ref{ass:holder} and \ref{ass:nearmanifold} hold, and set the candidate class
$$\mathcal{F}_{nm} = \nn\left(d, \xc\label{c680} L \log L,  \xc\label{c681} W \log W, \jdb{b1}\right).
$$ 
with $\jc{c680}, \jc{c681}$ not depending on $L, W$.  
If we specify the neural network in a flexible way satisfying 
$$
LW=\lfloor \xc\label{c682}  (nm)^{d_{\delta}/(4s+2d_{0})} (\log nm)^{4d_{\delta}/(2s+d_{\delta}) } \rfloor
$$ with some constant $\jc{c682}$ not depending on $m$ and $n$, 
and $d_\delta=O\left(d_{\mathcal{M}} \log (d / \delta) / \delta^2\right)$ is an integer such that $d_{\mathcal{M}} \leq d_\delta<d$ for any $\delta \in(0,1)$, 
then 
$$
\pef \leq \xc\label{c683} ( n^{-1}+ (nm)^{-2s/(2s+d_{\delta})}(\log nm)^{16s/(2s+d_{\delta})} ), 
$$
where $\jc{c683}$ is a constant free of $n$ and $m$.
\end{stheorem}

Another effective method to measure the intrinsic complexity of a set is the Minkowski dimension. 
Given a set $\Sigma \subseteq \mathbb{R}^{d}$, define the Minkowski dimension as
$$
\operatorname{Mdim}(\Sigma):=\lim _{\varepsilon \rightarrow 0} \frac{\log \mathcal{N}\left(\varepsilon,\Sigma, L^2 \right)}{\log (1 / \varepsilon)},
$$
if the limit on the right side exists.
The Minkowski dimension characterizes the limiting behavior of the covering number when the radius of the balls converges to zero, which provides an intuitive measure for the complexity of a set. It is worth noting that the Minkowski dimension does not depend on smoothness, thus it can measure the intrinsic dimension of some highly non-smooth sets. 

\begin{sassumption}[Minkowski dimension] \label{ass:f0minkowski}
The predictor $\x$ is supported on a set $A\subset\Omega$ with Minkowski dimension $d^{*}$ less than $d$. 
\end{sassumption}

\begin{stheorem}\label{thm:minkup}
Consider the repeated measures model \eqref{equ:model} and the estimator $\hat{f}_{nm}$ obtained by \eqref{equ:estimator}. Suppose that Assumptions \ref{ass:fobound}, 
\ref{ass:subexp}, \ref{ass:holder} and \ref{ass:f0minkowski} hold, and set the candidate class
$$\mathcal{F}_{nm} = \nn\left(d, \xc\label{c80} L \log L,  \xc\label{c81} W \log W, \jdb{b1}\right).
$$ 
with $\jc{c80}, \jc{c81}$ not depending on $L, W$.  
If we specify the neural network in a flexible way satisfying 
$$
LW=\lfloor \xc\label{c82}  (nm)^{d_{0}/(4s+2d_{0})} (\log nm)^{4d_{0}/(2s+d_{0}) } \rfloor
$$ with some constant $\jc{c82}$ not depending on $m$ and $n$, and $d \geq d_0 \geq \kappa d^* / \delta^2=O\left(d^* / \delta^2\right)$ for $\delta \in(0,1)$ and some constant $\kappa>0$, then 
$$
\pef \leq \xc\label{c83} ( n^{-1}+ (nm)^{-2s/(2s+d_{0})}(\log nm)^{16s/(2s+d_{0})} ), 
$$
where $\jc{c83}$ is a constant free of $n$ and $m$.
\end{stheorem}

Note that the convergence rate depends on $d_{\delta}=O(d_{\mathcal{M}}\log(d))$ under the approximate manifold assumption and $d_0=O\left(d^* / \delta^2\right)$ under the low Minkowski dimension assumption, respectively, which are smaller than $d$ but still greater than the intrinsic dimension. These error bounds may be sharpened given improved approximation theory that remains an open problem in theoretical deep learning.

\section{Results on Regression Splines}
\label{app:spline}

We present here a typical regression spline formulation in the form of the cardinal B-splines, which is also used extensively in later proofs.
Let the univariate function $B(x)=I(x \in[0,1])$ be the indicator function on for the interval $[0,1]$, and the cardinal B-splines of order $r$ are defined recursively by
$$
 B_{1}(x)=B(x), \quad B_{r}(x)=B_{r-1}(r) * B(x),
$$
where $*$ denotes convolution, i.e. 
$f * g(x):=\int f(x-t) g(t) \mathrm{d} t$. 
Then an explicit expression of $B_{r}(x)$ is given by 
\begin{equation}\label{equ:301}
B_{r}(x)=\frac{1}{(r-1) !} \sum_{i=0}^{r}(-1)^{i}
\binom{r}{i}
(x-i)_{+}^{r-1},        
\end{equation}
which is a piecewise polynomial of degree $(r-1)$ with the support  $[0,r]$ and $\|B_{r}\|_{\infty}=1$.
In multivariate situation, for $\bk=(k_{1}, k_{2}, ..., k_{d})$ and $\bl=(l_{1}, l_{2}, ..., l_{d})$, we define 
\begin{equation}\label{equ:302}
\mathcal{B}^{d}_{\bk,\bl}=\prod_{i=1}^{d} B_{r}( k_{i} x_{i}-l_{i}).    
\end{equation}
as the multivariate cardinal B-splines, where $r$ is suppressed for brevity.
For the isotropic case $k\in \mathbb{N}$, we slightly abuse the notation and let
$
\mathcal{B}^{d}_{k,\bl}=\prod_{i=1}^{d} B_{r}( k x_{i}-l_{i})
$ when all components of $\bk$ equal to $k$.
The following approximation result on splines for H\"older functions is well known. 

\begin{slemma}\label{lem:splineholder}
Suppose that $f \in \mathcal{C}^{s}\left(\Omega\right)$ and $\| f  \|_{\mathcal{C}^{s}}\leq \jdb{holdernorm}$. 
Let $k$ be a positive integer, 
there are $\left\{ a_{k,\bl} \right\}$ with $j\in \mathbb{Z}_{+}$ and multi-index 
$\bl \in \prod_{i=1}^{d} \left\{-r+1,-r+2, \ldots, k-1,
k\right\},
$
such that
$$
\left\|f-  \sum_{\bl} a_{k,\bl}  \mathcal{B}^{d}_{k,\bl}  \right\|_{\infty} \leq \xc\label{c360} k^{-s} \| f  \|_{\mathcal{C}^{s}},  
$$
where $a_{k,\bl}\leq \xc\label{c361}\| f  \|_{\mathcal{C}^{s}}$ for any $k$ and $\bl$,   
$\jc{c360}$ and $\jc{c361}$ are constants no depending on $k$ and $f$. 
\end{slemma}

 Note that the theory in this paper requires the functions in $\mathcal{F}_{nm}$ to be bounded. A simple way is to let $\mathcal{F}_{nm}$ be the truncated version of the functional space  $\left\{ \sum_{\bl} \alpha_{k,\bl}  \mathcal{B}^{d}_{k,\bl} : a_{k,\bl}\in \mathbb{R} \right\}$ for suitable $k$, which brings no significant change to the VC dimension of $\mathcal{F}_{nm}$. 

Another remedy is to set $\mathcal{F}_{nm}$ as $\left\{ \sum_{\bl} \alpha_{k,\bl}  \mathcal{B}^{d}_{k,\bl} : a_{k,\bl}\leq \jc{c361}\| f  \|_{\mathcal{C}^{s}} \right\}$. 
However, the upper bound $\jc{c361}\| f  \|_{\mathcal{C}^{s}} $ is usually unknown, which can be circumvented by relaxation to an increasing sequence, i.e.,
\begin{equation}\label{equ:splineinc}
    \mathcal{F}_{nm}=\left\{ \sum_{\bl} \alpha_{k,\bl}  \mathcal{B}^{d}_{k,\bl} : a_{k,\bl}\leq \log n \right\}.
\end{equation}
Recall the definition of the regression spline estimator as in \eqref{equ:estimator}, while $\log n$ does not introduce much computational complexity to the optimization. 
Lemma \ref{lem:splineholder} implies the approximation error is of order $O(k^{-2s})$, 
and a direct application of Corollary \ref{cor:cor1} yields the following result.
\begin{stheorem}
    Consider the repeated measures model \eqref{equ:model} and the estimator $\hat{f}_{nm}$ obtained by \eqref{equ:estimator}. Suppose that Assumptions \ref{ass:fobound}, 
\ref{ass:subexp} and \ref{ass:holder} hold, and set the candidate class
as \eqref{equ:splineinc}, then
$$
\pef \leq \xc\label{c690} (k^{-2s} + n^{-1}(\log n)^2 + (nm)^{-1} k^d (\log nm)^3),
$$
for sufficiently large $n$, where \jc{c690} is a constant free of $n$, $m$ and $k$. 
\end{stheorem}
Hence we choose $k\asymp O( (nm)^{1/(2s+d)}(\log nm)^{-3/(2s+d)} ) $ for trade-off, leading to the prediction error  $O( n^{-1}(\log n)^2 + (nm)^{-2s/(2s+d)}(\log nm)^{6s/(2s+d)} )$ which matches the lower bound established in Theorem \ref{thm:low1}.

\section{Additional Numerical Results} \label{sec: anr}
In this study, the Adam optimizer with a maximum iteration limit of 200 epochs is used for network training.  We only employ the early stopping technique with patience $10$ for regularization. The network configurations are setting to be $(L,W) = (1,50),(2,100),$$(3,200),$$(4,400),$ $(5,600),(5,800)$. Each example trains all networks with different depths and widths, then we select the model with the minimal validation error. Note that the proposed theory in this work assumes the order of $LW$. In the simulation study, $LW$ ranges from $50$ to $4000$, which is in the typical selection range of tuning parameters.

\begin{figure}[!ht]
\spacingset{1.36}
\scriptsize
    \centering
    \includegraphics[width=\textwidth]{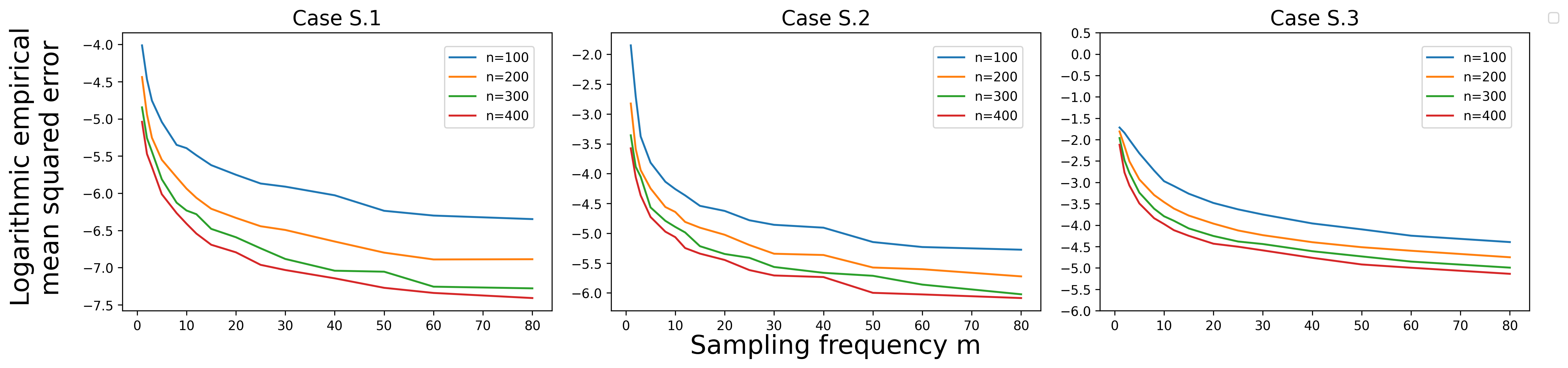}
	\caption{Finite sample mean squared error of DNN estimators for $n=100,200,300,400$ and sampling frequencies $m=1,2,3,5,8,10,12,15,20,25,30,40,50,60,80$. Each subplot portrays sampling frequencies ($m$) on the horizontal axis and logarithmic empirical mean squared error on the vertical axis. }\label{fig:simualtions}
\end{figure}

For all simulations, the $d$-dimensional covariate $X$ is generated from the uniform distribution on $\Omega \subseteq [0,1]^d$.
Recalling the model setting $Y_{ij}=\fo(\x_{ij})+U_{i}(\x_{ij})+\epsilon_{ij}$, we take the individual stochastic part  $U_{i}(\xx)=\sum_{k=1}^{\infty}\sum_{l=1}^{d}\frac{\sqrt{3}\xi_{ikl}}{\sqrt{d}k}\cos(k\pi\xx_{l})$ with independent scores $\xi_{ikl} \sim N(0,0.1^2)$, and measurement errors $\epsilon_{ij} \sim N(0,0.01^2)$.
Now, we add three examples, 
\begin{align*}
& \text { Case S.}1:  \fo(x) =  x_1^2x_2^3 + \log(1+x_3) + \sqrt{1+x_4x_5} +\exp(x_5/2),\Omega = [0,1]^5,\\
& \text { Case S.}2: \fo(x) =  (x_1^2x_2^3 + \log(1+x_3) + \sqrt{1+x_4x_5} +\exp(x_5/2))^2 - 12 ,\Omega = [0,1]^5, \\
& \text { Case S.}3: \fo(x) = \sum_{ k_1,...,k_5\geq 1} a^{ \max\{k_1,2k_2,2k_3,2k_4,2k_5\}} \prod_{l=1}^{5}\cos(2\pi l^{-1} b^{k_l} \xx_l), \Omega = [0,1]^5. 
\end{align*}
Case S.1 and Case S.2 have been previously considered by \citet{deepcox} for the hierarchical model.
Case S.3 serves as an anisotropic smooth counterpart of Case 1 because the coefficient decay rate in the last four components is higher.
Results of the DNN estimator are shown in Figure \ref{fig:simualtions}, and the comparisons with the regression spline method in Section \ref{app:spline}, RKHS method in \citet{cai2011optimal}, and local linear method in \citet{zhang2016sparse}, are in Table \ref{tab:bijiaos}. 
Here, the RKHS method uses the comment Laplace kernel, and the local linear regression uses the comment Epanechnikov Kernel; the tuning parameters are selected by the validation set.

\begin{figure}[!ht]
\spacingset{1.36}
\scriptsize
    \centering
    \includegraphics[width=0.6\textwidth]{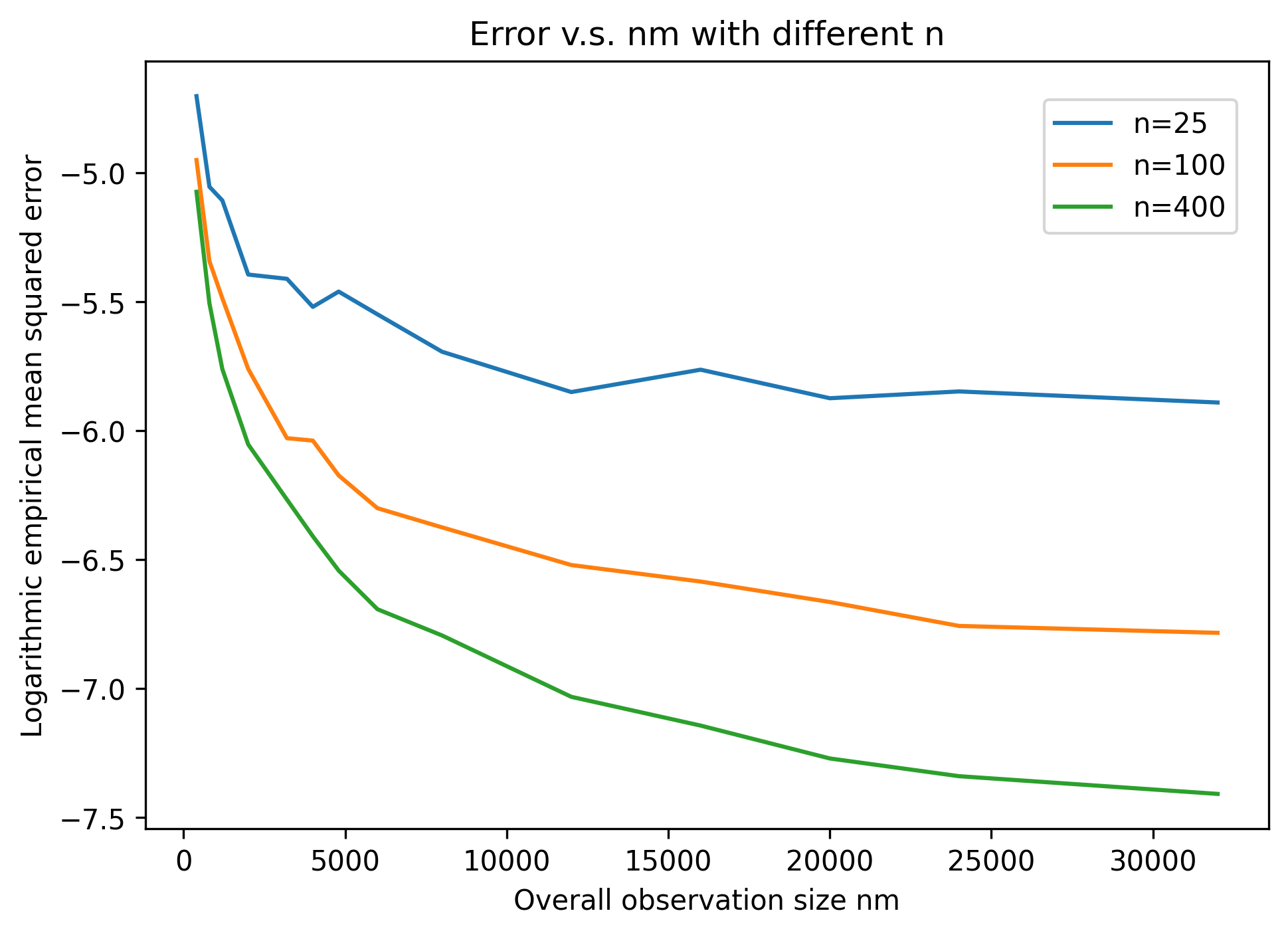}
	\caption{Finite sample mean squared error of the  DNN estimators for sample size $n=25,100,400$ and overall observation size $nm=400,800,1200,2000,3200,4000,4800,6000,8000,12000,16000,20000,24000,32000$. Each subplot portrays the overall observation size ($nm$) on the horizontal axis and the logarithmic empirical mean squared error on the vertical axis. }\label{fig:simualtionss}
\end{figure}

\begin{table}[ht]
\spacingset{1.36}
\centering
\scriptsize
    \caption{Finite sample mean square error of the DNN estimator, naive regression spline estimator (for cases with $d\leq 5$), RKHS estimator and local linear regression estimator. }
    \setlength{\tabcolsep}{0.03in}
	\label{tab:bijiaos}
 \begin{tabular}{ccccccccccccccccccc}
		\hline
		  \multicolumn{2}{c}{}& \multicolumn{5}{c}{Case S.1}&&\multicolumn{5}{c}{Case S.2}&&\multicolumn{5}{c}{Case S.3}\\
    \cline{3-7}\cline{9-13} \cline{15-19}
    \multicolumn{2}{c}{Sample size $n$} & \multicolumn{2}{c}{$100$}& & \multicolumn{2}{c}{$300$} & &\multicolumn{2}{c}{$100$}& & \multicolumn{2}{c}{$300$} & &\multicolumn{2}{c}{$100$}& & \multicolumn{2}{c}{$300$} \\
    \cline{3-4} \cline{6-7}\cline{9-10} \cline{12-13} \cline{15-16}\cline{18-19}
    \multicolumn{2}{c}{Sampling frequency  $m$} 
    & $10$ &$30$&& $10$ &$30$&& $10$ &$30$&& $10$ &$30$&& $10$ &$30$&& $10$ &$30$\\
    \hline
    &DNN  & 0.0045 & 0.0027 && 0.0020 & 0.0010& 
    & 0.0141 & 0.0078 && 0.0075 & 0.0038 & 
    & 0.0513 & 0.0235 && 0.0226 & 0.0118 \\
    & Spline & 0.1143 & 0.0177 && 0.0182 & 0.0045& 
    & 0.1108 & 0.0178 && 0.0183 & 0.0045& 
    & 0.2603 & 0.0691 && 0.0690 & 0.0315 \\
    & RKHS  & 0.0074 & 0.0039 && 0.0037 & 0.0019& 
    & 0.0339 & 0.0159 && 0.0159 & 0.0079 & 
    & 0.0819 & 0.0493 && 0.0488 & 0.0280 \\
    & Local Linear  & 0.7918 & 0.3556 && 0.3438 & 0.1099 & 
    & 1.1696 & 0.7856 && 0.6618 & 0.2783 & 
    & 0.1615 & 0.1320 && 0.1200 & 0.0941 \\
    \hline
	\end{tabular}
\end{table}
Similarly with Figure \ref{fig:simualtion} in the paper, Figure \ref{fig:simualtions} also reflects the phase transition as $n$ and $m$ grow.
It is observed that the DNN estimators learn these functions significantly better than Case 1 in the paper, due to their hierarchical structure and higher smoothness. Table \ref{tab:bijiaos} then demonstrates the outstanding performance of the DNN estimators compared to other methods. 
Finally, to further verify the rate more comprehensively, we add the verification for the influence of total observations $nm$ and sample size $n$. In particular, Figure \ref{fig:simualtionss} corresponds to Case S.1 with the sample size $n=25,100,400$ and the overall observation size $nm=400,800,1200,2000,3200,4000,4800,6000,8000,12000,16000,20000,24000,32000$. The figure indicates that larger values of $n$ lead to better estimation even when $nm$ is the same, which provides more comprehensive empirical validation. 

\section{Technical Lemmas}\label{app:lem}

This supplement contains the lemmas that are crucial for proof.
Lemmas \ref{lem:lema1}-\ref{lem:lema2} are crucial findings in the proof of Theorem \ref{thm:thm1} without involving cluster dependence.
Lemmas \ref{lem:bloc}-\ref{lem:bspl} play significant roles in the approximation theory of fully connected neural networks, while Lemmas \ref{lem:ltc}-\ref{lem:quj} serve as fundamental components for the empirical process. 
Additionally, Lemma \ref{lem:lowerlinearlemma} is to prove the suboptimality of the linear estimator class for unknown hierarchical structures.

\begin{slemma}\label{lem:lema1}
    Under the same model and assumptions as Theorem \ref{thm:thm1}, we have 
    $$
    \begin{aligned}
         &\left| \chf -\frac{1}{n m} \sum_{i=1}^{n} \sum_{j=1}^{m}\left(\hat f\left(\x_{i j}\right)- \fo\left(\x_{i j}\right)\right)^{2}\right| \\
        \leq & \frac{1}{2}\chf  +\xc\label{c801}\jdb{b1}^2r^{*}+\frac{\xc\label{c802}\jdb{b1}^2t}{nm}, 
    \end{aligned}
    $$
    with probability at least $1-\exp\{-t\}$, where \jc{c801} and \jc{c802} are universal constants.
\end{slemma}

\begin{slemma}\label{lem:lema2}
    Under the same model and assumptions as Theorem \ref{thm:thm1}, we have 
    $$
    \begin{aligned}
    & \left| 
     \frac{1}{n m} \sum_{i=1}^{n} \sum_{j=1}^{m} \left(\mathcal{T}_{\beta_{\epsilon}} \epsilon_{ij}-\mathbb{E}\left[ \mathcal{T}_{\beta_{\epsilon}} \epsilon_{ij} \right]\right)\left(\hat{f}\left(\x_{i j}\right)-\fo\left(\x_{i j}\right)\right)
    \right|\\
     \leq &\frac{1}{16}\chf  +\xc\label{c803}\beta_{\epsilon}^2r^{*}+\frac{\xc\label{c804}(\jdb{b1}\beta_{\epsilon}+ \jdb{b3}^{2}  ) t }{nm}.
     \end{aligned}
      $$
    with probability at least $1-\exp\{-t\}$, where
    $\mathcal{T}_{\beta_{\epsilon}}$ is the truncation operator, 
    \jc{c803} and \jc{c804} are universal constants.
\end{slemma}

\begin{slemma}\label{lem:bloc}
Given positive integer $L,W$ and $K\leq L^2W^2$, assume positive number $ \delta \leq \frac{1}{3 K}$. Then there is some ReLU neural network function $f_{\nn}\in \nn\left(1,  4L+5 ,  4W+3 \right)$ that satisfies 
$$
f_{\nn}(\xx)=k, \quad \text {when } x \in\left[\frac{k}{K}, \frac{k+1}{K}-\delta \cdot I_{\{k \leq K-2\}}\right], 
$$
for $k=0,1, ..., K-1$.
\end{slemma}

\begin{slemma}[Proposition 4.4 in \citet{lu2021deep}]\label{lem:pof}
Given positive integer $L, W$ and $\lceil s \rceil $. 
Suppose that  $y_{l}=0$ or $1$ for $l=0,1,2,...,  L^{2}W^{2}-1$. Then  there exists a ReLU neural network function $f_{\nn}\in \nn\left(1,  5(L+2)\log_{2}(4L) ,  16 \lceil s \rceil(W+1) \log_{2}(8W)\right)$  that satisfies 
$$
\left|f_{\nn}(l)-y_{l}\right| \leq (LW)^{-2 \lceil s \rceil}, 
$$ for $l=0,1,2, ... , L^{2}W^{2}-1$
and 
$0 \leq f_{\nn}(x) \leq 1$ for $x \in \mathbb{R}$. 

\end{slemma}

\begin{slemma}\label{lem:bspl}
Suppose that $\mathcal{B}^{d}_{\bk,\bl}$ is defined by \eqref{equ:302}.
Then, for any $L,W\geq 3$, there exists a ReLU neural network function $f_{\nn}\in \nn\left(d, \xc\label{c531} L , \xc\label{c532} W , 1\right)$ satisfying
$$
|f_{\nn}(\xx)-\mathcal{B}^{d}_{ 1,\bo}|\leq \xc\label{c533}W^{-\xc\label{c534}L},
$$
for every $\xx\in \left\{ \xx: x_{i}\in [0, r/k_{i}], 1\leq i\leq d \right\} $ and $\jc{c531}, \jc{c532}, \jc{c533}, \jc{c534}$ are all constants independent of $L$, $W$. 
\end{slemma}

Compared with Lemma 1 in \citet{suzuki2018adaptivity}, this result has more flexible choices of the depth and width, since we use the approximation construction of multivariate product given by \citet{lu2021deep}.

The following lemmas are classic results of the empirical process, and their proofs will be omitted. 
For mathematical rigor, we define 
$$\mathbb{E}\sup_{f\in \mathcal{F}}V(f)=\sup \left\{ \mathbb{E}\sup_{f \in \mathcal{F}^\ast} V(f): \mathcal{F}^\ast \subset \mathcal{F}, \mathcal{F}^\ast \ \mbox{is finite}\right\}$$ for a random process $V$ as in \citet{talagrand2014process} when the index set $\mathcal{F}$ is uncountable.

\begin{slemma}[Comparison Theorem; see Theorem 4.12 in \citet{ledoux1991probability}] \label{lem:ltc}
Let $\left\{\Phi_{i}, 1\leq i \leq n\right\}$ be functions that are $\xb\label{b103}$-Lipschitz and satisfy $\Phi_{i}(0)=0$.
Assume that $\left\{\sigma_{i}, 1\leq i \leq n\right\}$ are Rademacher random variables not depending on others. 
Then for any bounded subset $\Omega'$ in $\mathbb{R}^n$, we have
$$
\mathbb{E}\left[\sup_{ (t_1,...,t_n)\in \Omega'} \left|\frac{1}{n}\sum_{i=1}^{n} \sigma_{i} \Phi_{i}\left(t_{i}\right)\right|\right] \leq 2 \jb{b103} \cdot \mathbb{E}\left[\sup_{ (t_1,...,t_n)\in \Omega' } \left|\frac{1}{n} \sum_{i=1}^{n} \sigma_{i} t_{i}\right|\right],
$$
where the expectation is taken with respect to $\sigma_{i}$.
\end{slemma}

\begin{slemma}[Talagrand's inequality]\label{lem:tal}
Let $\left\{W_{i,f}, 1\leq i \leq n, f\in \mathcal{F}\right\}$ be  random variables independent cross $i$. Assume that for any $1\leq i \leq n$ and $f\in \mathcal{F}$, we have $|W_{i,f}|\leq \xb\label{b101}$, $\mathbb{E}\left[W_{i,f}\right]=0$
and $\mathbb{E}\left[W_{i,f}^2\right]\leq \xb\label{b102}$. 
Then with probability at least $1-\exp\left\{-t\right\}$, the following inequality holds:
$$
\sup_{f \in \mathcal{F}} \left| \frac{1}{n} \sum_{i=1}^{n} W_{i,f} \right| \leq 2\mathbb{E}\left[ \sup_{f \in \mathcal{F}} \left| \frac{1}{n} \sum_{i=1}^{n} W_{i,f} \right| \right] + \sqrt{\frac{8 \jb{b102} t}{n}}+\frac{35\jb{b101}t}{n}.
$$
\end{slemma}

Lemma \ref{lem:tal} could be derived from Theorem 4 of \citet{pascalta}. 
Unlike the usual versions of Talagrand's inequality with better constants, it does not impose the requirement that $W_{i}, 1\leq i \leq n$ to be identically distributed.

\begin{slemma}[Dudley’s entropy integral]\label{lem:chain}
Let $\mathcal{F}$ be a function class and suppose that $\sup_{f\in \mathcal{F}}\sum_{i=1}^n(f(X_{i})-\fo(X_{i}))^2/n\leq \xb^2\label{b106}$ with given $\mathcal{\x}=\{\x_{i}, 1\leq i \leq n\}$. Then we have 
$$
\mathbb{E}\left[ \sup_{f\in \mathcal{F}}
\left| \frac{1}{n}\sum_{i=1}^n \sigma_{i} (f(X_{i})-\fo(X_{i})) \right|
\right]
\leq \inf _{0 \leq t \leq \jb{b106}}\left\{4 t+\frac{12}{\sqrt{n}} \int_{t}^{\jb{b106}} \sqrt{\log \mathcal{N}\left(\varepsilon, \left.\mathcal{F}\right|_{\mathcal{X}}, L^2 \right)} d \varepsilon\right\},
$$
where the expectation is taken with respect to $\sigma_{i}$.
\end{slemma}

\begin{slemma}[Corollary 2.2 in \citet{bartlett2005local}]\label{lem:quj}
Let $\mathcal{F}$ be a function class and suppose that  $\|f\|_{\infty}\leq 1$ for any $f\in\mathcal{F}$. Suppose that positive numbers $r$ and $t$ satisfy
$$
r \geq 10  \mathbb{E}\left[ 
\sup_{f\in \mathcal{F}(r)}\left| \frac{1}{n}\sum_{i=1}^{n} \sigma_{i} \left(f(\x_{i})-\fo(\x_{i})\right) \right|
\right] +\frac{11   t}{n},
$$
where $\left\{X_{i}, 1\leq i \leq n\right\}$ are i.i.d. random variables. Then with probability at least $1-\exp\left\{-t\right\}$, one holds: 
$$
\mathcal{F}(r) \subseteq\left\{f \in \mathcal{F}: \sum_{i=1}^{n}\left(f(\x_{i})-\fo(\x_{i})\right)^2 \leq 2 r\right\}.
$$
\end{slemma}

\begin{slemma}\label{lem:lowerlinearlemma}
Consider the repeated measurement model \eqref{equ:model} and assume that $\mathcal{P}_{\x}$ has a uniform distribution on $\Omega = [0,1]^d$. 
Given a partition $\mathcal{A}$ of domain $\Omega$ with $|\mathcal{A}|\asymp 2^K$ where $K \asymp \log (nm) $, and each $A\in\mathcal{A}$ has the same measure $2^{-K}$. 
For a constant $C$, define the event $\mathcal{D}$ by
$$
\max_{A\in \mathcal{A}}\left|\left\{\x_{ij} \mid \x_{ij} \in A, 1\leq i \leq n, 1\leq j\leq m \right\}\right| \leq C 2^{-K}nm . 
$$
Consider the functional class $\mathcal{F}^{\circ}$ defined on domain $\Omega$. 
Assume that there are positive numbers $\xb\label{xb:delta}, \xb\label{xb:bf}$ and $\xb\label{xb:bcpp}$ such that the following conditions holds.  

(i) $\forall A \in \mathcal{A}$, there is a function $g \in \mathcal{F}^{\circ}$ satisfying $g(\xx) \geq \frac{1}{2}\jb{xb:delta}\jb{xb:bf}$ for any $\xx \in A$.

(ii) There is a $K'>0$ such that $ \sum_{i=1}^{n}\sum_{j=1}^{m} g(\x_{ij})^2 \leq nm\jb{xb:bcpp}\jb{xb:delta}^2 2^{-K^{\prime}}$ for any $g \in \mathcal{F}^{\circ}$ on the event $\mathcal{D}$. 
Then we have 
$$
\inf_{\tilde{f}_{L}}\sup_{\mathcal{P}_{Z,\epsilon} \in \mathcal{Q}_{Z,\epsilon}} \mathcal{E}_{\fo}(\tilde{f}_{L})  
\gtrsim \frac{1}{n} + \min\left( \frac{\jb{xb:bf}^2 2^{K'}}{\jb{xb:bcpp}nm} +  \jb{xb:delta}^2 \jb{xb:bf}^2 2^{-K} \right), 
$$
where the infimum is taken over all linear estimators $\tilde{f}_{L}$ based on the available data set $\left\{ \left(Y_{ij},\x_{ij} \right): 1\leq i \leq n, 1\leq j \leq m \right\}$ with form $
\tilde{f}_{L}(\xx)=\sum_{i=1}^n\sum_{j=1}^m Y_{ij} \varphi_{ij}\left(\xx, X_{11}, \ldots, X_{nm}\right)$ and the supremum is taken over all distributions $\mathcal{P}_{Z,\epsilon}$ satisfying the same assumptions in Theorem  \ref{thm:lowlinear}.
\end{slemma}

\section{Proofs of Theorems and Corollaries} \label{app:pftc}
In this section, we provide proofs for the theorems and corollaries presented earlier in the paper.
Due to space economy, we write  $\cf$ for the integral $\int_{\Omega}\left(f(\xx)-\fo(\xx)\right)^2d\mathcal{P}_{\x}$, $\mathcal{F}=\mathcal{F}_{nm}$, and abbreviate $\phi=\phi_{nm}$, $r^{*}=r^{*}_{nm}$, etc.

\begin{proof}[Proof of Theorem \ref{thm:thm1}]\label{prf:thm1}
The proof is divided into five steps. 
In the first step, we bound the prediction error  $\pefh$ by 
$ \hmfhf := \mathbb{E}[\sum_{i=1}^{n}\sum_{j=1}^{m}(\hat{f}(\x_{ij})-\fo(\x_{ij}))^2/nm] $. 
According to Step 2, we bound $\hmfhf$ through several parts, which will be further analyzed in Step 3 and Step 4, respectively. 
Then the desired result follows from the conclusion in Step 5. 

Step 1. 
Lemma \ref{lem:lema1} implies that
$$
\chf \leq \frac{2}{n m}\cehf+2\jc{c801}\jdb{b1}^2r^{*}
+\frac{2\jc{c802}\jdb{b1}^2t}{nm}, 
$$
holds with probability at least $1-\exp\{-t\}$. 
By taking the expectation of the above result, we can further obtain
\begin{equation}\label{equ:106}
\pefh \leq \frac{2}{n m}\mathbb{E}\left[\cehf\right]+2\jc{c801}\jdb{b1}^2r^{*}
+\frac{2\jc{c802}\jdb{b1}^2}{nm}.     
\end{equation}

Step 2. 
We now focus on the first term on the right side of \eqref{equ:106}. Recalling that $\hat{f}$ minimize the empirical loss, we have 
\begin{equation}\label{equ:107}
\hmfhf := 
\frac{1}{n m}\mathbb{E}\left[ \sum_{i=1}^{n}\sum_{j=1}^{m} \left( Y_{ij} - \hat f(\x_{ij})\right)^2   \right] \leq \frac{1}{n m}\mathbb{E}\left[ \sum_{i=1}^{n}\sum_{j=1}^{m} \left( Y_{ij} - f(\x_{ij})\right)^2   \right]
\end{equation}
holds for any $f\in \mathcal{F}$. Substituting $Y_{ij}=\fo(\x_{ij})+U_{i}(\x_{ij})+\epsilon_{ij}$ into the above inequality, a bit of algebra shows 
$$ 
\begin{aligned}
\hmfhf
\leq \int_{\Omega}\left(f-\fo\right)^2d\mathcal{P}_{\x} + \frac{2}{n m}\mathbb{E}\left[\sum_{i=1}^{n}\sum_{j=1}^{m}U_{i}(\x_{ij})\hat f(\x_{ij}) \right]+\frac{2}{n m}\mathbb{E}\left[\sum_{i=1}^{n}\sum_{j=1}^{m}\epsilon_{ij}\hat f(\x_{ij}) \right].
\end{aligned}
$$
Because $f$ is arbitrary and $\mathbb{E}\left[ U_{i}(\x_{ij}) \fo(\x_{ij}) \right]=\mathbb{E}\left[ \epsilon_{ij} \fo(\x_{ij}) \right]=0$, we have 
\begin{equation}\label{equ:112}
\begin{aligned}
\hmfhf \leq & \inf_{f\in\mathcal{F}} \int_{\Omega}\left(f-\fo\right)^2d\mathcal{P}_{\x} + \frac{2}{n m}\mathbb{E}\left[\sum_{i=1}^{n}\sum_{j=1}^{m}U_{i}(\x_{ij})\left(\hat f(\x_{ij})-\fo(\x_{ij})\right) \right]\\
&+\frac{2}{n m}\mathbb{E}\left[\sum_{i=1}^{n}\sum_{j=1}^{m}\epsilon_{ij}\left(\hat f(\x_{ij})-\fo(\x_{ij})\right)\right].
\end{aligned}   
\end{equation}

Step 3. In this step, we mainly bound $\frac{1}{n m}\mathbb{E}\left[\sum_{i=1}^{n}\sum_{j=1}^{m}U_{i}(\x_{ij})\left(\hat f(\x_{ij})-\fo(\x_{ij})\right) \right]$ in \eqref{equ:112}.
Recalling the truncation operator $\mathcal{T}_{\beta}$:  $\mathcal{T}_{\beta}f(x)=f(x)I(|f(x)|\leq\beta)+\operatorname{sign}(f)\beta I(|f(x)|>\beta)$, the triangle inequality implies
\begin{equation}\label{equ:113}
 \begin{aligned}
& \frac{1}{n m}\mathbb{E}\left[\sum_{i=1}^{n}\sum_{j=1}^{m}U_{i}(\x_{ij})\left(\hat f(\x_{ij})-\fo(\x_{ij})\right) \right]\\
\leq & \mathbb{E}\left[ \left|\frac{1}{n} \sum_{i=1}^{n}\left(\frac{1}{m} \sum_{j=1}^{m} \left(\mathcal{T}_{\beta_{U}} U_{i}\right)\left(\x_{i j}\right)\left(\hat{f}\left(\x_{i j}\right)-\fo\left(\x_{i j}\right)\right)-\int_{\Omega} \left(\mathcal{T}_{\beta_{U}} U_{i}\right)\left(\hat{f}-\fo\right) d \mathcal{P}_{\x}\right)\right| \right]\\ 
&+\mathbb{E}\left[ \left|\frac{1}{nm} \sum_{i=1}^{n} \sum_{j=1}^{m} \left(U_{i} - \mathcal{T}_{\beta_{U}} U_{i}\right)\left(\x_{i j}\right)\left(\hat{f}\left(\x_{i j}\right)-\fo\left(\x_{i j}\right)\right)-\int_{\Omega} \left(U_{i} - \mathcal{T}_{\beta_{U}} U_{i}\right)\left(\hat{f}-\fo\right) d \mathcal{P}_{\x}\right| \right]\\ 
&+\mathbb{E}\left[\left|\frac{1}{n} \sum_{i=1}^{n} \int_{\Omega} U_{i}\left(\hat{f}-\fo\right) d \mathcal{P}_{\x}\right|\right]\\
:=& \left(\uppercase\expandafter{\romannumeral1}\right)+\left(\uppercase\expandafter{\romannumeral2}\right)+\left(\uppercase\expandafter{\romannumeral3}\right).
\end{aligned}   
\end{equation}

Step 3.1. 
To control $\left(\uppercase\expandafter{\romannumeral1}\right)$, we first study 
$$
\sup_{f\in \mathcal{F}}
\left| 
\frac{1}{n m} \sum_{i=1}^{n} \sum_{j=1}^{m} \frac{\left(\mathcal{T}_{\beta_{U}} U_{i}\right)\left(\x_{i j}\right)\left(f\left(\x_{i j}\right)-\fo\left(\x_{i j}\right)\right)-\int_{\Omega} \left(\mathcal{T}_{\beta_{U}} U_{i}\right)(f-\fo) d \mathcal{P}_{\x}}{\cf  +r} 
\right|.
$$
By the symmetrization result of the Rademacher process, we have
$$
\begin{aligned}
&\mathbb{E}\left[ \sup_{f\in \mathcal{F}} \left| 
\frac{1}{n m} \sum_{i=1}^{n} \sum_{j=1}^{m} \frac{\left(\mathcal{T}_{\beta_{U}} U_{i}\right)\left(\x_{i j}\right)\left(f\left(\x_{i j}\right)-\fo\left(\x_{i j}\right)\right)-\int_{\Omega} \left(\mathcal{T}_{\beta_{U}} U_{i}\right)(f-\fo) d \mathcal{P}_{\x}}{\cf  +r} 
\right| \right]\\
\leq& 2 \mathbb{E}\left[\sup _{f \in \mathcal{F}}\left|\frac{1}{n m} \sum_{i=1}^{n} \sum_{j=1}^{m} \frac{\sigma_{i j} \left(\mathcal{T}_{\beta_{U}} U_{i}\right)\left(\x_{i j}\right)\left(f\left(\x_{i j}\right)-\fo\left(\x_{i j}\right)\right)}{\cf+r}\right|\right],
\end{aligned}
$$
where $\sigma_{ij}$'s are i.i.d. Rademacher random variables. 
Then a peeling of the function class implies that it is less than 
$$
\begin{aligned}
\leq &  2\mathbb{E}\left[\sup _{f \in \mathcal{F}(r)}\left|\frac{1}{n m} \sum_{i=1}^{n} \sum_{j=1}^{m} \frac{\sigma_{i j} \left(\mathcal{T}_{\beta_{U}} U_{i}\right)\left(\x_{i j}\right)\left(f\left(\x_{i j}\right)-\fo\left(\x_{i j}\right)\right)}{\cf  +r}\right|\right]\\
&+\sum_{k=1}^{\infty}2\mathbb{E}\left[\sup _{f \in \mathcal{F}(4^{k}r) \backslash \mathcal{F}(4^{k-1}r)}\left|\frac{1}{n m} \sum_{i=1}^{n} \sum_{j=1}^{m} \frac{\sigma_{i j} \left(\mathcal{T}_{\beta_{U}} U_{i}\right)\left(\x_{i j}\right)\left(f\left(\x_{i j}\right)-\fo\left(\x_{i j}\right)\right)}{\cf  +r}\right|\right], 
\end{aligned}
$$
which could be further bounded by 
\begin{equation}\label{equ:180}
\begin{aligned}
& \frac{2}{r}\mathbb{E}\left[\sup _{f \in \mathcal{F}(r)}\bigg|\frac{1}{n m} \sum_{i=1}^{n} \sum_{j=1}^{m} \sigma_{i j} \left(\mathcal{T}_{\beta_{U}} U_{i}\right)\left(\x_{i j}\right)\left(f\left(\x_{i j}\right)-\fo\left(\x_{i j}\right)\right)\bigg|\right]\\
&+\sum_{k=1}^{\infty} \frac{2}{(4^{k-1}+1)r} \mathbb{E}\left[\sup _{f \in \mathcal{F}(4^{k}r) }\bigg|\frac{1}{n m} \sum_{i=1}^{n} \sum_{j=1}^{m} \sigma_{i j} \left(\mathcal{T}_{\beta_{U}} U_{i}\right)\left(\x_{i j}\right)\left(f\left(\x_{i j}\right)-\fo\left(\x_{i j}\right)\right) \bigg|\right].\\
\end{aligned}
\end{equation}
Taking the expectation conditional on ${U_i, X_{ij} , 1 \leq i \leq n, 1 \leq j \leq m}$ and utilizing Lemma \ref{lem:ltc}, we conclude that
$$
\begin{aligned}
&\mathbb{E}\left[ \sup _{f \in \mathcal{F}(r')}\bigg|\frac{1}{n m} \sum_{i=1}^{n} \sum_{j=1}^{m} \sigma_{i j} \left(\mathcal{T}_{\beta_{U}} U_{i}\right)\left(\x_{i j}\right)\left(f\left(\x_{i j}\right)-\fo\left(\x_{i j}\right)\right)\bigg| \right]\\
=& \mathbb{E}\left[\mathbb{E}\left[ \sup _{f \in \mathcal{F}(r')}\left.\bigg|\frac{1}{n m} \sum_{i=1}^{n} \sum_{j=1}^{m} \sigma_{i j} \left(\mathcal{T}_{\beta_{U}} U_{i}\right)\left(\x_{i j}\right)\left(f\left(\x_{i j}\right)-\fo\left(\x_{i j}\right)\right)\bigg|\right|U,X \right]\right]\\
\leq &  2\beta_{U} \mathbb{E}\left[ \sup _{f \in \mathcal{F}(r')}\bigg|\frac{1}{n m} \sum_{i=1}^{n} \sum_{j=1}^{m} \sigma_{i j} \left(f\left(\x_{i j}\right)-\fo\left(\x_{i j}\right)\right)\bigg|\right]\\
\leq & 2\beta_{U} \phi(r'),
\end{aligned}
$$
for any $r'\geq r^{*}$.
Combining the previous several inequalities, we thus obtain
\begin{align}\label{equ:120}
&\mathbb{E}\left[ \sup_{f\in \mathcal{F}} \left| 
\frac{1}{n m} \sum_{i=1}^{n} \sum_{j=1}^{m} \frac{\left(\mathcal{T}_{\beta_{U}} U_{i}\right)\left(\x_{i j}\right)\left(f-\fo\right)\left(\x_{i j}\right)-\int_{\Omega} \left(\mathcal{T}_{\beta_{U}} U_{i}\right)(f-\fo) d \mathcal{P}_{\x}}{\cf  +r} 
\right| \right] \notag \\
 &\leq \frac{4\beta_{U} \phi(r)}{r}+\sum_{k=1}^{\infty} \frac{4\beta_{U} \phi(4^k r)}{(4^{k-1}+1)r} 
\leq  \frac{4\beta_{U} \phi(r)}{r}+\sum_{k=1}^{\infty} \frac{2^{k+2} \beta_{U} \phi( r)}{(4^{k-1}+1)r} 
\leq  \frac{16\beta_{U} \phi(r)}{r}.
\end{align}
Furthermore, one can verify that
\begin{equation}\label{equ:121}
\mathbb{E}\left[
\frac{\left(\mathcal{T}_{\beta_{U}} U_{i}\right)\left(\x\right)\left(f\left(\x\right)-\fo\left(\x\right)\right)-\int_{\Omega} \left(\mathcal{T}_{\beta_{U}} U_{i}\right)(f-\fo) d \mathcal{P}_{\x}}{\cf  +r} 
\right]=0,
\end{equation}
\begin{equation}\label{equ:122}
\left| 
\frac{\left(\mathcal{T}_{\beta_{U}} U_{i}\right)\left(\x\right)\left(f\left(\x\right)-\fo\left(\x\right)\right)-\int_{\Omega} \left(\mathcal{T}_{\beta_{U}} U_{i}\right)(f-\fo) d \mathcal{P}_{\x}}{\cf  +r} 
\right|\leq \frac{4\jdb{b1}\beta_{U}}{r},
\end{equation}
and
\begin{equation}\label{equ:123}
\mathbb{E}\left[\left| 
\frac{\left(\mathcal{T}_{\beta_{U}} U_{i}\right)\left(\x\right)\left(f\left(\x\right)-\fo\left(\x\right)\right)-\int_{\Omega} \left(\mathcal{T}_{\beta_{U}} U_{i}\right)(f-\fo) d \mathcal{P}_{\x}}{\cf  +r} 
\right|^2\right]\leq \frac{\beta_{U}^2}{4r}.
\end{equation}
By combining with \eqref{equ:120}, \eqref{equ:121}, \eqref{equ:122} and \eqref{equ:123} and  checking the conditions  required by Lemma \ref{lem:tal}, we can conclude that, with probability at least $1-\exp\{-t\}$,
$$
\begin{aligned}
& \left| 
\frac{1}{n m} \sum_{i=1}^{n} \sum_{j=1}^{m} \frac{\left(\mathcal{T}_{\beta_{U}} U_{i}\right)\left(\x_{i j}\right)\left(\hat f\left(\x_{i j}\right)-\fo\left(\x_{i j}\right)\right)-\int_{\Omega} \left(\mathcal{T}_{\beta_{U}} U_{i}\right)\left(\hat{f}-\fo\right) d \mathcal{P}_{\x}}{\chf  +r} 
\right|\\
\leq & \sup_{f\in \mathcal{F}}
\left| 
\frac{1}{n m} \sum_{i=1}^{n} \sum_{j=1}^{m} \frac{\left(\mathcal{T}_{\beta_{U}} U_{i}\right)\left(\x_{i j}\right)\left(f\left(\x_{i j}\right)-\fo\left(\x_{i j}\right)\right)-\int_{\Omega} \left(\mathcal{T}_{\beta_{U}} U_{i}\right)(f-\fo) d \mathcal{P}_{\x}}{\cf  +r} 
\right|\\
\leq & \frac{32\beta_{U} \phi(r)}{r}+\sqrt{ \frac{2\beta_{U}^2 t}{nmr} }+ \frac{140\jdb{b1}\beta_{U} t}{nmr}.
\end{aligned}
$$
We determine universal constants $\xc\label{c107}$ and $\xc\label{c108}$ such that $r_{U}(t)=\jc{c107}\beta_{U}^2r^{*}+\frac{\jc{c108}(\jdb{b1}+\beta_{U})\beta_{U} t }{nm}$ satisfies 
$$
\frac{32\beta_{U} \phi(r_{U})}{r_{U}}+\sqrt{ \frac{2\beta_{U}^2 t}{nmr_{U}} }+ \frac{140\jdb{b1}\beta_{U} t}{nmr_{U}} \leq \frac{1}{32}.
$$
Therefore, with probability at least $1-\exp\{-t\}$,
$$
\begin{aligned}
& \left| 
\frac{1}{n m} \sum_{i=1}^{n} \sum_{j=1}^{m} \left(\mathcal{T}_{\beta_{U}} U_{i}\right)\left(\x_{i j}\right)\left(\hat f\left(\x_{i j}\right)-\fo\left(\x_{i j}\right)\right)-\int_{\Omega} \left(\mathcal{T}_{\beta_{U}} U_{i}\right)\left(\hat{f}-\fo\right) d \mathcal{P}_{\x}
\right|\\
\leq &\frac{1}{32}\chf  +\frac{\jc{c107}\beta_{U}^2r^{*}}{32}+\frac{\jc{c108}(\jdb{b1}+\beta_{U})\beta_{U} t }{32nm}.
\end{aligned}
$$
Furthermore, we can bound the expectation
\begin{equation}\label{equ:125}
\begin{aligned}
\left(\uppercase\expandafter{\romannumeral1}\right)
\leq \frac{1}{32}\pefh  
+\jc{c109}\beta_{U}^2r^{*}+\frac{\jc{c110}(\jdb{b1}+\beta_{U})\beta_{U}  }{nm}, 
\end{aligned}    
\end{equation}
where $\xc\label{c109}$ and $\xc\label{c110}$ are universal constants.

Step 3.2. 
In this substep, we focus on $\left(\uppercase\expandafter{\romannumeral2}\right)$. Applying the absolute value inequality, we obtain
$$
\begin{aligned}
\left(\uppercase\expandafter{\romannumeral2}\right)\leq & 2\mathbb{E}\left[ \left|\left(U_{1} - \mathcal{T}_{\beta_{U}} U_{1}\right)\left(\x_{11}\right)\left(\hat{f}\left(\x_{i j}\right)-\fo\left(\x_{i j}\right)\right)\right| \right]\\
\leq &4\jdb{b1}\mathbb{E}\left[ \left|U_{11}\right|I(\left|U_{11}\right|>\beta_{U} ) \right].
\end{aligned}
$$
Applying the facts
\begin{equation}\label{equ:140}
\left|U_{11}\right|\leq 2\jdb{b2} \exp\left\{\frac{|U_{11}|}{2\jdb{b2} }\right\} 
\ \text{ and }\ 
I(\left|U_{11}\right|>\beta_{U} ) \leq \exp\left\{\frac{|U_{11}|-\beta_{U}}{2\jdb{b2} }\right\},
\end{equation}
we have 
\begin{equation}\label{equ:126}
\left(\uppercase\expandafter{\romannumeral2}\right)\leq 8\jdb{b1} \jdb{b2} \mathbb{E}\left[\exp\left\{\frac{|U_{11}|}{\jdb{b2} } -  \frac{\beta_{U}}{2\jdb{b2} } \right\}\right] = 8\jdb{b1} \jdb{b2} \exp\left\{- \frac{\beta_{U}}{2\jdb{b2} } \right\}.
\end{equation}

Step 3.3. Now we consider $\left(\uppercase\expandafter{\romannumeral3}\right)$. 
 Using the GM-QM inequality,
 $$
\begin{aligned}
\left(\uppercase\expandafter{\romannumeral3}\right)
= &\mathbb{E}\left[\left| \int_{\Omega}\left( \frac{1}{n}\sum_{i=1}^{n}U_{i}\right) \left(\hat{f}-\fo\right) d \mathcal{P}_{\x}\right|\right]\\
\leq & 8 \mathbb{E}\left[ \int_{\Omega}\left( \frac{1}{n}\sum_{i=1}^{n}U_{i}\right)^2 d \mathcal{P}_{\x}\right]+\frac{1}{32}\mathbb{E}\left[ \int_{\Omega} \left(\hat{f}-\fo\right)^ 2 \mathcal{P}_{\x}\right].
\end{aligned}
$$
Because $\mathbb{E}\left[ U(X)^2 \right]\leq \mathbb{E}\left[\jdb{b2}^2 \exp\left\{U(X)/\jdb{b2} \right\}\right] \leq \jdb{b2}^2 $, we have $\mathbb{E}\left[ \int_{\Omega}\left( \sum_{i=1}^{n}U_{i}/n\right)^2 d \mathcal{P}_{\x}\right]\leq \jdb{b2}^2 /n$. Substituting this into the earlier inequality yields
\begin{equation}\label{equ:127}
\begin{aligned}
\left(\uppercase\expandafter{\romannumeral3}\right)
\leq  \frac{8\jdb{b2}^2 }{n} + \frac{1}{32}\pefh.
\end{aligned}    
\end{equation}

Substituting \eqref{equ:125}, \eqref{equ:126} and \eqref{equ:127} into \eqref{equ:113} and setting $\beta_{U}=2\jdb{b2} \log n$, we have 
\begin{equation}\label{equ:128}
\begin{aligned}
&\frac{1}{n m}\mathbb{E}\left[\sum_{i=1}^{n}\sum_{j=1}^{m}U_{i}(\x_{ij})\left(\hat f(\x_{ij})-\fo(\x_{ij})\right) \right] \\
\leq & \frac{1}{16}\pefh+ \frac{\xc (\jdb{b1}^2+\jdb{b2}^2 )}{n}+ \frac{5\jc{c110}(\jdb{b1}^2+\jdb{b2}^2 \log^2 n)}{nm}+4\jc{c109}\jdb{b2}^2 r^{*}\log^2 n.
\end{aligned}
\end{equation}

Step 4. 
In this step, we aim to control $\frac{1}{n m}\mathbb{E}\left[\sum_{i=1}^{n}\sum_{j=1}^{m}\epsilon_{ij}\left(\hat f(\x_{ij})-\fo(\x_{ij})\right)\right]$. 
The proof in this step is similar to Step 3. The main gap is that $\epsilon_{ix}$ is not integrable with respect to $\xx$. 
Thus it cannot be given by a direct extension of Step 3. 
Using the triangle inequality, we have
\begin{align}\label{equ:129}
& \frac{1}{n m}\mathbb{E}\left[\sum_{i=1}^{n}\sum_{j=1}^{m}\epsilon_{ij}\left(\hat f(\x_{ij})-\fo(\x_{ij})\right) \right] \notag \\ 
\leq & \mathbb{E}\left[ \left|\frac{1}{n} \sum_{i=1}^{n}\left(\frac{1}{m} \sum_{j=1}^{m} \left(\mathcal{T}_{\beta_{\epsilon}} \epsilon_{ij} - \mathbb{E}\left[ \mathcal{T}_{\beta_{\epsilon}} \epsilon_{ij} \right]\right)\left(\hat{f}\left(\x_{i j}\right)-\fo\left(\x_{i j}\right)\right)
\right)\right| \right]\\ 
&+\left|\mathbb{E}\left[ \frac{1}{nm} \sum_{i=1}^{n}\sum_{j=1}^{m}  \left(\epsilon_{ij} - \mathcal{T}_{\beta_{\epsilon}} \epsilon_{ij}+\mathbb{E}\left[ \mathcal{T}_{\beta_{\epsilon}} \epsilon_{ij} \right]\right)\left(\hat{f}\left(\x_{i j}\right)-\fo\left(\x_{i j}\right)\right) \right]\right| \notag \\ 
:=& \left(\uppercase\expandafter{\romannumeral4}\right)+\left(\uppercase\expandafter{\romannumeral5}\right). \notag
\end{align}

Step 4.1.
Lemma \ref{lem:lema2} implies that, 
with probability at least $1-\exp\{-t\}$, 
    \begin{align*}
    & \left| 
     \frac{1}{n m} \sum_{i=1}^{n} \sum_{j=1}^{m} \left(\mathcal{T}_{\beta_{\epsilon}} \epsilon_{ij}-\mathbb{E}\left[ \mathcal{T}_{\beta_{\epsilon}} \epsilon_{ij} \right]\right)\left(\hat{f}\left(\x_{i j}\right)-\fo\left(\x_{i j}\right)\right)
    \right|\\
     \leq &\frac{1}{16}\chf  +\jc{c803}\beta_{\epsilon}^2r^{*}+\frac{\jc{c804}(\jdb{b1}\beta_{\epsilon}+ \jdb{b3}^{2}  ) t }{nm}.
     \end{align*}
Further, by taking the expectation of the previous result, we have
\begin{equation}\label{equ:135}
\begin{aligned}
\left(\uppercase\expandafter{\romannumeral4}\right) \leq  \frac{1}{16}\pefh  +\jc{c803}\beta_{\epsilon}^2r^{*}+\frac{\jc{c804}(\jdb{b1}\beta_{\epsilon}+ \jdb{b3}^{2}  )  }{nm}.
\end{aligned}    
\end{equation}

Step 4.2.
Next we bound $\left(\uppercase\expandafter{\romannumeral5}\right)$. Note that $\mathbb{E}\left[\epsilon_{ij}\right]=0$ yields $\left|\mathbb{E}\left[ \mathcal{T}_{\beta_{\epsilon}} \epsilon_{ij} \right]\right|=\left|\mathbb{E}\left[ \epsilon_{ij} I(|\epsilon_{ij}|>\beta_{\epsilon}) \right]\right|$. Then 
\begin{equation}\label{equ:136}
 \begin{aligned}
\left(\uppercase\expandafter{\romannumeral5}\right)
\leq & 2\jdb{b1}\mathbb{E}\left[ \left| \epsilon_{ij} - \mathcal{T}_{\beta_{\epsilon}} \epsilon_{ij}+\mathbb{E}\left[ \mathcal{T}_{\beta_{\epsilon}} \epsilon_{ij} \right] \right| \right]
\leq  4\jdb{b1} \mathbb{E}\left[\left| \epsilon_{11}\right| I(|\epsilon_{11}|>\beta_{\epsilon})  \right]\\
\leq & 4\jdb{b1} \mathbb{E}\left[ 2\jdb{b3} \exp\left\{\frac{|\epsilon_{11}|}{2\jdb{b3} }\right\}\exp\left\{\frac{|\epsilon_{11}|-\beta_{\epsilon}}{2\jdb{b3}}\right\}  \right]
\leq  8\jdb{b1}  \jdb{b3}  \exp\left\{\frac{-\beta_{\epsilon}}{2\jdb{b3} }\right\}.
\end{aligned}   
\end{equation}
Substituting \eqref{equ:135} and \eqref{equ:136} into \eqref{equ:129} and setting $\beta_{\epsilon}=2\jdb{b3} \log nm$, we have 
\begin{equation}\label{equ:137}
\begin{aligned}
&\frac{1}{n m}\mathbb{E}\left[\sum_{i=1}^{n}\sum_{j=1}^{m}\epsilon_{ij}\left(\hat f(\x_{ij})-\fo(\x_{ij})\right) \right] \\
\leq & \frac{1}{16}\pefh+ \frac{\jc{c118}(\jdb{b1}^2+\jdb{b3}^{2} \log ^2 nm)}{nm}+\jc{c119}\jdb{b3}^{2} r^{*}\log^2 nm,
\end{aligned}
\end{equation}
where $\xc\label{c118}$ and $\xc\label{c119}$ are universal constants. 

Step 5. 
Plugging \eqref{equ:128}\eqref{equ:137} into \eqref{equ:112} yields
\begin{equation}\label{equ:138}
\begin{aligned}
\hmfhf
\leq & \inf_{f\in\mathcal{F}} \int_{\Omega}\left(f-\fo\right)^2d\mathcal{P}_{\x} 
+\frac{1}{4}\pefh+\frac{\xc(\jdb{b1}^2+\jdb{b2}^2 )}{n}
\\&+ \frac{\xc(\jdb{b1}^2+\jdb{b2}^2 \log n+\jdb{b3}^{2} \log nm)}{nm}
+\xc r^{*}(\jdb{b2}^2 \log n+\jdb{b3}^{2} \log nm)).\\
\end{aligned}   
\end{equation}
Substituting \eqref{equ:138} into \eqref{equ:106} and noting that all constants $c_{k}$ used in this proof are universal, the desired result follows. 
Thus we complete the proof.
\end{proof}

\begin{remark} \label{rem:rem1}

It is known that, in practice, the exact empirical risk minimizer may not be attainable due to optimization error, while the non-convex function class could also trap the gradient-descent-based optimization methods into local minima. Thus 
it is meaningful to study the theoretical performance of estimators when they cannot strictly reach the global minimum. 
Like in \citet{schmidtaos}, we define the optimization error of an estimator $\hat{f}_{nm}$ by
$$
\delf
:=\mathbb{E}\left[\frac{1}{nm} \sum_{i=1}^{n}\sum_{j=1}^{m}\left(Y_{ij}-\hat{f}_{nm}\left(\x_{ij}\right)\right)^{2}-\inf _{f \in \mathcal{F}_{nm}} \frac{1}{n} \sum_{i=1}^{n}\sum_{j=1}^{m}\left(Y_{ij}-f\left(\x_{ij}\right)\right)^{2}\right],
$$
which measures the difference between the estimator $\hat{f}_{nm}$ and the global empirical risk minimum over the chosen class $\mathcal{F}_{nm}$.
It is easy to show that $\delf \geq 0$ and $\delf = 0$ if and only if $\hat{f}_{nm}$ satisfies \eqref{equ:estimator}. Consider the optimization error $\delf$ defined above and 
under the same model and assumptions as Theorem \ref{thm:thm1}, there are 
universal positive constants $\xc\label{c005}$, $\xc\label{c006}$ and $\xc\label{c007}$ such that 
$$
\begin{aligned}
&\jc{c005}\delf - \jc{c006}\left(\frac{\jdb{b1}^2+\jdb{b2}^2}{n} + (r^{*}_{nm} + \frac{1}{nm})\iota_{nm}
\right)\\
&\leq \pef \\
&\leq 
\jc{c007}\left(\inf_{f\in\mathcal{F}_{nm}} \int_{\Omega}\left(f-\fo\right)^2d\mathcal{P}_{\x} + \delf 
+\frac{\jdb{b1}^2+\jdb{b2}^2}{n}+(r^{*}_{nm} + \frac{1}{nm})\iota_{nm}
\right),\\
\end{aligned}
$$
where $\iota_{nm}=\jdb{b1}^2+\jdb{b2}^2(\log n)^2+\jdb{b3}^2(\log nm)^2$.

We first prove the upper bound. 
It requires a reconstruction of Step 2 in the proof of Theorem \ref{thm:thm1}. 
Obviously, 
$$
\inf _{f \in \mathcal{F}_{nm}} \frac{1}{n} \sum_{i=1}^{n}\sum_{j=1}^{m}\left(Y_{ij}-f\left(\x_{ij}\right)\right)^{2}
\leq \frac{1}{n m} \sum_{i=1}^{n}\sum_{j=1}^{m} \left( Y_{ij} - f(\x_{ij})\right)^2 
$$
holds for any $f\in \mathcal{F}$. 
Therefore we have
$$
\frac{1}{n m}\mathbb{E}\left[ \sum_{i=1}^{n}\sum_{j=1}^{m} \left( Y_{ij} - \hat f(\x_{ij})\right)^2   \right] \leq \frac{1}{n m}\mathbb{E}\left[ \sum_{i=1}^{n}\sum_{j=1}^{m} \left( Y_{ij} - f(\x_{ij})\right)^2   \right] + \delf. 
$$
Because $f$ is arbitrary and $\mathbb{E}\left[ U_{i}(\x_{ij})\hat \fo(\x_{ij}) \right]=\mathbb{E}\left[ \epsilon_{ij}\hat \fo(\x_{ij}) \right]=0$, we have 
$$
\begin{aligned}
& \frac{1}{n m}\mathbb{E}\left[\cehf\right]\\
\leq & \inf_{f\in\mathcal{F}} \int_{\Omega}\left(f-\fo\right)^2d\mathcal{P}_{\x} + \frac{2}{n m}\mathbb{E}\left[\sum_{i=1}^{n}\sum_{j=1}^{m}U_{i}(\x_{ij})\left(\hat f(\x_{ij})-\fo(\x_{ij})\right) \right]\\
&+\frac{2}{n m}\mathbb{E}\left[\sum_{i=1}^{n}\sum_{j=1}^{m}\epsilon_{ij}\left(\hat f(\x_{ij})-\fo(\x_{ij})\right)\right] + \delf.
\end{aligned}   
$$
The rest is the same as the proof of Theorem \ref{thm:thm1}.

For the lower bound.  
Let    $$\tilde{f}_{nm}\in \arg\min_{f\in \mathcal{F}_{nm}}\frac{1}{nm}\sum_{i=1}^{n}\sum_{j=1}^{m}\left(Y_{ij}-f(X_{ij})\right)^{2}$$ be a exact empirical risk minimizer.
A few algebra shows 
$$
\begin{aligned}
& \frac{1}{n m}\mathbb{E}\left[\cehf\right]-  \frac{1}{n m}\mathbb{E}\left[\cetf\right]\\
= & \delf 
+ \frac{2}{n m}\mathbb{E}\left[\sum_{i=1}^{n}\sum_{j=1}^{m}(U_{i}(\x_{ij})+\epsilon_{ij})(\hat f(\x_{ij})-\fo(\x_{ij})) \right]\\
&-\frac{2}{n m}\mathbb{E}\left[\sum_{i=1}^{n}\sum_{j=1}^{m}(U_{i}(\x_{ij})+\epsilon_{ij})(\hat f(\x_{ij})-\fo(\x_{ij})) \right].
\end{aligned}
$$
Using the results of Step 3 and Step 4 in the Proof of Theorem \ref{thm:thm1}, we have 
\begin{equation}\label{equ:151}
    \begin{aligned}
& \frac{1}{n m}\mathbb{E}\left[\cehf\right]-  \frac{1}{n m}\mathbb{E}\left[\cetf\right]\\
\geq & \delf 
- \frac{1}{4}\left(  \mathbb{E}\left[\chf\right] +  \mathbb{E}\left[\ctf\right]  \right)\\
&- \xc\left( \frac{\jdb{b1}^2+\jdb{b2}}{n} +  (r^{*}_{nm} + \frac{1}{nm})\iota_{nm} \right)
\end{aligned}.
\end{equation}
Similar to the proof of Lemma \ref{lem:lema1} and Step 1 in the proof of Theorem \ref{thm:thm1}, we can derive that
$$
\mathbb{E}\left[\ctf\right] \leq \frac{2}{n m}\mathbb{E}\left[\cetf\right]+2\jc{c801}\jdb{b1}^2r^{*}
+\frac{2\jc{c802}\jdb{b1}^2}{nm}. 
$$
and 
$$
\mathbb{E}\left[\ctf\right] \geq \frac{2}{3n m}\mathbb{E}\left[\cetf\right]-\frac{2\jc{c801}\jdb{b1}^2r^{*}}{3}
-\frac{2\jc{c802}\jdb{b1}^2}{3nm}.   
$$
Substituting the above two formulas and their counterparts about $\hat{f}$ together into the \eqref{equ:151}. 
Noting that the term $\mathbb{E}\left[\cetf\right]$ is nonnegative, the desired lower bound follows. 
    
\end{remark}

\begin{proof}[Proof of Corollary \ref{cor:cor1} ]

(i) Let  $\mathcal{G}=\left\{g=f/\jdb{b1}:f\in \mathcal{F}\right\}$, then every $g\in \mathcal{G}$ satisfies $\|g\|_{\infty}\leq 1$. Denote $\go=\fo/\jdb{b1}$. Then the estimation problem becomes 
$$
\hat{g}_{nm}=\arg\min_{g\in \mathcal{G}_{nm}}\sum_{i=1}^{n}\sum_{j=1}^{m}\left(Y_{ij}/\jdb{b1}-g(X_{ij})\right)^{2}.
$$
Denote $\jdb{b2} '=\jdb{b2} /\jdb{b1}$ and $\jdb{b3} '=\jdb{b3} /\jdb{b1}$. By Theorem \ref{thm:thm1}, we have 
\begin{equation}\label{equ:145}
\begin{aligned}
  \mathcal{E}_{\go}(\hat{g}_{nm})
\leq & 
\jc{c001}\left(\inf_{g\in\mathcal{G}} \int_{\Omega}\left(g-\go\right)^2d\mathcal{P}_{\x} 
+\frac{1+\jdb{b2}'^{2} }{n}\right.
\\ & \quad \quad + \left.\left(r^{*}_{nm}(\mathcal{G}) + \frac{1}{nm}\right)(1+\jdb{b2}'^2 \log n+\jdb{b3}'^{2} \log nm)\right).\\
\end{aligned}
\end{equation}
Let
$$
r_1=\sup\left\{ r: 10\mathbb{E}\left[ 
\sup_{g\in \mathcal{G}(r)}\left| \frac{1}{nm}\sum_{i=1}^{n}\sum_{j=1}^{m} \sigma_{ij} \left(g(\x_{ij})-\go(\x_{ij})\right) \right|
\right] +\frac{11 \log nm}{nm} \geq r \right\}.
$$
Because $
\mathbb{E}\left[ 
\sup_{g\in \mathcal{G}(r)}\left| \frac{1}{nm}\sum_{i=1}^{n}\sum_{j=1}^{m} \sigma_{ij} \left(g(\x_{ij})-\go(\x_{ij})\right) \right|
\right]
$ is monotonic non-decreasing with respect to $r$, simple analysis shows that $r_1$ satisfies
\begin{equation}\label{equ:165}
r_1=10\mathbb{E}\left[ 
\sup_{g\in \mathcal{G}(r_1)}\left| \frac{1}{nm}\sum_{i=1}^{n}\sum_{j=1}^{m} \sigma_{ij} \left(g(\x_{ij})-\go(\x_{ij})\right) \right|
\right] +\frac{11 \log nm}{nm}.
\end{equation}
(In fact, we can deduce a contradiction if replacing the equal sign in the above equation with either the greater or less.)
 For any $r\geq r_1$, consider
$$
\begin{aligned}
\mathcal{R}_{nm}\left(\mathcal{G}(r)\right)=&\mathbb{E}\left[\sup_{g\in \mathcal{G}(r)}\left| \frac{1}{nm}\sum_{i=1}^{n}\sum_{j=1}^{m} \sigma_{ij} \left(g(\x_{ij})-\go(\x_{ij})\right)\right| \right]\\
\leq &  \mathbb{E}\left[\sup_{g\in \hat{\mathcal{G}}(2r)}\left| \frac{1}{nm}\sum_{i=1}^{n}\sum_{j=1}^{m} \sigma_{ij} \left(g(\x_{ij})-\go(\x_{ij})\right)\right|\right]\\
&+ \mathbb{E}\left[\sup_{g\in \mathcal{G}(r)\backslash \hat{\mathcal{G}}(2r)}\left| \frac{1}{nm}\sum_{i=1}^{n}\sum_{j=1}^{m} \sigma_{ij} \left(g(\x_{ij})-\go(\x_{ij})\right)\right| \right].
\end{aligned}
$$
where $ \hat{\mathcal{G}}(2r)=\left\{ g\in\mathcal{G}: \sum_{i=1}^{n}\sum_{j=1}^{m}\left( \hat g(\x_{ij})-\go(\x_{ij})\right)^2/nm \leq 2r \right\}$ is a random set depending on samples.
 Applying Lemma \ref{lem:chain}, we have  
$$
\begin{aligned}
\sup_{g\in \hat{\mathcal{G}}(2r)}\left| \frac{1}{nm}\sum_{i=1}^{n}\sum_{j=1}^{m} \sigma_{ij} \left(g(\x_{ij})-\go(\x_{ij})\right)\right| \leq& \frac{4}{nm}+\frac{12}{\sqrt{nm}} \int_{1/nm}^{\sqrt{2r}} \sqrt{\log \mathcal{N}\left(\varepsilon, \left.\mathcal{G}\right|_{\mathcal{X}}, L^2 \right)} d \varepsilon \\
\leq & \frac{4}{nm}+\frac{12\sqrt{2r}}{\sqrt{nm}}\sqrt{\log \mathcal{N}\left(1/nm, \left.\mathcal{G}\right|_{\mathcal{X}}, L^2 \right)}.
\end{aligned}
$$
Hence, the first term is bounded by 
$$
\frac{4}{nm}+\frac{12\sqrt{2r}}{\sqrt{nm}}\mathbb{E}\left[\sqrt{\log \mathcal{N}\left(1/nm, \left.\mathcal{G}\right|_{\mathcal{X}}, L^2 \right)}\right].
$$
Also, the second term could be bounded by
$$
\begin{aligned}
& \mathbb{E}\left[\sup_{g\in \mathcal{G}(r)\backslash \hat{\mathcal{G}}(2r)}\left| \frac{1}{nm}\sum_{i=1}^{n}\sum_{j=1}^{m} \sigma_{ij} \left(g(\x_{ij})-\go(\x_{ij})\right)\right| \right]\\
\leq & \mathbb{P}\left(\left\{ \exists g\in\mathcal{G}(r):\  \frac{1}{nm}\sum_{i=1}^{n}\sum_{j=1}^{m} \sigma_{ij} \left(g(\x_{ij})-\go(\x_{ij})\right) > 2r \right\}\right)\cdot 2 \\
\leq & \frac{2 }{nm}, 
\end{aligned}
$$
where the last inequality is given by Lemma \ref{lem:quj} and \eqref{equ:165}. 
Hence,
$$
\begin{aligned}
\mathcal{R}_{nm}\left(\mathcal{G}(r)\right)&\leq \frac{6}{nm}+\frac{12\sqrt{2r}}{\sqrt{nm}}\mathbb{E}\left[\sqrt{\mathcal{N}\left(1/nm, \left.\mathcal{G}\right|_{\mathcal{X}}, L^2 \right)}\right]\\
\end{aligned}
$$
is bounded by a sub-root function for $r\geq r_1$. In addition, recalling \eqref{equ:165}, we have
$$
\begin{aligned}
r_1 =& 10\mathbb{E}\left[ 
\sup_{g\in \mathcal{G}(r_1)}\left| \frac{1}{nm}\sum_{i=1}^{n}\sum_{j=1}^{m} \sigma_{ij} \left(g(\x_{ij})-\go(\x_{ij})\right) \right|
\right] +\frac{11 \log nm}{nm}\\
\leq & \frac{60}{nm}+\frac{120\sqrt{2r_1}}{\sqrt{nm}}\mathbb{E}\left[\sqrt{\log \mathcal{N}\left(1/nm, \left.\mathcal{G}\right|_{\mathcal{X}}, L^2 \right)}\right] + \frac{11 \log nm}{nm}\\
\leq & \frac{60}{nm}+\frac{r_1}{2}+\frac{120^2}{nm}\mathbb{E}\left[\log \mathcal{N}\left(1/nm, \left.\mathcal{G}\right|_{\mathcal{X}}, L^2 \right)\right] + \frac{11 \log nm}{nm}.
\end{aligned}
$$
That implies, there exists a constant $\xc\label{c120}$
$$
r_1 \leq \frac{\jc{c120}}{nm} \left( \mathbb{E}\left[ \log \mathcal{N}\left(1/nm, \left.\mathcal{G}\right|_{\mathcal{X}}, L^2 \right) \right]  + \log nm \right).
$$
Substituting $r^{*}_{nm}(\mathcal{G}) \leq r_1$ to \eqref{equ:145}, we have 
$$
\begin{aligned}
  \mathcal{E}_{\go}(\hat{g}_{nm})&
\leq 
\xc\left(\inf_{g\in\mathcal{G}} \int_{\Omega}\left(g-\go\right)^2d\mathcal{P}_{\x} 
+\frac{1+\jdb{b2}'^2 }{n}\right.
\\& + \left.\frac{1}{nm}\left(  \mathbb{E}\left[ \log \mathcal{N}\left(1/nm, \left.\mathcal{G}\right|_{\mathcal{X}}, L^2 \right) \right] + \log nm  \right)(1+\jdb{b2}'^2\log n+\jdb{b3}'^{2} \log nm)\right).\\
\end{aligned}
$$
Recalling that $\jdb{b1}\mathcal{G}=\mathcal{F}$ yields $\mathcal{N}\left(1/nm, \left.\mathcal{G}\right|_{\mathcal{X}}, L^2 \right) = \mathcal{N}\left(\jdb{b1}/nm, \left.\mathcal{F}\right|_{\mathcal{X}}, L^2 \right)$, the desired result follows.

(ii) Applying Theorem 9.4 in \citet{gyorfi2002distribution} and the relationship of covering number and packing number (see, for example, Lemma 9.2 in \citet{gyorfi2002distribution}), we have
$$
\log \mathcal{N}\left(\frac{\jdb{b1}}{nm}, \left.\mathcal{F}_{nm}\right|_{\mathcal{X}}, L^2 \right)
\leq  \vc(\mathcal{F}_{nm})\log ( 2e(nm)^2 \log 3e(nm)^2 ) +\log 3.
$$
Recalling that $\log nm\geq 1$, then (i) yields the desired result. 
\end{proof}

\begin{proof}[Proof of Corollary \ref{cor:networksize} ]

Theorem 6 in \citet{Bartlett2019nearly} tells us that the VC dimension of a ReLU neural with length $L$ and width $W$ can be bounded by
$$
\vc\leq \xc\label{c188} L^2W^2\log(LW),
$$
where \jc{c188} is a constant. 
Combining this conclusion and (ii) of Corollary \ref{cor:cor1} implies the desired result immediately.

\end{proof}

\begin{proof}[Proof of Theorem \ref{thm:holdernn}]
By Corollary 3.1 of \citet{huangjian}, we know there is a neural network $f'_{\nn}$ with width $\xc L\log L$ and depth $\xc W\log W$ that can uniformly approximate $\fo$ with an accuracy of  $\xc \left(LW\right)^{-2 s / d}$, where the constants do not depend on $L$ and $W$.
Using a truncation that could be realized by ReLU units easily, we can substitute the approximation error into Corollary \ref{cor:networksize}. 
Specifying $L, W$ satisfying 
$$
LW=\lfloor \jc{c55} (nm)^{ d/(4s+2d) } (\log nm)^{  - 4d/(2s+d) } \rfloor
$$ 
and the desired result follows.
\end{proof}

\begin{proof}[Proof of Theorem \ref{thm:low1}]

Given the sample sizes $n$ and $m$, we define  
$$J=\lfloor  \log_{2}(nm) /(2s+d) \rfloor$$ and $\xc\label{c555}$ to be a constant determined later.  
To consider B-splines with disjoint support sets,  we introduce a multi-index set 
\begin{equation}\label{equ:332}
\mathcal{V}= 
    \left\{  \bv=(v_{1}, v_{2},...,v_{d}) : -r+1 \leq  v_{i} \leq 2^{J} \text{ and } v_{i} = rt+1 \text{ for some } t\in \mathbb{Z} \right\}.
\end{equation}
One can verify that, for $\bv, \bv' \in \mathcal{V}, \bv\neq \bv'$, the support sets of $\mathcal{B}^{d}_{2^J,\bv}$ and $\mathcal{B}^{d}_{2^J,\bv'}$ are disjoint. 
Then we consider the following set 
\begin{equation}\label{equ:333}
    \left\{  \sum_{\bv} a_{J,\bv} 2^{-Js} \jc{c555} \mathcal{B}^{d}_{2^J,\bv} : a_{J,\bv}\in\{0,1\}, \bv \in  \mathcal{V} \right\}.
\end{equation}
It can be shown that the cardinality of $\mathcal{V}$ is given by $|\mathcal{V}|=( \lceil 2^{J} /r \rceil +1)^{d}\geq 2^{Jd} /r^d$.
Moreover, for any $\ba = \left(a_{J,\bv}: v \in \mathcal{V} \right)$,   the function $f_{\ba}=\sum_{\bv}  a_{J,\bv} 2^{-Js} \jc{c555} \mathcal{B}^{d}_{2^J,\bv}$ in the set defined by \eqref{equ:333}   belongs to $\mathcal{C}^{s}$, since it is a linear combination of  B-splines with order $r$ and disjoint support sets. In fact, 
$$\|f_{\ba}\|_{\mathcal{C}^{s}} \leq \xc \jc{c555} 2^{-Js} \|   \mathcal{B}^{d}_{2^J,\bo} \|_{\mathcal{C}^{s}} \leq \xc\label{c556} \jc{c555} < \infty.$$ 
Note that $f_{\ba}-f_{\ba'} = \sum_{\bv} (a_{J,\bv} - a_{J,\bv}' )  2^{-Js} \jc{c555} \mathcal{B}^{d}_{2^J,\bv}$ for two different functions $f_{\ba}$ and $f_{\ba'}$ from the set \eqref{equ:333}.
By the construction of $\mathcal{V}$, we have 
$$
\left\| f_{\ba}-f_{\ba'} \right\|_{L^2(\Omega)}^{2} = 2^{-2Js-Jd}  \xc\label{c550} \jc{c555}^2 \ham (\ba,\ba') ,
$$
where $\ham$ is the standard Hamming distance and $\jc{c550}=\| \mathcal{B}^{d}_{1,\bo} \|_{L^2(\Omega)}^{2}$.   
Using the Varshamov-Gilbert lemma, there is a subset $\mathcal{A} \subseteq \{0,1\}^{ |\mathcal{V}| }$ of cardinality $|\mathcal{A}|\geq \exp\{ |\mathcal{V}| /8 \}$, such that, for any $\ba, \ba' \in \mathcal{A}, \ba\neq \ba'$, $\ham(\ba, \ba') \geq |\mathcal{V}| /8 $ holds. 

Note that proving the lower bound only requires consideration of a specific case. Thus, we can first assume $Z$ is a Gaussian process and $\epsilon$ is a Gaussian noise with variance $\jdb{b3}$. 
Then,  the joint distribution of $(Y_{11}, Y_{12}, ..., Y_{1m},...,Y_{n1}, Y_{n2}, ..., Y_{nm})$ is multivariate normal distribution with mean 
\begin{equation}\label{equ:345}
    \mu_{\x,f_{\ba}}= (f_{\ba}(\x_{11}),f_{\ba}(\x_{11}),...,f_{\ba}(\x_{1m}),..., f_{\ba}(\x_{n1}),f_{\ba}(\x_{n1}),...,f_{\ba}(\x_{nm}) )
\end{equation}
and covariance matrix 
\begin{equation}\label{equ:346}
    C_{\x} =\diag\left\{ \left(\cov(Z_{1}(\x_{1j}, \x_{1k}))\right)_{1\leq j,k \leq m} \right\}_{1\leq i \leq n}  + \jdb{b3} I,
\end{equation}
conditional on  $\fo=f_{\ba}$ and measurement points $\left\{\x_{ij},  1\leq j\leq m  \right\} $. Hence, we can bound the Kullback-Leibler divergence of two probability measures from $\mathcal{P}_{f_{\ba}}$ to $\mathcal{P}_{f_{\ba'}}$ by 
$$
\begin{aligned}
   \kl(\mathcal{P}_{f_{\ba}},\mathcal{P}_{f_{\ba'}})= & \mathbb{E}\left[ ( \mu_{\x,f_{\ba}}-\mu_{\x,f_{\ba'}} ) C_{\x}^{-1} ( \mu_{\x,f_{\ba}}-\mu_{\x,f_{\ba'}} ) \right]\\
   \leq & \mathbb{E}\left[ ( \mu_{\x,f_{\ba}}-\mu_{\x,f_{\ba'}} ) (\jdb{b3} I)^{-1} ( \mu_{\x,f_{\ba}}-\mu_{\x,f_{\ba'}} ) \right]\\
   =& nm \jdb{b3}^{-1} \left\| f_{\ba}-f_{\ba'} \right\|_{L^2(\Omega)}^{2} \\
   =&  2^{-2Js-Jd} \xc \jc{c555}^2 nm \jdb{b3}^{-1} \ham (\ba,\ba')\\
   \leq & \xc\label{c551} \jc{c555}^2 \jdb{b3}^{-1}  (nm)^{d/(2s+d)},
\end{aligned}
$$
where $\jc{c551}$ is a constant independent of $n$ and $m$. 
Also, 
$$
\begin{aligned}
\left\| f_{\ba}-f_{\ba'} \right\|_{L^2(\Omega)}^{2} = &  2^{-2Js-Jd}\jc{c550}  \jc{c555}^2 \ham (\ba,\ba')\\
\geq & \xc\label{c552} \jc{c555}^2 (nm)^{-2s/(2s+d)},
\end{aligned}
$$
where $\jc{c552}$ is also a constant independent of $n$ and $m$. 
Applying Fano's inequality yields 
$$
\begin{aligned}
 \max_{a\in\mathcal{A}} \mathbb{E}_{f_{\ba}} \big[\int_{\Omega}(\tilde{f}-f_{\ba}&)^2d\mathcal{P}_{\x}\big] 
=
\max_{a\in\mathcal{A}}\mathbb{E}_{f_{\ba}}\left[ \left\| \tilde{f}-f_{\ba} \right\|_{L^2(\Omega)}^{2} \right]\\
\geq & \frac{\jc{c552} \jc{c555}^2 (nm)^{-2s/(2s+d)}}{4}\left( 
1-\frac{\jc{c551} \jc{c555}^2 \jdb{b3}^{-1}  (nm)^{d/(2s+d)} + \log 2}{2^{-d-3}r^{-d} (nm)^{d/(2s+d)} }
\right).
\end{aligned}
$$
Now we specify the value of $\jc{c555}$ such that $\jc{c556} \jc{c555} \leq 1$ and $\jc{c552} \jc{c555}^2 \leq 2^{-d-5}r^{-d} $. Then we have 
$$
1-\frac{\jc{c551} \jc{c555}^2 \jdb{b3}^{-1}  (nm)^{d/(2s+d)} + \log 2}{2^{-d-3}r^{-d} (nm)^{d/(2s+d)} } \geq \frac{1}{2}, 
$$
for $nm$ large enough. Taking into account the finite cases where $nm$ is small, we know that there is a constant $\xc\label{c557}$ independent of $n$ and $m$ such that 
$$
\max_{a\in\mathcal{A}} \mathbb{E}_{f_{\ba}} \left[\int_{\Omega}(\tilde{f}-f_{\ba})^2d\mathcal{P}_{\x}\right] 
\geq \jc{c557} (nm)^{-2s/(2s+d)}.
$$
for any estimator $\tilde{f}$. This implies 
\begin{equation}\label{equ:336}
    \inf_{\tilde{f}}\sup_{\mathcal{P}_{Z,\epsilon} \in \mathcal{Q}_{Z,\epsilon}} \mathbb{E}_{\fo} \left[\int_{\Omega}(\tilde{f}-\fo)^2d\mathcal{P}_{\x}\right] 
\geq \jc{c557} (nm)^{-2s/(2s+d)}.
\end{equation}

To show the remaining $1/n$ term in the lower bound,  it now suffices to consider $Z_{i}$ as constant functions over $\xx \in \Omega$ and the measurement noise $\epsilon_{ij}$ as zero. Assume that $\var(Z_{i}(\x)|\x)=\jdb{b2}$. Then the problem becomes the mean estimation from $n$ i.i.d. observations. It is well known that
the optimal rate is $\jdb{b2}/n$. It implies
\begin{equation}\label{equ:337}
    \inf_{\tilde{f}}\sup_{\mathcal{P}_{Z,\epsilon} \in \mathcal{Q}_{Z,\epsilon}} \mathbb{E}_{\fo} \left[\int_{\Omega}(\tilde{f}-\fo)^2d\mathcal{P}_{\x}\right] 
\geq \jdb{b2}n^{-1}.
\end{equation}
Then combining \eqref{equ:336} and \eqref{equ:337} completes the proof. 
\end{proof}

\begin{proof}[Proof of Theorem \ref{lem:hie}]
The proof is shown by induction on $l$.

First, the case $l=1$ is given by Corollary 3.1 in \citet{huangjian}.
Assume that the result holds for cases $<l$.
Now we consider the case for $l$. 
Suppose that $f \in \mathcal{H}^{l,\mathcal{G}}$ with $\mathcal{G}=((s,K);\mathcal{G}_{1},...,\mathcal{G}_{K} )$. 
Then we have $\gamma_{i}=\gamma(\mathcal{G}_{i})\leq \gamma/\min(1,s)$ and $s/K \leq \gamma$. 
By the definition of hierarchical composition model, we know $f(x)$ is equal to some $g\left(h_{1}(x), \ldots, h_{K}(x)\right)$, where
 $h_{i} \in \mathcal{H}^{l-1, \mathcal{G}_{i}}(\Omega)$ and $g\in \mathcal{C}^{s}$. 
Since $h_{i} \in \mathcal{H}^{l-1, \gamma_{i}}(\Omega)$, 
there exist positive numbers $M_l$ such that $h_{i}(x)\leq M_l$ for each $x\in \Omega$ and $i: 1\leq i\leq K$. Note that 
$$g\left(h_{1}(x), \ldots, h_{K}(x)\right)=g\left( 2M_l \cdot  \frac{h_{1}(x) + M_l }{2M_l}-\frac{1}{2}, \ldots, 2M_l \cdot  \frac{h_{K}(x) + M_l }{2M_l}-\frac{1}{2}\right)$$ 
and functions $g\left( 2M_l \cdot*  -1/2, \ldots, 2M_l \cdot *   -1/2\right)\in \mathcal{C}^{s}([0,1]^K)$, 
$h_{i}/M_l \in \mathcal{H}^{l-1, \gamma_{i}}(\Omega)$, without loss of generality, we assume that $0\leq h_{i}(x)\leq 1$ for each $i$ and $x\in \Omega$. 

Given any $L,W\geq 2$. By the inductive hypothesis, there exist neural networks 
$$\check h_{i} \in \nn\left(d, \xc\label{c515} L \log L, \xc\label{c516} W \log W, 1\right),$$
and a constant $\xc\label{c510}$ such that
\begin{equation}\label{equ:210}
|\check h_{i}(x)-h_{i}(x)|\leq \jc{c510} \left(LW\right)^{-2 \gamma_{i}},    
\end{equation}
for every $x\in \Omega$. 
Also, according to Corollary 3.1 in \citet{huangjian},  there exist a $\check g_{i} \in \nn\left(d, \xc\label{c513} L \log L, \xc\label{c514} W \log W, \jdb{b1}\right)$ and a constant $\xc\label{c511}$ such that 
\begin{equation}\label{equ:211}
|\check g(x)-g(x)|\leq \jc{c511}\left(LW\right)^{-2 s/k} ,     
\end{equation} for each $x\in [0,1]^K$. Then
$$
\begin{aligned}
& \left|\check g\left(\check h_{1}(x), \ldots, \check h_{K}(x)\right)-g\left(h_{1}(x), \ldots, h_{K}(x)\right)\right|\\
\leq & \left|\check g\left(\check h_{1}(x), \ldots, \check h_{K}(x)\right)-g\left(\check h_{1}(x), \ldots, \check h_{K}(x)\right)\right|\\
&\quad +\left| g\left(\check h_{1}(x), \ldots, \check h_{K}(x)\right)-g\left(h_{1}(x), \ldots, h_{K}(x)\right)\right|.
\end{aligned}
$$
By \eqref{equ:211}, the first term is bounded by $\jc{c511}\left(LW\right)^{-2 s/K}\leq \jc{c511}\left(LW\right)^{-2 \gamma}$. 
Additionally, by recalling \eqref{equ:210} and $g\in \mathcal{C}^{s}$, the second term is bounded by $$\xc (\sum_{i=1}^{K} \left(LW\right)^{-2 \gamma_{i}})^{\min(s,1)}\leq \xc  \left(LW\right)^{-2 \min_{i} \gamma_{i} \min(s,1)} \leq \xc \left(LW\right)^{-2 \gamma }. $$
To obtain the final result, we only need to check the size of $\check g \circ \left(\check h_{1}, \ldots, \check h_{K}\right)$. The width is bounded by $\max\left(  \jc{c516} W \log W, \jc{c514} W \log W \right)= \xc W \log W $ and the depth is bounded by $  \jc{c515} L \log L + \jc{c513} L \log L=\xc L \log L$. This completes the proof.

\end{proof}

\begin{proof}[Proof of Theorem \ref{thm:hieup}]
Combining Theorem \ref{thm:anisoholdernn} and Corollary \ref{cor:networksize} yields the desired result immediately. 
\end{proof}

\begin{proof}[Proof of Theorem \ref{thm:hielow}]
By the assumption, we know that there is a node $G^{\iota}_{1}=(\sso,\ddo)$ in the tree $\mathcal{G}$ such that $\sso^{*}/\ddo=\gamma$ and $\ddo\leq d$.
Assume that $G_{1}^{\iota-1}, G_{2}^{\iota-1}, ...,G_{\ddo}^{\iota-1}$ are children nodes of $G^{\iota}_{1}$ in the tree $\mathcal{G}$.
Let  $\mathcal{G}_{1}^{\iota}$ and $\mathcal{G}_{1}^{\iota-1}, \mathcal{G}_{2}^{\iota-1}, ...,\mathcal{G}_{\ddo}^{\iota-1}$ be the maximal sub-trees whose roots are $G_{1}^{\iota}$ and $G_{1}^{\iota-1}, G_{2}^{\iota-1}, ...,G_{\ddo}^{\iota-1}$, respectively. 
Let $h_{k}(\xx)=x_{k}$ be the $k$-th component of $\xx$. 
It is easy to verify that $h_{k} \in \mathcal{H}^{\iota-1, \mathcal{G}_{k}}(\Omega)$ for $1\leq k\leq \ddo$.
Therefore, $$
g\left(h_{1}(\xx),h_{2}(\xx), ..., h_{\ddo}(\xx)\right) 
\in \mathcal{H}^{\iota, \mathcal{G}_{1}^{\iota}}(\Omega)
$$ for any $g: [0,1]^{\ddo} \rightarrow \mathbb{R}$ satisfying $\|g\|_{\mathcal{C}^{\sso}}\leq M$. 
Denote $\ds=\prod_{(s,K)\in\mathcal{A}(G_{1}^{\iota})}\min(1,s)$ where $\mathcal{A}(G_{1}^{\iota})$ denotes the ancestor nodes of $G_{1}^{\iota}$.
Iteratively applying the fact that $m(\xx)=(h_{k}(\xx))^{\min (1,s)} \in \mathcal{C}^{s}$ yields
\begin{equation}\label{equ:347}
    f=(g\left(h_{1}(\xx),h_{2}(\xx), ..., h_{K}(\xx)\right))^{\ds} 
\in \mathcal{H}^{l, \mathcal{G}}(\Omega).
\end{equation}

Given the sample sizes $n$ and $m$, let  $J=\lfloor  \log_{2}(nm) /(2s_{0}^{*}+ \ddo ) \rfloor$ and $\xc\label{c569}$ be a constant to specific later.  
Let 
$$
\mathcal{V}= 
    \left\{  \bv=(v_{1}, v_{2},...,v_{\ddo}) : -r+1 \leq  v_{i} \leq 2^{J} \text{ and } v_{i} = rt+1 \text{ for some } t\in \mathbb{Z} \right\}.
$$
We consider the following set 
$$
    \left\{ f \text{ is of the form } \eqref{equ:347} \text{ with } g= \sum_{\bv} a_{J,\bv} 2^{-J\sso} \jc{c569} \mathcal{B}^{\ddo}_{2^J,\bv} , a_{J,\bv}\in\{0,1\}, \bv \in  \mathcal{V} \right\}.
$$
Note that the support sets of $\mathcal{B}^{\ddo}_{2^J,\bv}$ and $\mathcal{B}^{\ddo}_{2^J,\bv'}$ are disjoint for $\bv, \bv' \in \mathcal{V}, \bv\neq \bv'$, hence $\|g_{a}\|_{\mathcal{C}^{\sso}} \leq 2\|  2^{-J\sso} \jc{c569} \mathcal{B}^{\ddo}_{2^J,\bv}  \|_{\mathcal{C}^{\sso}}   \leq \xc\label{c570}\jc{c569}$.  
Moreover, we could rewrite the above set as
$$
  \left\{ f= \sum_{\bv} a_{J,\bv} 2^{-Js^{*}} \jc{c569}^{\ds} (\mathcal{B}^{\ddo}_{2^J,\bv})^{\ds}: a_{J,\bv}\in\{0,1\}, \bv \in  \mathcal{V} \right\}.
$$
The remaining discussion to obtain the $(nm)^{-\gamma/(2\gamma+1)}$ term in the lower bound is entirely similar to Theorem \ref{thm:low1}. 
Additionally, the proof of $n^{-1}$ term is also the same. 
Finally, by combining the results, we complete the proof. 
\end{proof}

\begin{proof}[Proof of Theorem \ref{thm:aniso}]
Define the high order difference operator $\Delta_{\boldsymbol{h}}^{r}(f)$ by 
$$\Delta_{\boldsymbol{h}}^{r}(f)(x)= \sum_{i=0}^{r} \binom{r}{i} (-1)^{r-i} f(x+i \boldsymbol{h}),
$$
when $x+i h$ are all within the domain. Otherwise, define $\Delta_{\boldsymbol{h}}^{r}(f)(x)=0$ if some $x+i h$ are out of the domain. 
For $p \in(0, \infty]$, we define the $r$-th modulus of smoothness of $f$ by 
$$
\omega_{r, p}(f, \boldsymbol{t})=\sup _{|h_{j}| \leq t_{j}}\left\|\Delta_{\boldsymbol{h}}^{r}(f)\right\|_{p}, 
$$
when $\boldsymbol{t}=(t_{1}, t_{2},...,t_{d})$ is a positive component vector. 
We also define the $i$-th partial modulus of smoothness by
$$
\omega_{r, p}^{(i)}(f, t)=\omega_{r, p}(f, t \boldsymbol{e}_{i}),
$$
when $t$ is a positive real number.

Step 1. 
We first approximate the anisotropic $\boldsymbol{s}$-Hölder function by B-spline quasi-interpolant representation. 
By the definition of modulus and $(4.76)$ in \citet{ginebook}, 
$$
\begin{aligned}
\sup_{t>0} \frac{ \omega_{r, \infty }^{(i)}(f, t) }{t^{s_{i}}}
\leq & 2^{r-\lfloor s_{i}\rfloor -1}  \sup_{t>0} \frac{ \omega_{\lfloor s_{i}\rfloor +1, \infty }^{(i)}(f, t) }{t^{s_{i}}}\\
\leq & 2^{r-\lfloor s_{i}\rfloor -1}  \sup_{t>0} \frac{ \omega_{1, \infty }^{(i)}(\partial^{\lfloor s_{i}\rfloor}_{i} f, t) }{t^{s_{i}-\lfloor s_{i}\rfloor}}\\
= & 2^{r-\lfloor s_{i}\rfloor -1}  \sup _{\substack{x, x' \in \Sigma \\  x' = x + te_{i}, t\neq 0}} \frac{\left|\partial^{\lfloor s_{i}\rfloor}_{i} f(x)-\partial^{\lfloor s_{i}\rfloor}_{i} f(x')\right|}{\|x-x'\|^{s_{i}-\lfloor s_{i} \rfloor}}.
\end{aligned}
$$
Therefore, there exists $\xc\label{c535}$ not depending on $f$ such that
\begin{equation}\label{equ:305}
    \|f\|_{\infty} + \sum_{i=1}^d \sup_{t>0} \frac{ \omega_{r, \infty }^{(i)}(f, t) }{t^{s_{i}}} \leq \jc{c535} \|f\|_{\mathcal{C}^{\boldsymbol{s}}} <\infty. 
\end{equation}
Denote $\sM=\max_{1\leq j \leq d}s_{j} $, $\sm=\min_{1\leq j \leq d}s_{j} $,  $s_{i}'= \sm/s_{i}$ 
and the harmonic mean
$$
\tilde{s}=\frac{d}{\frac{1}{s_{1}}+\frac{1}{s_{2}}+...+\frac{1}{s_{d}}}.
$$
Let 
$$
\bk(j)= \left(  2^{\lfloor j s_{1}^{\prime}\rfloor}, 2^{\lfloor j s_{2}^{\prime}\rfloor}, ...,  2^{\lfloor j s_{d}^{\prime}\rfloor}  \right),
$$ for each positive integer $j$. 
By Theorem 3.2.6 in \citet{LeisnerPhD}, there are real numbers $\left\{ a_{j,\bl} \right\}$ for $j\in \mathbb{Z}_{+}$ and multi-index $\bl \in \prod_{i=1}^{d} \left\{-r+1,-r+2, \ldots, 2^{\left\lfloor j s_{i}^{\prime}\right\rfloor}-1,2^{\left\lfloor j s_{i}^{\prime}\right\rfloor}\right\},
$
such that 
$$
\left\|f-  \sum_{\bl} a_{j,\bl}  \mathcal{B}^{d}_{\bk(j),\bl}  \right\|_{\infty} \leq \xc\label{c536} \omega_{r, \infty }\left(f,(2^{-j s_{1}^{\prime}},2^{-j s_{2}^{\prime}}, ..., 2^{-j s_{d}^{\prime}})\right).
$$
By the definition of $\omega_{r, \infty }$ and $\omega_{r, \infty}^{(i)}$, we know $$\omega_{r, \infty }\left(f,(2^{-j s_{1}^{\prime}},2^{-j s_{2}^{\prime}}, ..., 2^{-j s_{d}^{\prime}})\right)\leq \sum_{i=1}^{d} \omega_{r, \infty}^{(i)}(f, 2^{-j s_{i}^{\prime}}). $$
Further, using \eqref{equ:305}, we have
$$
\begin{aligned}
\sum_{i=1}^{d} \omega_{r, \infty}^{(i)}(f, 2^{-j s_{i}^{\prime}})
\leq  \sum_{i=1}^{d} \xc 2^{-j s_{i}^{\prime}s_{i}} \|f\|_{\mathcal{C}^{\boldsymbol{s}}}=\xc 2^{-j\sm } \|f\|_{\mathcal{C}^{\boldsymbol{s}}}.
\end{aligned}
$$
Combining several inequalities derived above yields that there is a constant $\xc\label{c538}$ independent of $f$ and $j$ such that 
\begin{equation}\label{equ:304}
    \left\|f-  \sum_{\bl} a_{j,\bl}  \mathcal{B}^{d}_{\bk(j),\bl}  \right\|_{\infty} \leq \jc{c538} 2^{-j\sm } \|f\|_{\mathcal{C}^{\boldsymbol{s}}}
\end{equation}
holds. In addition, by 
Lemma 3.2.2 and 
Lemma 3.2.5 in \citet{LeisnerPhD}, it follows that one exists $\xc\label{c539}$ and $\xc\label{c540}$ such that
\begin{equation}\label{equ:306}
    \left| a_{j,\bl} \right| \leq 
\jc{c539}\left\| \sum_{\bl} a_{j,\bl}  \mathcal{B}^{d}_{\bk(j),\bl}  \right\|_{\infty} \leq \jc{c540}  \|f\|_{\mathcal{C}^{\boldsymbol{s}}}
\end{equation}
for each $j$ and $\bl$.

Step 2.
Given $L$ and $W$, 
let $J=\lfloor \frac{2\tilde{s}}{d\sm } \log_{2}(LW) \rfloor$. We will  approximate $\sum_{\bl} a_{J,\bl}  \mathcal{B}^{d}_{\bk(J),\bl} $ by a fully connected ReLU neural network with the desired size. 
Since $J$ has been determined, let $\delta\in \left(0,\min_{1\leq i \leq d }2^{-\lfloor J s_{i}^{\prime}\rfloor}/3 \right]$ be a small number to be determined latter. 
Define a trifling region $\Upsilon(J,\delta) \subseteq \Omega$ by
$$
\Upsilon(J,\delta):=\bigcup_{t=1}^{d}\left\{\xx=(x_{1}, x_{2}, ...,x_{d}): \xx\in \Omega \text{ and } x_{t}\in \cup_{k=1}^{2^{J s_{t}^{\prime}}-1}\left(\frac{k}{2^{J s_{t}^{\prime}}}-\delta, \frac{k}{2^{J s_{t}^{\prime}}}\right) \right\}.
$$
Note that 
$$
\begin{aligned}
2^{J s_{i}^{\prime}} 
\leq (LW)^{\frac{2\tilde{s}}{d\sm } s_{i}^{\prime} } = (LW)^{\frac{2\tilde{s}}{d s_{i}} } \leq (LW)^{2},
\end{aligned}
$$
where the last inequality is due to $1/s_{i}\leq \sum_{i=1}^{d}1/s_{i}=d/\tilde{s}$. Hence, applying Lemma \ref{lem:bloc}, we get 
$\check{f}_{1,i}\in \nn\left(1,  4L+5 ,  4W+3 \right)$  satisfying
$$
\check{f}_{1,i}(x)=\kappa_{i}, \quad \text {when } x \in\left[\frac{\kappa_{i}}{2^{J s_{i}^{\prime}}}, \frac{\kappa_{i} +1}{2^{J s_{i}^{\prime}}}-\delta \cdot I_{\{\kappa_{i} \leq 2^{J s_{i}^{\prime}}-2\}}\right] ,
$$
for $\kappa_{i}=0,1, ..., 2^{J s_{i}^{\prime}}-1$ and $i=1,2,...,d$. Let 
$$
\check{f}_{2}(\xx)=(\check{f}_{1,1}(\xx),\check{f}_{1,2}(\xx),...,\check{f}_{1,d}(\xx)).
$$
Intuitively, this function divides $\Omega\backslash \Upsilon(J,\delta)$ into $\prod_{i=1}^{d} 2^{\lfloor J s_{i}^{\prime}\rfloor}$ blocks, i.e. for any $\xx=(x_{1},x_{2},...,x_{d}) \in \Omega\backslash \Upsilon(J,\delta)$, $\check{f}_{2}(\xx)=(\kappa_{1},\kappa_{2},...,\kappa_{d})$ if and only if $$\xx \in \prod_{i=1}^{d} \left[\frac{\kappa_{i}}{2^{J s_{i}^{\prime}}}, \frac{\kappa_{i}+1}{2^{J s_{i}^{\prime}}}-\delta \cdot I_{\{\kappa_{i} \leq 2^{J s_{i}^{\prime}}-2\}}\right]. $$

Note that the domain of a B-spline is bounded.  
Let $$\mathcal{U}:=\left\{ \bu= (u_{1},u_{2},...u_{d})\in \mathbb{Z}^{d}: 0\leq u_{i}\leq r-1 \text{ for } 1\leq i \leq d \right\}. $$
For  any $\xx \in \Omega\backslash \Upsilon(J,\delta)$, 
assume that 
$\check{f}_{2}(\xx)=(\kappa_{1},\kappa_{2},...,\kappa_{d}):=\bka$. 
 Then we know that $\mathcal{B}^{d}_{\bk(J),\bl}(\xx)\neq 0$ if and only if $\bka - \bl \in \mathcal{U}$. 
We should approximate the coefficients of these terms for each $\xx$. 
To apply Lemma \ref{lem:pof}, we should first construct a function $\check{f}_{3}(\bka )=\sum_{i=1}^{d}\kappa_{i}2^{J(i-1)}$ that maps all possible values of $\bka$ onto a subset of $\mathbb{Z}$ one by one.
Because $\bka$  only takes $\prod_{i=1}^{d} 2^{\lfloor J s_{i}^{\prime}\rfloor}\leq \prod_{i=1}^{d} (LW)^{\frac{2\tilde{s}}{d s_{i}} } = (LW)^{2}$ different values, by Lemma \ref{lem:pof}, there are 
$\check{f}_{4,\bu}\in \nn\left(1,  5(L+2)\log_{2}(4L) ,  16 \lceil \sM \rceil(W+1) \log_{2}(8W)\right)$  satisfying
\begin{equation}\label{equ:307}
    \left|\check{f}_{4,\bu} \circ \check{f}_{3} (\bka ) - a_{j,\bka+\bu} \right| \leq \xc (LW)^{-2 \lceil \sM \rceil}, 
\end{equation}
for $\bu\in \mathcal{U}$. 

By Lemma \ref{lem:bspl}, there is a $\check{f}_{5}\in \nn\left(d, \jc{c531} L , \jc{c532} W , 1\right)$ satisfying
$
|\check{f}_{5}(\xx)-\mathcal{B}^{d}_{ 1,\bo}|\leq \jc{c533}W^{-\jc{c534}L},
$
for every $\xx\in \left\{ \xx: x_{i}\in [0, r/k_{i}], 1\leq i\leq d \right\}.$
Also, by  Lemma 4.2 in \citet{lu2021deep}, there is a function $\check{f}_{6}$ implemented by a ReLU FNN with depth $L$ and width $9W + 1$ could approximate product function in a determined bounded domain and suitable precision.
Let 
$$ 
\check{f}_{7}(\xx) = \sum_{\bu \in\mathcal{U}}
\check{f}_{6}\left(
\check{f}_{4,\bu} \circ \check{f}_{3} \circ \check{f}_{2}(\xx), 
\check{f}_{5}( 2^{\bk(J)} \xx - \check{f}_{2}(\xx)-\bu )
\right).
$$ 
The approximation error could be bounded by 
$$
\begin{aligned}
&\left|\check{f}_{7}(\xx)- \sum_{\bl} a_{J,\bl}  \mathcal{B}^{d}_{\bk(J),\bl} \right|
=\left|\check{f}_{7}(\xx)- \sum_{\bu \in\mathcal{U}} a_{J,\bka+\bu}  \mathcal{B}^{d}_{\bk(J),\bka+\bu} \right|\\
\leq & \sum_{\bu \in\mathcal{U}} \left|\check{f}_{6}\left(
\check{f}_{4,\bu} \circ \check{f}_{3} \circ \check{f}_{2}(\xx), 
\check{f}_{5}( 2^{\bk(J)} \xx - \check{f}_{2}(\xx)-\bu ) 
\right)-  a_{J,\bka+\bu}  \mathcal{B}^{d}_{\bk(J),\bka+\bu} \right|\\
\leq &  \sum_{\bu \in\mathcal{U}}  \left|\check{f}_{6}\left(
\check{f}_{4,\bu} \circ \check{f}_{3} \circ \check{f}_{2}(\xx), 
\check{f}_{5}( 2^{\bk(J)} \xx - \check{f}_{2}(\xx)-\bu ) 
\right)- \check{f}_{4,\bu} \circ \check{f}_{3} \circ \check{f}_{2}(\xx)\times 
\check{f}_{5}( 2^{\bk(J)} \xx - \check{f}_{2}(\xx)-\bu )\right|\\
&+\sum_{\bu \in\mathcal{U}}  \left|
\check{f}_{4,\bu} \circ \check{f}_{3} \circ \check{f}_{2}(\xx)\times 
\check{f}_{5}( 2^{\bk(J)} \xx - \check{f}_{2}(\xx)-\bu ) -  a_{J,\bka+\bu}  \mathcal{B}^{d}_{\bk(J),\bka+\bu} \right|.
\end{aligned}
$$
By  \eqref{equ:306}, \eqref{equ:307} and  Lemma 4.2 in \citet{lu2021deep}, the first term can be bounded by $\xc\label{c541} W^{-L} $. The second term is bounded by
$$
\begin{aligned}
& \sum_{\bu \in\mathcal{U}}  \left|
\check{f}_{4,\bu} \circ \check{f}_{3} \circ \check{f}_{2}(\xx)\times 
\check{f}_{5}( 2^{\bk(J)} \xx - \check{f}_{2}(\xx)-\bu ) -  a_{J,\bka+\bu}  \mathcal{B}^{d}_{\bk(J),\bka+\bu} \right|\\
\leq & \sum_{\bu \in\mathcal{U}}  \left|
\check{f}_{4,\bu} \circ \check{f}_{3} \circ \check{f}_{2}(\xx)
 -  a_{J,\bka+\bu}   \right| \cdot 
| \check{f}_{5}( 2^{\bk(J)} \xx - \check{f}_{2}(\xx)-\bu )|\\
& +\sum_{\bu \in\mathcal{U}}  \left|
\check{f}_{5}( 2^{\bk(J)} \xx - \check{f}_{2}(\xx)-\bu ) -    \mathcal{B}^{d}_{\bk(J),\bka+\bu} \right| \cdot | a_{J,\bka+\bu}|\\
\leq & \xc \label{c542} (LW)^{-2 \lceil \sM \rceil}+ \xc \label{c543} W^{-\jc{c534}L}.
\end{aligned}
$$
Therefore, 
$$
\left|\check{f}_{7}(\xx)- \sum_{\bl} a_{J,\bl}  \mathcal{B}^{d}_{\bk(J),\bl} \right| \leq \jc{c541} W^{-L} + \jc{c542} (LW)^{-2 \lceil \sM \rceil}+ \jc{c543} W^{-\jc{c534}L} \leq \xc\label{c545} (LW)^{-2\tilde{s}/d}, 
$$
for $\xx \in \Omega\backslash \Upsilon(J,\delta)$.
Recalling \eqref{equ:304} and $J=\lfloor \frac{2\tilde{s}}{d\sm } \log_{2}(LW) \rfloor$, we have 
$$\left|f(x)-   \sum_{\bl} a_{J,\bl} \mathcal{B}^{d}_{\bk(J),\bl} (x) \right|\leq  \jc{c538} 2^{-J\sm } \|f\|_{\mathcal{C}^{\boldsymbol{s}}}\leq \xc (LW)^{-2 \tilde{s}/d }.$$
Combining the previous two results by triangle inequality, it follows that 
$$
\left| \check{f}_{7}(\xx)-f(\xx)\right| \leq \xc\label{c546} (LW)^{-2 \tilde{s}/d },
$$
where $\xx \in \Omega\backslash \Upsilon(J,\delta)$ and $\jc{c546}$ is a constant independent of $L$ and $W$. 
And the size of an implementation can be determined by Figure \ref{fig:101} easily.

\begin{figure}
    \centering
    
    \begin{tikzpicture}
		[L1Node/.style={rounded corners,draw=black!100,
		minimum size=20pt},
		L2Node/.style={rounded corners,draw=pink!50,fill=pink!20,very thick, minimum size=20pt}]
		
		\foreach \x in {1}  
		\node[L1Node] (a_1) at (0,0){$\xx$};
		
		\node[L1Node] (a_2) at (1.3,1.3) {$\xx$};
		\node[L1Node] (a_3) at (1.3,-0.9) {$\bka$};
		
		\node[L1Node] (a_4) at (3,1.3) {$\xx$};
		\node[L1Node] (a_5) at (3,-0.15) {$\bka$};
		\node[L1Node] (a_6) at (3,-1.65) {$\check{f}_{3}\left(\bka \right)$};
		
		\node[L1Node] (a_7) at (6,1.8) {$\check{\mathcal{B}}^{d}_{\bk(J),\bka+\bu_{1}}( \xx )$};
		\node[L1Node] (a_8) at (6,0.25) {$\check{\mathcal{B}}^{d}_{\bk(J),\bka+\bu_{r^d}}( \xx )$};
		
		\node[L1Node] (a_9) at (6,-0.9) {$\check{a}_{J,\bka+\bu_{1}}$};
		\node[L1Node] (a_10) at (6,-2.2) {$\check{a}_{J,\bka+\bu_{r^d}}$};
		
		\node[L1Node] (a_11) at (10,1) {$\check{f}_{6}(\check{a}_{J,\bka+\bu_{1}}, \check{\mathcal{B}}^{d}_{\bk(J),\bka+\bu_{1}}( \xx ) )$};
		\node[L1Node] (a_12) at (10,-1.4) {$\check{f}_{6}(\check{a}_{J,\bka+\bu_{r^d}}, \check{\mathcal{B}}^{d}_{\bk(J),\bka+\bu_{r^d}}( \xx ) )$};
		
		\node[L1Node] (a_13) at (13,0) {$\check{f}_{7}(\xx)$};

		\draw[-stealth](a_1)--(a_2);
		\draw[-stealth](a_1)--(a_3);
		\draw[-stealth](a_2)--(a_4);
		\draw[-stealth](a_3)--(a_5);
		\draw[-stealth](a_3)--(a_6);
		\draw[-stealth](a_4)--(a_7);
		\draw[-stealth](a_4)--(a_8);
		\draw[-stealth](a_5)--(a_7);
		\draw[-stealth](a_5)--(a_8);
		\draw[-stealth](a_6)--(a_9);
		\draw[-stealth](a_6)--(a_10);
		\draw[-stealth](a_7)--(a_11);
		\draw[-stealth](a_8)--(a_12);
		\draw[-stealth](a_9)--(a_11);
		\draw[-stealth](a_10)--(a_12);
		\draw[-stealth](a_11)--(a_13);
		\draw[-stealth](a_12)--(a_13);
		
		\node[](w_1) at (0.7,-0.3){$\check{f}_{2}$};
		\node[](w_2) at (2.1,-1.1){$\check{f}_{3}$};
		\node[](w_3) at (4.5,1.35){$\check{f}_{5,\bu_{1}}$};
		\node[](w_4) at (4.3,0.34){$\check{f}_{5,\bu_{r^d}}$};
		\node[](w_5) at (4.4,-1.1){$\check{f}_{4,\bu_{1}}$};
		\node[](w_6) at (4.42,-1.68){$\check{f}_{4,\bu_{r^d}}$};
		\node[](w_7) at (6,1.15){$\vdots$};
		\node[](w_8) at (6,-1.5){$\vdots$};
		\node[](w_9) at (10,0){$\vdots$};
		\node[](w_9) at (11.6,0){add};

	\end{tikzpicture}
    
    \caption{An illustration of the sub-network architecture implementing $\check{f}_{7}(\xx)$. Here, we denote $\check{f}_{5}( 2^{\bk(J)} \xx - \check{f}_{2}(\xx)-\bu )$, $\check{f}_{5}( 2^{\bk(J)} \xx -  \check{f}_{2}(\xx)-\bu )$ and $\check{f}_{4,\bu} \circ \check{f}_{3} (\bka ) $ as $\check{f}_{5,\bu}(\xx)$, $\check{\mathcal{B}}^{d}_{\bk(J),\bka+\bu}(\xx)$ and  $\check{a}_{J,\bka+\bu}$, respectively.}
    
    \label{fig:101}
\end{figure}

Step 3. 
The median function of three arguments could be implemented by a ReLU neural network with width $14$ and depth $2$ because
$$
\begin{aligned}
\operatorname{mid}\left(t_{1}, t_{2}, t_{3}\right)=&\relu\left(t_{1}+t_{2}+t_{3}\right)-\relu\left(-t_{1}-t_{2}-t_{3}\right)\\&-\max \left\{t_{1}, t_{2}, t_{3}\right\}-\min \left\{t_{1}, t_{2}, t_{3}\right\},
\end{aligned}
$$
where 
$$
\begin{aligned}
\max \left\{t_{1}, t_{2}, t_{3}\right\}& =\max \left\{\max \left\{t_{1}, t_{2}\right\}, t_{3}\right\}\\ & =\max \left\{\max \left\{t_{1}, t_{2}\right\}, \relu\left(t_{3}\right)-\relu\left(-t_{3}\right)\right\},
\end{aligned}
$$
$$
\begin{aligned}
\max \left\{t_{1}, t_{2}\right\}=\frac{1}{2}(& \relu\left(t_{1}+t_{2}\right)-\relu\left(-t_{1}-t_{2}\right)\\&+\relu\left(t_{1}-t_{2}\right)+\relu\left(t_{2}-t_{1}\right)).
\end{aligned}
$$ 
and the min function can be implemented similarly.

Set $\check{f}_{8,0}=\check{f}_{7}$, then we inductively define
$$
\check{f}_{8,i}(x):=\operatorname{mid}\left(\check{f}_{8,i-1}\left(x-\delta e_{i}\right), \check{f}_{8,i-1}(x), \check{f}_{8,i-1}\left(x+\delta e_{i}\right)\right), 
$$
for $i=1, 2,..., d$. 
We also introduce 
$$
\Upsilon_{i}(J,\delta):=\bigcup_{t=i+1}^{d}\left\{\xx=(x_{1}, x_{2}, ...,x_{d}): \xx\in \Omega \text{ and } x_{t}\in \cup_{k=1}^{2^{J s_{t}^{\prime}}-1}\left(\frac{k}{2^{J s_{t}^{\prime}}}-\delta, \frac{k}{2^{J s_{t}^{\prime}}}\right) \right\},
$$
for $i=0,1,..., d$. It is trivial that $\Upsilon_{0}(J,\delta)=\Upsilon(J,\delta)$ and $\Upsilon_{d}(J,\delta)=\emptyset$.

To control the approximation error $\check{f}_{8,d}$ in $\Sigma$,  we prove the following stronger claim by induction 
\begin{equation}\label{equ:308}
    |\check{f}_{8,i}(\xx) - f(\xx)|\leq \jc{c546} (LW)^{-2 \tilde{s}/d } + i \|f\|_{\mathcal{C}^{\boldsymbol{s}}} \delta^{\min(\sm,1)}
\end{equation}
for $ \xx \in \Omega\backslash \Upsilon_{i}(J,\delta)$ and $i=0,1,..., d$.
In fact, 
when $i=0$, it is given by the construction in the previous step. 
Suppose that the claim holds for $<i$, consider the case for $i=t$. 
Recalling the definition of $\Upsilon_{i}(J,\delta)$, we know there are at most one of $\xx-\delta e_{i}$, $\xx$ and $\xx+\delta e_{i}$ is in $\Upsilon_{i}(J,\delta) $.
Hence, by inductive assumption, there is at most one of $ \check{f}_{8,i-1}(\xx-\delta e_{i}) -f(\xx-\delta e_{i})$, $ \check{f}_{8,i-1}(\xx) - f(\xx)$ and $  \check{f}_{8,i-1}(\xx+\delta e_{i})-f(\xx+\delta e_{i})$ is larger than 
$\jc{c546} (LW)^{-2 \tilde{s}/d } + (i-1) \|f\|_{\mathcal{C}^{\boldsymbol{s}}} \delta^{\min(\sm,1)}$.
Combining with 
$
\left|f\left(\xx \pm \delta e_{i}\right)-f(\xx)\right| \leq \|f\|_{\mathcal{C}^{\boldsymbol{s}}}  \delta^{\min(s_{i},1)}\leq \|f\|_{\mathcal{C}^{\boldsymbol{s}}}  \delta^{\min(\sm ,1)}  
$ 
yields
that at least two of 
$|\check{f}_{8,i-1}(\xx-\delta e_{i})- f(\xx) |$, $| \check{f}_{8,i-1}(\xx) -f(\xx)|$ and $|  \check{f}_{8,i-1}(\xx+\delta e_{i})|-f(\xx)$ are smaller than $$\jc{c546} (LW)^{-2 \tilde{s}/d } + i \|f\|_{\mathcal{C}^{\boldsymbol{s}}} \delta^{\min(\sm,1)}. $$
This implies that at least two of $\check{f}_{8,i}(\xx-\delta e_{i})$, $\check{f}_{8,i}(\xx)$ and $\check{f}_{8,i}(\xx+\delta e_{i})$ are in the interval 
$$
\left[ f(\xx)-\jc{c546} (LW)^{-2 \tilde{s}/d } - i \|f\|_{\mathcal{C}^{\boldsymbol{s}}} \delta^{\min(\sm,1)},  f(\xx)+\jc{c546} (LW)^{-2 \tilde{s}/d } + i \|f\|_{\mathcal{C}^{\boldsymbol{s}}} \delta^{\min(\sm,1)} \right].
$$
Therefore, middle value $\check{f}_{8,i}(x)= \operatorname{mid}\left( \check{f}_{8,i-1}(\xx-\delta e_{i}),\check{f}_{8,i}(\xx), \check{f}_{8,i-1}(\xx+\delta e_{i}) \right)$ belongs to the above interval, i.e. $$\left| \check{f}_{8,i}(x)-f(x)\right| \leq \jc{c546} (LW)^{-2 \tilde{s}/d } + i \|f\|_{\mathcal{C}^{\boldsymbol{s}}} \delta^{\min(\sm,1)}.$$
Hence, the case for $i$ holds. This completes the proof of the claim \eqref{equ:308} by induction. In particularly, let $i=d$ and it follows 
$$
|\check{f}_{8,d}(\xx) - f(\xx)|\leq \jc{c546} (LW)^{-2 \tilde{s}/d } + d \|f\|_{\mathcal{C}^{\boldsymbol{s}}} \delta^{\min(\sm,1)},
$$
for $\xx \in \Omega$.
Since the exact value of $\delta$ has not yet been determined, we take it small enough to $\delta^{\min(\sm,1)}\leq  (LW)^{-2\tilde{s}/d}$. Hence, 
$$
|\check{f}_{8,d}(\xx) - f(\xx)|\leq \jc{c546} (LW)^{-2 \tilde{s}/d } + d \|f\|_{\mathcal{C}^{\boldsymbol{s}}} \delta^{\min(\sm,1)}\leq \xc\label{c548} (LW)^{-2 \tilde{s}/d },
$$
for $\xx \in \Omega$ and $\jc{c548}$ which is a constant independent of $L$ and $W$. 
The recursive construction of the neural network makes the size verification simple, and the desired result follows.

\end{proof}

\begin{proof}[Proof of Theorem \ref{thm:anisoholdernn}]
By the proof of Theorem 6.2 in \citet{huangjian}, we know that for any $f \in \mathcal{C}^{s}\left(\Omega\right)$ satisfying $\|f\|_{\infty}\leq \jdb{b1}$ 
and any $L,W\geq 3$, there is some ReLU neural network function $f_{\nn}\in \nn\left(d, \xc\label{c503} L \log L, \xc\label{c504} W \log W, \jdb{b1}\right)$ satisfying
$$
|f_{\nn}(x)-f(x)| \leq \xc\label{c505}\left(LW\right)^{-2 s / d_{\mathcal{M}}},
$$
for every $x\in \mathcal{M}$.
Here, $\jc{c503}$ and $\jc{c504}$ are constants independent of $L$, $W$ and $f$. 
Combining the above conclusion and Corollary \ref{cor:networksize}, the desired result follows immediately. 

\end{proof}

\begin{proof}[Proof of Theorem \ref{thm:anisolow}]

This proof is similar to the proof of Theorem \ref{thm:low1}. We first recall some notations. 
Let $\sM=\max_{1\leq j \leq d}s_{j} $, $\sm=\min_{1\leq j \leq d}s_{j} $,  $s_{i}'= \sm/s_{i}$ 
and $\tilde{s}$ is the harmonic mean of $s_{i}$ for $1\leq i \leq d$. 
Denote 
$$
\bk(j)= \left(  2^{\lfloor j s_{1}^{\prime}\rfloor}, 2^{\lfloor j s_{2}^{\prime}\rfloor}, ...,  2^{\lfloor j s_{d}^{\prime}\rfloor}  \right),
$$
for each positive integer $j$. 
Given the sample sizes $n$ and $m$, let  $J=\lfloor  \log_{2}(nm) \cdot \ts /(2 \ts\sm+d\sm ) \rfloor$ and $\xc\label{c560}$ be a constant determined next.   
Then we define a multi-index set
\begin{equation}\label{equ:340}
\mathcal{V}= 
    \left\{  \bv=(v_{1}, v_{2},...,v_{d}) : -r+1 \leq  v_{i} \leq 2^{\lfloor J s_{i}^{\prime}\rfloor} \text{ and } v_{i} = rt+1 \text{ for some } t\in \mathbb{Z} \right\}.
\end{equation}
and the B-spline linear space
\begin{equation}\label{equ:341}
    \left\{  \sum_{\bv} a_{J,\bv} 2^{-J\sm} \jc{c560} \mathcal{B}^{d}_{\bk(J),\bv} : a_{J,\bv}\in\{0,1\}, \bv \in  \mathcal{V} \right\}.
\end{equation}
We know that $|\mathcal{V}|= \prod_{i=1}^{d} ( \lceil k_{i} /r \rceil +1)^{d}\geq 2^{Jd \sm/\ts } /(2r)^d \geq (nm)^{d/(2\ts+d) }2^{-d(\sm+\ts)/\ts}r^{-d}$ $=\xc\label{c561}(nm)^{d/(2\ts+d)}$, for a constant $\jc{c561}$.
Consider $x' = x + te_{i}$.
Note that 
$$
\begin{aligned}
& 2^{-J\sm}   \frac{\left|\partial^{\lfloor s_{i} \rfloor}_{i} \mathcal{B}^{d}_{\bk(J),\bv}(\xx)-\partial^{\lfloor s_{i} \rfloor}_{i} \mathcal{B}^{d}_{\bk(J),\bv}(\xx')\right|}{\|\xx-\xx'\|^{s_{i}-\lfloor s_{i} \rfloor}}\\
= &   2^{-J\sm} k_{i}^{s_{i}}  \frac{\left|\partial^{\lfloor s_{i} \rfloor}_{i} \mathcal{B}^{d}_{1,\bo}( \bk(J) \xx - \bv )-\partial^{\lfloor s_{i} \rfloor}_{i} \mathcal{B}^{d}_{1,\bo}( \bk(J) \xx' - \bv )\right|}{ k_{i}^{s_{i}-\lfloor s_{i} \rfloor} \|\xx-\xx'\|^{s_{i}-\lfloor s_{i} \rfloor}}\\
\leq &  \frac{\left|\partial^{\lfloor s_{i} \rfloor}_{i} \mathcal{B}^{d}_{1,\bo}( \bk(J) \xx - \bv )-\partial^{\lfloor s_{i} \rfloor}_{i} \mathcal{B}^{d}_{1,\bo}( \bk(J) \xx' - \bv )\right|}{ k_{i}^{s_{i}-\lfloor s_{i} \rfloor} \|\xx-\xx'\|^{s_{i}-\lfloor s_{i} \rfloor}}.
\end{aligned}
$$
It implies that 
$$
 \|   2^{-J\sm} \mathcal{B}^{d}_{\bk(J),\bv} \|_{\mathcal{C}^{\boldsymbol{s}}}\leq \|\mathcal{B}^{d}_{1,\bo}\|_{\mathcal{C}^{\boldsymbol{s}}}.
$$
Recall that the support sets of each B-spline disjoint. Therefore, we have 
$$
\|   \sum_{\bv} a_{J,\bv} 2^{-J\sm} \jc{c560} \mathcal{B}^{d}_{\bk(J),\bv}  \|_{\mathcal{C}^{\boldsymbol{s}}}\leq 2\jc{c560} \|\mathcal{B}^{d}_{1,\bo}\|_{\mathcal{C}^{\boldsymbol{s}}}=\jc{c560}\xc\label{c562}<\infty.
$$
for any $\ba = \left(a_{J,\bv}: v \in \mathcal{V} \right)$.
In addition, we have 
$$
\xc\label{c563} \jc{c560}^2 \ham (\ba,\ba')2^{-2J\sm -Jd\sm/\ts}\leq 
\left\| f_{\ba}-f_{\ba'} \right\|_{L^2(\Omega)}^{2} \leq 
\xc\label{c564} \jc{c560}^2 \ham (\ba,\ba')2^{-2J\sm -Jd\sm/\ts},
$$
where $\jc{c563}$ and $\jc{c564}$ are constants independent of $n$ and $m$. 
We still use the Varshamov-Gilbert bound, and know that there exists a subset $\mathcal{A} \subseteq \{0,1\}^{ |\mathcal{V}| }$ which satisfies $|\mathcal{A}|\geq \exp\{ |\mathcal{V}| /8 \}$, and $\ham(\ba, \ba') \geq |\mathcal{V}| /8 $ for $\ba, \ba' \in \mathcal{A}, \ba\neq \ba'$.

Similar with the previous proof, to show the $(nm)^{-2\ts/(\ts+d)}$ term in the lower bound, 
we assume that  $Z$ is a Gaussian process and the measurement noise $\epsilon$ is normal random variables with variance $\jdb{b3}$. 
Hence, using the KL divergence of two Gaussian probability measures, we have
$$
\begin{aligned}
   \kl(\mathcal{P}_{f_{\ba}},\mathcal{P}_{f_{\ba'}})= & \mathbb{E}\left[ ( \mu_{\x,f_{\ba}}-\mu_{\x,f_{\ba'}} ) C_{\x}^{-1} ( \mu_{\x,f_{\ba}}-\mu_{\x,f_{\ba'}} ) \right]\\
   \leq &  nm \jdb{b3}^{-1} \left\| f_{\ba}-f_{\ba'} \right\|_{L^2(\Omega)}^{2} \\
   \leq & \jc{c564} \jc{c560}^2 nm \jdb{b3}^{-1}  2^{-2J\sm -Jd\sm/\ts} |\mathcal{V}|\\
   \leq & \xc\label{c565} \jc{c560}^2 \jdb{b3}^{-1}  (nm)^{d/(2\ts+d)},
\end{aligned}
$$
where $ \mu_{\x,f_{\ba}}$,  $ \mu_{\x,f_{\ba'}}$  and $C_{\x}$ are given by \eqref{equ:345} and  \eqref{equ:346}, and  
$\jc{c565}$ is a constant independent of $n$ and $m$. 
We can also check that  
$$
\begin{aligned}
\left\| f_{\ba}-f_{\ba'} \right\|_{L^2(\Omega)}^{2} \geq  & \jc{c563}  \jc{c560}^2 \ham (\ba,\ba') 2^{-2Js-Jd}\\
\geq & \xc\label{c566} \jc{c560}^2 (nm)^{-2\ts/(2\ts+d)},
\end{aligned}
$$
with $\ba,\ba'\in \mathcal{A}$ and $\ba\neq \ba'$. 
Here $\jc{c566}$ is a constant independent of sample size.
Then it follows from Fano's inequality that
$$
\begin{aligned}
& \max_{a\in\mathcal{A}} \mathbb{E}_{f_{\ba}} \left[\int_{\Omega}(\tilde{f}-f_{\ba})^2d\mathcal{P}_{\x}\right] =
\max_{a\in\mathcal{A}}\mathbb{E}_{f_{\ba}}\left[ \left\| \tilde{f}-f_{\ba} \right\|_{L^2(\Omega)}^{2} \right]\\
\geq & \frac{\jc{c566} \jc{c560}^2 (nm)^{-2s/(2s+d)}}{4}\left( 
1-\frac{\jc{c565} \jc{c560}^2 \jdb{b3}^{-1}  (nm)^{d/(2\ts+d)} + \log 2}{2^{-3}\jc{c561} (nm)^{d/(2\ts+d)} }
\right), 
\end{aligned}
$$
for any estimator $\tilde{f}$. 
Let $\jc{c560}$ small enough such that  $\jc{c560}\jc{c562}\leq 1$ and $\jc{c565} \jc{c560}^2 \jdb{b3}^{-1} \leq 2^{-5} \jc{c561} $. Then 
$$
1-\frac{\jc{c565} \jc{c560}^2 \jdb{b3}^{-1}  (nm)^{d/(2\ts+d)} + \log 2}{2^{-3}\jc{c561} (nm)^{d/(2\ts+d)} } \geq \frac{1}{2}, 
$$
when $nm$ is larger than some fixed number. 
Therefore, there exists a constant $\xc\label{c567}$ free of $n$ and $m$ such that 
$$
\begin{aligned}
    \inf_{\tilde{f}}\sup_{\mathcal{P}_{Z,\epsilon} \in \mathcal{Q}_{Z,\epsilon} } \mathbb{E}_{\fo} \left[\int_{\Omega}(\tilde{f}-\fo)^2d\mathcal{P}_{\x}\right] \geq &
\inf_{\tilde{f}} \max_{a\in\mathcal{A}} \mathbb{E}_{f_{\ba}} \left[\int_{\Omega}(\tilde{f}-f_{\ba})^2d\mathcal{P}_{\x}\right] \\
\geq & \xc (nm)^{-2s/(2s+d)}.
\end{aligned}
$$
The infimum here is correct because Fano inequality holds for any estimator.
Moreover, after the same argument as \eqref{equ:337}, we have
$$
    \inf_{\tilde{f}}\sup_{\mathcal{P}_{Z,\epsilon} \in \mathcal{Q}_{Z,\epsilon}} \mathbb{E}_{\fo} \left[\int_{\Omega}(\tilde{f}-\fo)^2d\mathcal{P}_{\x}\right] 
\geq \jdb{b2}n^{-1}.
$$
Then the desired result follows.
\end{proof}
\begin{proof}[Proof of Theorem \ref{thm:lowlinear}]
Using Lemma \ref{lem:lowerlinearlemma}, the proof is given by the construction of partition $\mathcal{A}$ in (ii-b) of the proof of Theorem 5 in \citet{suzuki2021deep} with $\tilde{d} = 1$. 
\end{proof}

\section{Proofs of Lemmas}\label{app:prflem}

\begin{proof}[Proof of Lemma \ref{lem:lema1}]
    For the purpose of bounding $\chf$ with its empirical version sharply, we use the peeling technique combined with Talagrand's inequality. Assume that $r$ is a positive number. 
We first consider
$$
\left|  \frac{1}{n m} \sum_{i=1}^{n} \sum_{j=1}^{m} \frac{\chf  -\left(\hat f\left(\x_{i j}\right)- \fo\left(\x_{i j}\right)\right)^{2}}{\chf  +r} 
\right|,
$$
which is bounded by 
$$
\sup_{f\in \mathcal{F}}
\left| 
\frac{1}{n m} \sum_{i=1}^{n} \sum_{j=1}^{m} \frac{\cf  -\left(f\left(\x_{i j}\right)- \fo\left(\x_{i j}\right)\right)^{2}}{\cf  +r} 
\right|.
$$
As a result of the symmetry and convexity, we have
$$
\begin{aligned}
&\mathbb{E}\left[\sup _{f \in \mathcal{F}}\left|\frac{1}{n m} \sum_{i=1}^{n} \sum_{j=1}^{m} \frac{\cf  -\left(f\left(\x_{i j}\right)- \fo\left(\x_{i j}\right)\right)^{2}}{\cf  +r}\right|\right]\\
\leq & 2 \mathbb{E}\left[\sup _{f \in \mathcal{F}}\left|\frac{1}{n m} \sum_{i=1}^{n} \sum_{j=1}^{m} \frac{\sigma_{i j}\left(f\left(\x_{i j}\right)- \fo\left(\x_{i j}\right)\right)^{2}}{\cf  +r}\right|\right],
\end{aligned}
$$
which is usually called as the symmetrization conclusion of Rademacher process.
Let the real number $r$ be greater than $r^{*}$. 
Note that $\mathcal{F}$ is the union of $\mathcal{F}(r)$ and $\mathcal{F}(4^{k}r) \backslash \mathcal{F}(4^{k-1}r)$ for $\ k\geq 1$, we have 
\begin{align*}
& \mathbb{E}\left[\sup _{f \in \mathcal{F}}\left|\frac{1}{n m} \sum_{i=1}^{n} \sum_{j=1}^{m} \frac{\sigma_{i j}\left(f\left(\x_{i j}\right)- \fo\left(\x_{i j}\right)\right)^{2}}{\cf  +r}\right|\right]\\
\leq &  \mathbb{E}\left[\sup _{f \in \mathcal{F}(r)}\left|\frac{1}{n m} \sum_{i=1}^{n} \sum_{j=1}^{m} \frac{\sigma_{i j}\left(f\left(\x_{i j}\right)- \fo\left(\x_{i j}\right)\right)^{2}}{\cf  +r}\right|\right]\\
&+\sum_{k=1}^{\infty}\mathbb{E}\left[\sup _{f \in \mathcal{F}(4^{k}r) \backslash \mathcal{F}(4^{k-1}r)}\left|\frac{1}{n m} \sum_{i=1}^{n} \sum_{j=1}^{m} \frac{\sigma_{i j}\left(f\left(\x_{i j}\right)- \fo\left(\x_{i j}\right)\right)^{2}}{\cf  +r}\right|\right],
\end{align*}
which could be further bounded by
$$
\begin{aligned}
 & \frac{1}{r}\mathbb{E}\left[\sup _{f \in \mathcal{F}(r)}\left|\frac{1}{n m} \sum_{i=1}^{n} \sum_{j=1}^{m} \sigma_{i j}\left(f\left(\x_{i j}\right)- \fo\left(\x_{i j}\right)\right)^{2}\right|\right]\\
&+\sum_{k=1}^{\infty} \frac{1}{(4^{k-1}+1)r} \mathbb{E}\left[\sup _{f \in \mathcal{F}(4^{k}r) }\left|\frac{1}{n m} \sum_{i=1}^{n} \sum_{j=1}^{m} \sigma_{i j}\left(f\left(\x_{i j}\right)- \fo\left(\x_{i j}\right)\right)^{2} \right|\right].\\
\end{aligned}
$$
The last inequality stems from the definition of $\mathcal{F}(\cdot)$. 
As a reminder, the functions in $\mathcal{F}$ are bounded in infinity norm. 
By invoking Lemma \ref{lem:ltc}, we can further upper-bound the Rademacher complexity of  $\left\{ (f-\fo)^2: f\in \mathcal{F}(r') \right\}$ with
$$
\begin{aligned}
&\mathbb{E}\left[ \sup _{f \in \mathcal{F}(r')}\left|\frac{1}{n m} \sum_{i=1}^{n} \sum_{j=1}^{m} \sigma_{i j}\left(f\left(\x_{i j}\right)- \fo\left(\x_{i j}\right)\right)^{2}\right|
\right]\\
\leq & 4\jdb{b1} \mathbb{E}\left[ \sup _{f \in \mathcal{F}(r')}\left|\frac{1}{n m} \sum_{i=1}^{n} \sum_{j=1}^{m} \sigma_{i j}\left(f\left(\x_{i j}\right)- \fo\left(\x_{i j}\right)\right)\right|
\right]\\
\leq & 4 \jdb{b1} \phi(r'),
\end{aligned}
$$
for any $r'\geq r^{*}$. 
Hence, by substituting this bound into the last line of the earlier formula, we obtain
$$
\begin{aligned}
\frac{4 \jdb{b1} \phi(r)}{r}+\sum_{k=1}^{\infty} \frac{4 \jdb{b1} \phi(4^k r)}{(4^{k-1}+1)r} 
\leq \frac{4 \jdb{b1} \phi(r)}{r}+\sum_{k=1}^{\infty} \frac{2^{k+2} \jdb{b1} \phi( r)}{(4^{k-1}+1)r} 
\leq \frac{16 \jdb{b1} \phi(r)}{r},
\end{aligned}
$$
where the first inequality is because that $\phi$ is sub-root and the second inequality is derived by simple algebra.
Combining the above several results yields
\begin{equation}\label{equ:101}
\begin{aligned}
&\mathbb{E}\left[\sup _{f \in \mathcal{F}}\left|\frac{1}{n m} \sum_{i=1}^{n} \sum_{j=1}^{m} \frac{\cf  -\left(f\left(\x_{i j}\right)- \fo\left(\x_{i j}\right)\right)^{2}}{\cf  +r}\right|\right]\leq \frac{32 \jdb{b1} \phi(r)}{r}.
\end{aligned}    
\end{equation}
Then, to apply Talagrand's inequality, we shall verify the remaining conditions. It is easy to show 
\begin{equation}\label{equ:102}
    \mathbb{E}\left[ \frac{\cf  -\left(f\left(\x\right)- \fo\left(\x\right)\right)^{2}}{\cf  +r} \right] = 0,
\end{equation}
\begin{equation}\label{equ:103}
    \left|\frac{\cf  -\left(f\left(\x\right)- \fo\left(\x\right)\right)^{2}}{\cf  +r} \right|\leq \frac{4\jdb{b1}}{r},
\end{equation}
 and 
\begin{equation}\label{equ:104}
\begin{aligned}
 \mathbb{E}\left[ \left|\frac{\cf  -\left(f\left(\x\right)- \fo\left(\x\right)\right)^{2}}{\cf  +r} \right|^2 \right]
\leq \frac{\mathbb{E}\left[ \left(f\left(\x\right)- \fo\left(\x\right)\right)^{4} \right]}{4r \cf }
\leq  \frac{\jdb{b1}^2}{r}.
\end{aligned}
\end{equation}
The first inequality of \eqref{equ:104} is because both that the expectation of the numerator is zero  
and using AM–GM inequality for the denominator.
According, we can use Talagrand's concentration inequality now. By Lemma \ref{lem:tal} and \eqref{equ:101}, \eqref{equ:102}, \eqref{equ:103}, \eqref{equ:104}, with probability at least $1-\exp\{-t\}$, it holds 
\begin{equation}\label{equ:105}
\begin{aligned}
&\frac{1}{n m} \sum_{i=1}^{n} \sum_{j=1}^{m} \frac{\chf  -\left(\hat f\left(\x_{i j}\right)- \fo\left(\x_{i j}\right)\right)^{2}}{\chf  +r} \\
\leq &
\sup_{f\in \mathcal{F}}
\left| 
\frac{1}{n m} \sum_{i=1}^{n} \sum_{j=1}^{m} \frac{\cf  -\left(f\left(\x_{i j}\right)- \fo\left(\x_{i j}\right)\right)^{2}}{\cf  +r} 
\right|\\
\leq & 
\frac{64 \jdb{b1} \phi(r)}{r}+ \sqrt{\frac{8 \jdb{b1}^2 t}{nmr}} +\frac{140 \jdb{b1}^2 t}{nmr}.
\end{aligned}    
\end{equation}
Therefore, there are universal constants $\jc{c801}$ and $\jc{c802}$ such that $r_{1}(t)=2\jc{c801}\jdb{b1}^2r^{*}+\frac{2\jc{c802}\jdb{b1}^2t}{nm}$ satisfies
$$
\frac{64 \jdb{b1} \phi(r_{1}(t))}{r_{1}(t)}
+ \sqrt{\frac{8 \jdb{b1}^2 t}{nmr_{1}(t)}} +\frac{140 \jdb{b1}^2 t}{nmr_{1}(t)}
\leq \frac{64 \jdb{b1} \phi(2\jc{c801}\jdb{b1}^2r^{*})}{2\jc{c801}\jdb{b1}^2r^{*}} + \sqrt{\frac{4}{\jc{c802}}} +\frac{70 }{\jc{c802}}
\leq \frac{1}{2}.
$$
Substituting $r=r_{1}(t)$ into \eqref{equ:105}, the desired result follows. 
\end{proof}

\begin{proof}[Proof of Lemma \ref{lem:lema2}]
We use a similar method with Lemma \ref{lem:lema1} to consider 
$$
\left|
\frac{
\left(\mathcal{T}_{\beta_{\epsilon}} \epsilon_{ij}-\mathbb{E}\left[ \mathcal{T}_{\beta_{\epsilon}} \epsilon_{ij} \right]\right)\left(\hat{f}\left(\x_{i j}\right)-\fo\left(\x_{i j}\right)\right)
}{\chf  +r}\right|,
$$
which is bounded by
$$
\sup_{f\in \mathcal{F}}
\left| 
 \frac{
\left(\mathcal{T}_{\beta_{\epsilon}} \epsilon_{ij}-\mathbb{E}\left[ \mathcal{T}_{\beta_{\epsilon}} \epsilon_{ij} \right]\right)\left(f\left(\x_{i j}\right)-\fo\left(\x_{i j}\right)\right)
}{\cf  +r}
\right|.
$$
After some discussion similar to the proof of Lemma \ref{lem:lema1}, we obtain,  
with probability at least $1-\exp\{-t\}$,
$$
\begin{aligned}
& \left| 
\frac{1}{n m} \sum_{i=1}^{n} \sum_{j=1}^{m} \frac{
\left(\mathcal{T}_{\beta_{\epsilon}} \epsilon_{ij}-\mathbb{E}\left[ \mathcal{T}_{\beta_{\epsilon}} \epsilon_{ij} \right]\right)\left(\hat{f}\left(\x_{i j}\right)-\fo\left(\x_{i j}\right)\right)
}{\chf  +r} 
\right|\\
\leq & \sup_{f\in \mathcal{F}}
\left| 
\frac{1}{n m} \sum_{i=1}^{n} \sum_{j=1}^{m} \frac{
\left(\mathcal{T}_{\beta_{\epsilon}} \epsilon_{ij}-\mathbb{E}\left[ \mathcal{T}_{\beta_{\epsilon}} \epsilon_{ij} \right]\right)\left(\hat{f}\left(\x_{i j}\right)-\fo\left(\x_{i j}\right)\right)
}{\chf  +r} 
\right|\\
\leq & \frac{32\beta_{\epsilon} \phi(r)}{r}+\sqrt{ \frac{2 \jdb{b3}^{2}  t}{nmr} }+ \frac{140\jdb{b1}\beta_{\epsilon} t}{nmr}.
\end{aligned}
$$
The second term on the right side is different from its counterpart in step 3 of the proof of Theorem \ref{thm:thm1} because the variance estimation is smaller due to the independent of $X_{ij}$ and $\epsilon_{ij}$. 
Choosing constants $\xc\label{c111}$ and $\xc\label{c112}$ such that $r_{3}(t)=\jc{c111}\beta_{\epsilon}^2r^{*}+\frac{\jc{c112}(\jdb{b1}\beta_{\epsilon}+\jdb{b3}^{2} ) t }{nm}$ satisfy
$$
 \frac{32\beta_{\epsilon} \phi(r)}{r}+\sqrt{ \frac{2 \jdb{b3}^{2}  t}{nmr} }+ \frac{140\jdb{b1}\beta_{\epsilon} t}{nmr} \leq \frac{1}{16}.
$$
Therefore, with probability at least $1-\exp\{-t\}$,
$$
\begin{aligned}
& \left| 
\frac{1}{n m} \sum_{i=1}^{n} \sum_{j=1}^{m} \left(\mathcal{T}_{\beta_{\epsilon}} \epsilon_{ij}-\mathbb{E}\left[ \mathcal{T}_{\beta_{\epsilon}} \epsilon_{ij} \right]\right)\left(\hat{f}\left(\x_{i j}\right)-\fo\left(\x_{i j}\right)\right)
\right|\\
\leq &\frac{1}{16}\chf  +\frac{\jc{c111}\beta_{\epsilon}^2r^{*}}{16}+\frac{\jc{c112}(\jdb{b1}\beta_{\epsilon}+ \jdb{b3}^{2}  ) t }{16nm}.
\end{aligned}
$$
This completes the proof. 
\end{proof}

\begin{proof}[Proof of Lemma \ref{lem:bloc}]
The main proof of Step 1 can be found in \citet{lu2021deep}. 

Step 1. We first consider the case $K=L^2 W^2$.
Applying Lemma  2.2 in \citet{SHEN201974} or Lemma 5.2 in \citet{lu2021deep} yields that there is a ReLU neural network
$$
\check{f}_{1} \in \nn\left(\text{inputdim}=1, \text{widthvec}=  [2W; 4WL-1]  \right)
$$
satisfying
$$
\check{f}_{1}\left(\frac{k}{L W^2}\right)=\check{f}_{1}\left(\frac{k+1}{L W^2}-\delta\right)=k
$$
for $k=0,1,...,L W^2-2$, 
$$
\check{f}_{1}\left(\frac{L W^2-1}{L W^2}\right)=\check{f}_{1}\left(1\right)=L W^2-1,\quad \check{f}_{1}\left(2\right)=0,
$$
and being linear on each interval between the nodes above.
These imply that 
$$
\check{f}_{1}\left(x\right)=k, \quad \text {when } x \in\left[\frac{k}{L W^2}, \frac{k+1}{L W^2}-\delta \cdot I_{\{k \leq L W^2-2\}}\right],
$$
for $k=0,1,...,L W^2-2$, and 
$$
\check{f}_{1}\left(x\right)=L W^2-1, \quad \text {when } x \in\left[\frac{L W^2-1}{L W^2}, 1\right].
$$
To represent the shallow and wide neural network $\check{f}_{1}$ with a deep and narrow one, we use Proposition 2.2 in \citet{shen2020} and get there is a neural network 
$$
\check{f}_{2}\in \nn\left( 1, 2L+1, 4W+2 \right)
$$
such that  $\check{f}_{2}(\xx)$ is 
always equal to $\check{f}_{1}(\xx)$.

Reconsider the construction of $\check{f}_{2}(\xx)$ and set $W=1$, we know there is a neural network 
$$
\check{f}_{3}\in \nn\left( 1, 2L+1, 6 \right)
$$
such that
$$
\check{f}_{3}\left(x\right)=k, \quad \text {when } x \in\left[\frac{k}{L}, \frac{k+1}{L}-\delta \cdot I_{\{k \leq L-2\}}\right]
$$
for $k=0,1,..., L-2$, and 
$$
\check{f}_{3}\left(x\right)= L-1, \quad \text {when } x \in\left[\frac{ L-1}{L}, 1\right].
$$

After constructing $\check{f}_{2}$ and $\check{f}_{3}$, we analyze the original problem. Let $$
\check{f}_{4}=L\check{f}_{2}(\xx)+\check{f}_{3}\left(L W^2(x-\check{f}_{2}(\xx))\right), 
$$
it is easy to verify 
$$
\check{f}_{4}(\xx)=k, \quad \text {when } x \in\left[\frac{k}{L^2 W^2}, \frac{k+1}{L^2 W^2}-\delta \cdot I_{\{k \leq L^2 W^2-2\}}\right] ,
$$
for $k=0,1, \cdots, L^2 W^2-1$.
And now we should only specify the width and depth of the desired network. To this goal, we decompose the function $\check{f}_{4}$. The first composited network inputs $x$ and outputs $(\check{f}_{2}(x),x)$. Then the second composited network inputs $(\check{f}_{2}(x),x)$ and outputs $(\check{f}_{2}(x),L^2 W(x-\check{f}_{2}(x))$. Using the results of the previous layer, the last composited network calculates $\check{f}_{4}=L\check{f}_{2}(x)+\check{f}_{3}\left(L^2 W(x-\check{f}_{2}(x))\right)$. This yields the width of $\check{f}_{4}$ is less than $ 4W+2+1 \leq 4W+3 $ and depth $\leq 2L+1+2+2L+1+2= 4L+6$. This completes the proof for $K=L^2 W^2$.

Step 2. 
When $K<L^2 W^2$, Let 
$$
\check{f}_{5}(x)= \check{f}_{4}'\left(\frac{Kx}{L^2 W^2}\right).
$$
The construction of  $ \check{f}_{4}'$ uses the previous step, 
and the only difference is that $\delta'$ is taken as $\frac{K\delta}{L^2 W^2}$. 
One can check that $\check{f}_{5}$ satisfies the requirements. The desired result follows.
\end{proof}

\begin{proof}[Proof of Lemma \ref{lem:bspl}]

Step 1.
We first consider the case $d=1$. 
Applying Lemma 5.3 in \citet{lu2021deep}, we know there is a ReLU neural network
$$
\check{f}_{1} \in \nn\left(1, 7(r-1)(r-2)L, 9W+r+7  \right)
$$
satisfying
$$
\left| \check{f}_{1}(x)-x^{r-1} \right| \leq 9 (r-2)(W+1)^{-7(r-1)L},
$$
for $x\in[0,1]$. Let $\check{f}_{2}=r^{r-1}\check{f}_{1}\left(\frac{x}{r}\right)$. Then 
$$
\left| \check{f}_{2}(x)-x^{r-1} \right| \leq 9 (r-2)r^{r-1}(W+1)^{-7(r-1)L},
$$
for $x\in[0,r]$. Using $\check{f}_{2}$ and recalling $\eqref{equ:301}$, 
let  
$$
\check{f}_{3}(x)=\frac{1}{(r-1) !} \sum_{i=0}^{r}(-1)^{i} \binom{r}{i}\check{f}_{2} \circ \operatorname{ReLU}(x-i).
$$
In order to make the composition in the next step more convenient, we truncate the function 
$$
\check{f}_{4}(x)=
\operatorname{ReLU} \left(\check{f}_{3}(x) \right) - \operatorname{ReLU} \left(\check{f}_{3}(x) -1\right).
$$
Then the approximation error is bounded by
$$
\begin{aligned}
\left|\check{f}_{4}(x)- B_{r}(x) \right|\leq & \left|\check{f}_{3}(x)- B_{r}(x) \right|\\
\leq & \frac{1}{(r-1) !} \sum_{i=0}^{r} \binom{r}{i}\left|\check{f}_{2} \circ \operatorname{ReLU}(x-i) - (x-i)_{+}^{r-1} \right|\\
< & \frac{9\cdot 2^{r}r^{r}}{(r-1) !} (W+1)^{-7(r-1)L}.
\end{aligned}
$$

Step 2. 
Now we focus on the multivariate case for $\mathcal{B}^{d}_{\bk,\bo}(x)$ with 
$x=(x_{1}, x_{2}, ..., x_{d})$. By Lemma 5.3 in \citet{lu2021deep}, we know there is a ReLU neural network
$$
\check{f}_{5} \in \nn\left(d, 7d(d-1)L, 9W+d+8  \right)
$$
satisfying 
$$
\left|\check{f}_{5}(x_{1}, x_{2}, ..., x_{d})-x_{1}x_{2}\cdots x_{d}\right|\leq 9(d-1)(N+1)^{-7dL},
$$
for $x_{i}\in [0,1]$. 
Let 
$$
\check{f}_{6}(x)=\check{f}_{5}\left( \check{f}_{4}(x_{1}), \check{f}_{4}(x_{2}), ... ,\check{f}_{4}(x_{d})  \right).
$$
Hence, we have
$$
\begin{aligned}
& \left|\check{f}_{6}(x)- \mathcal{B}^{d}_{1,\bo}(x) \right|\\
= & \left|\check{f}_{5}\left( \check{f}_{4}(x_{1}), \check{f}_{4}(x_{2}), ... ,\check{f}_{4}(x_{d})  \right) - \prod_{i=1}^{d} B_{r}(  x_{i}) \right|\\
\leq & \left|\check{f}_{5}\left( \check{f}_{4}(x_{1}), \check{f}_{4}(x_{2}), ... ,\check{f}_{4}(x_{d})  \right) - \prod_{i=1}^d \check{f}_{4}(x_{i}) \right|  
+ \left| \prod_{i=1}^d \check{f}_{4}(x_{i}) - \prod_{i=1}^d B_{r}(x_{i}) \right| .
\end{aligned}
$$
By the construction of $\check{f}_{5}$, the first term is bounded by $9(d-1)(W+1)^{-7dL}$. Then the second term is bounded by
$$
\begin{aligned}
& \left| \prod_{i=1}^d \check{f}_{4}(x_{i} ) - \prod_{i=1}^d B_{r}(x_{i} ) \right|\\
\leq & \sum_{j=0}^{d-1} \left| \prod_{i=1}^{d-j} \check{f}_{4}(x_{i} )\prod_{i=d-j+1}^d B_{r}( x_{i} ) - \prod_{i=1}^{d-j-1} \check{f}_{4}(x_{i} )\prod_{i=d-j}^d B_{r}(  x_{i} ) \right|\\
\leq & \sum_{j=0}^{d-1} \left| \check{f}_{4}(x_{d-j})-  B_{r}(x_{d-j})  \right|\\
\leq & \frac{9\cdot 2^{r}r^{r}d}{(r-1) !} (W+1)^{-7(r-1)L}.
\end{aligned}
$$
The first inequality is given by triangle inequality, the second inequality is due to $|\check{f}_{4}|, |B_{r}|\leq 1$, and the last inequality is implied from the approximation error of $\check{f}_{4}$.
Hence, it follows that the approximation error satisfies the requirement. 
Finally, it is easy to check the size of neural network $\check{f}_{6}$. This completes the proof. 
\end{proof}

\begin{proof}[Proof of Lemma \ref{lem:lowerlinearlemma}]
The proof is divided into two parts, each proving a term on the right-hand side. The proof of the $nm$ terms comes from Theorem 1 in \citet{zhang2002wavelet}, and the $n$ terms are the same as the other lower bound in this paper. 
To prove the $nm$ term, we assume $Z$ is a Gaussian process and $\epsilon$ is a Gaussian noise with variance $\jdb{b3}$. 
For simplicity, we write
$$
\varphi_{ij}\left(x\right)=\varphi_{ij}\left(x, X_{11}, \ldots, X_{nm}\right).
$$
By \citet{zhang2002wavelet},  we know that $P(\mathcal{D})\geq 1-o(1)$. Denote $I_{\mathcal{D}}$ as indicator function of event $\mathcal{D}$. 
Recalling the definition of linear estimator $\tilde{f}_{L}$, we have
\begin{equation}\label{equ:681}
\begin{aligned}
\mathcal{E}_{\fo}(\tilde{f}_{L})=&\mathbb{E}\left[ \left\| \sum_{i=1}^n\sum_{j=1}^m Y_{ij} \varphi_{ij}\left(\xx\right) - \fo(\xx) \right\|_{L^2(\Omega)}^{2}\right]\\
\geq & \mathbb{E}\left[ \left\| \sum_{i=1}^n\sum_{j=1}^m \fo(\x_{ij}) \varphi_{ij}\left(x\right) - \fo(\xx) \right\|_{L^2(\Omega)}^{2}I_{\mathcal{D}}\right] \\ & + \jdb{b3} \sum_{i=1}^{n} \sum_{j=1}^m \mathbb{E}\left[ \left\| \varphi_{ij}(\xx) \right\|_{L^2(\Omega)}^{2} I_{\mathcal{D}}\right]. \\
\end{aligned}
\end{equation}
Because $\mathcal{A}$ is a partition of $\Omega$ and $|\mathcal{A}|\geq \xc\label{c601}^{-1} 2^K$ for some constant $\jc{c601}$, there exists $A\in\mathcal{A}$ such that 
\begin{equation}\label{equ:682}
\sum_{i=1}^{n} \sum_{j=1}^m \mathbb{E}\left[ \left\| \varphi_{ij}(\xx) \right\|_{L^2(A)}^{2} I_{\mathcal{D}} \right] \leq \jc{c601} 2^{-K} \jdb{b3}^{-1} \mathcal{E}_{\fo}(\tilde{f}_{L}). 
\end{equation}
Fix $A$ and use condition (i), there is a function $g \in \mathcal{F}^{\circ}$ satisfying $g(\xx) \geq \frac{1}{2}\jb{xb:delta}\jb{xb:bf}$ for all $\xx \in A$, which further implies
\begin{equation}\label{equ:683}
\mathbb{E}\left[ \left\| g(\xx) \right\|_{L^2(A)}I_{\mathcal{D}}\right]  \geq  \frac{1}{2}\jb{xb:delta}\jb{xb:bf} 2^{-K/2} P(\mathcal{D}).
\end{equation}
Denote the random function
$$
h(\xx) = \sum_{i=1}^n\sum_{j=1}^m g(\x_{ij}) \varphi_{ij}\left(x\right). 
$$
Then, we have 
$$
\begin{aligned}
\mathbb{E}\left[ \left\| h(\xx) \right\|_{L^2(\Omega)}^{2} I_{\mathcal{D}} \right] = & \mathbb{E}\left[  \int_{\Omega}\left( \sum_{i=1}^n\sum_{j=1}^m g(\x_{ij}) \varphi_{ij}\left(x\right) \right)^2 d\xx \cdot I_{\mathcal{D}} \right] \\ 
\leq & \mathbb{E}\left[ \left( \sum_{i=1}^n\sum_{j=1}^m g^2 (\x_{ij}) \right)   \int_{\Omega}  \sum_{i=1}^n\sum_{j=1}^m \varphi_{ij}^2\left(x\right)  d\xx \cdot I_{\mathcal{D}} \right] \\
\leq & \jc{c601} 2^{-(K^{\prime}+K)} nm \jb{xb:bcpp}\jb{xb:delta}^2  \jdb{b3}^{-1} \mathcal{E}_{\fo}(\tilde{f}_{L}), 
\end{aligned}
$$
where the first inequality is given by Cauchy inequality and the second inequality is due to \eqref{equ:682} and condition (ii).  
This implies 
\begin{equation}\label{equ:684}
\mathbb{E}\left[ \left\| h(\xx) \right\|_{L^2(\Omega)} I_{\mathcal{D}} \right] \leq \sqrt{\mathbb{E}\left[ \left\| h(\xx) \right\|_{L^2(\Omega)}^{2} I_{\mathcal{D}} \right]} \leq \sqrt{ \jc{c601} 2^{-(K^{\prime}+K)} nm \jb{xb:bcpp}\jb{xb:delta}^2  \jdb{b3}^{-1} \mathcal{E}_{\fo}(\tilde{f}_{L}) }.
\end{equation}
$$
\begin{aligned}
\sup_{\mathcal{P}_{Z, \epsilon} \in \mathcal{Q}_{Z,\epsilon}} \mathcal{E}_{\fo}(\tilde{f}_{L})  \geq & \mathbb{E}\left[ \left\| h\left(x\right) - g(\xx) \right\|_{L^2(\Omega)}^{2}I_{\mathcal{D}}\right]\\
\geq & \mathbb{E}\left[ \left\| h\left(x\right) - g(\xx) \right\|_{L^2(A)}^{2}I_{\mathcal{D}}\right]\\
\geq & \mathbb{E}\left[ \left(
\left\| g(\xx) \right\|_{L^2(A)}I_{\mathcal{D}} - \left\| h\left(x\right) \right\|_{L^2(A)}I_{\mathcal{D}}
\right)^{2}\right]\\
\geq &  \left( \mathbb{E}\left[ \left\| g(\xx) \right\|_{L^2(A)}I_{\mathcal{D}}\right] - 
\mathbb{E}\left[\left\| h\left(x\right) \right\|_{L^2(A)}I_{\mathcal{D}}\right]
\right)^{2}\\
\end{aligned}
$$
where the first inequality is implied from \eqref{equ:681}, the third inequality is due to the Triangle inequality, and the last inequality is given by Jensen inequality. Recalling \eqref{equ:683} and \eqref{equ:684}, we have 
\begin{equation}\label{equ:689}
\begin{aligned}
 & \mathbb{E}\left[ \left\| g(\xx) \right\|_{L^2(A)}I_{\mathcal{D}}\right] - \mathbb{E}\left[\left\| h\left(x\right) \right\|_{L^2(A)}I_{\mathcal{D}}\right] \\&\quad\quad\quad\geq \frac{1}{2}\jb{xb:delta}\jb{xb:bf} 2^{-K/2} P(\mathcal{D}) - \sqrt{ \jc{c601} 2^{-(K^{\prime}+K)} nm \jb{xb:bcpp}\jb{xb:delta}^2  \jdb{b3}^{-1} \mathcal{E}_{\fo}(\tilde{f}_{L})} . 
\end{aligned}   
\end{equation}
If the right side of \eqref{equ:689} $\leq \frac{1}{4}\jb{xb:delta}\jb{xb:bf} 2^{-K/2} P(\mathcal{D})$, i.e. $$\frac{1}{4}\jb{xb:delta}\jb{xb:bf} 2^{-K/2} P(\mathcal{D}) \leq  \sqrt{ \jc{c601} 2^{-(K^{\prime}+K)} nm \jb{xb:bcpp}\jb{xb:delta}^2  \jdb{b3}^{-1} \mathcal{E}_{\fo}(\tilde{f}_{L})},$$ we have 
\begin{equation}\label{equ:685}
\mathcal{E}_{\fo}(\tilde{f}_{L}) \geq \frac{\jb{xb:bf}^2 2^{K'} \jdb{b3}P(\mathcal{D})}{16\jc{c601} nm \jb{xb:bcpp} } . 
\end{equation}
If the right side of \eqref{equ:689} $\geq \frac{1}{4}\jb{xb:delta}\jb{xb:bf} 2^{-K/2} P(\mathcal{D})>0$, we have 
\begin{equation}\label{equ:686}
\begin{aligned}
\sup_{\mathcal{P}_{Z, \epsilon} \in \mathcal{Q}_{Z,\epsilon}} \mathcal{E}_{\fo}(\tilde{f}_{L})  
\geq &  \left( \mathbb{E}\left[ \left\| g(\xx) \right\|_{L^2(A)}I_{\mathcal{D}}\right] - 
\mathbb{E}\left[\left\| h\left(x\right) \right\|_{L^2(A)}I_{\mathcal{D}}\right]
\right)^{2}\\
\geq & \frac{1}{16}\jb{xb:delta}^2 \jb{xb:bf}^2  2^{-K} P(\mathcal{D})^2 . 
\end{aligned}  
\end{equation}
Combining the discussion in \eqref{equ:685} and \eqref{equ:686} and noting that $P(\mathcal{D}) \geq 1-o(1)$, we have 
\begin{equation}\label{equ:687}
    \inf_{\tilde{f}_{L}}\sup_{\mathcal{P}_{Z,\epsilon} \in \mathcal{Q}_{Z,\epsilon}} \mathcal{E}_{\fo}(\tilde{f}_{L})  
\gtrsim \min\left( \frac{\jb{xb:bf}^2 2^{K'}}{\jb{xb:bcpp}nm} +  \jb{xb:delta}^2 \jb{xb:bf}^2 2^{-K} \right). 
\end{equation}
Moreover, to show the remaining $1/n$ term in the lower bound,  it now suffices to consider $Z_{i}$ as constant functions over $\xx \in \Omega$ and the measurement noise $\epsilon_{ij}$ as zero. Assume that $\var(Z_{i}(\x)|\x)=\jdb{b2}$. Then the problem becomes the mean estimation from $n$ i.i.d. observations. It is well known that
the rate 
\begin{equation}\label{equ:688}
    \inf_{\tilde{f}}\sup_{\mathcal{P}_{Z,\epsilon} \in \mathcal{Q}_{Z,\epsilon}} \mathcal{E}_{\fo}(\tilde{f}_{L}) 
\geq \jdb{b2}n^{-1}.
\end{equation}
Then combining \eqref{equ:687} and \eqref{equ:688} completes the proof. 

\end{proof}

\bibliographystyle{agsm}
\bibliography{Bibliography-MM-MC}

@article{schmidtaos,
author = {Johannes Schmidt-Hieber},
title = {{Nonparametric regression using deep neural networks with ReLU activation function}},
volume = {48},
journal = {The Annals of Statistics},
number = {4},
publisher = {Institute of Mathematical Statistics},
pages = {1875 -- 1897},
year = {2020},
URL = {https://doi.org/10.1214/19-AOS1875}
}

@article{lu2021deep,
  title={Deep network approximation for smooth functions},
  author={Lu, Jianfeng and Shen, Zuowei and Yang, Haizhao and Zhang, Shijun},
  journal={SIAM Journal on Mathematical Analysis},
  volume={53},
  number={5},
  pages={5465--5506},
  year={2021},
  publisher={SIAM}
}

@book{ginebook, 
place={Cambridge}, 
series={Cambridge Series in Statistical and Probabilistic Mathematics}, 
title={Mathematical Foundations of Infinite-Dimensional Statistical Models}, 
publisher={Cambridge University Press}, 
address={New York},
author={Giné, Evarist and Nickl, Richard}, 
year={2015}, 
collection={Cambridge Series in Statistical and Probabilistic Mathematics},
edition=1,
DOI={10.1017/CBO9781107337862}, 
}

@article{pascalta,
author = {Pascal Massart},
title = {{About the constants in Talagrand's concentration inequalities for empirical processes}},
volume = {28},
journal = {The Annals of Probability},
number = {2},
publisher = {Institute of Mathematical Statistics},
pages = {863 -- 884},
keywords = {Concentration inequalities, concentration of measure, Deviation inequalities, Empirical processes},
year = {2000},
//DOI = {10.1214/aop/1019160263},
URL = {https://doi.org/10.1214/aop/1019160263}
}

@book{ledoux1991probability,
  title={Probability in Banach Spaces: Isoperimetry and Processes},
  author={Ledoux, Michel and Talagrand, Michel},
  volume={23},
  year={1991},
  publisher={Springer Science \& Business Media},
  address={Berlin, Heidelberg},
  edition=1
}

@book{talagrand2014process,
  title={Upper and Lower Bounds for Stachastic Processes},
  author={Talagrand, Michel},
  year={2014},
  publisher={Springer Berlin, Heidelberg},
  edition=1
}

@article{bartlett2005local,
  title={Local Rademacher complexities},
  author={Bartlett, Peter L and Bousquet, Olivier and Mendelson, Shahar},
  journal={The Annals of Statistics},
  volume={33},
  number={4},
  pages={1497--1537},
  year={2005},
  publisher={Institute of Mathematical Statistics},
  //mrnumber={MR2166554},
  //DOI={10.1214/009053605000000282}
}

@article{kohler2021rate,
  title={On the rate of convergence of fully connected deep neural network regression estimates},
  author={Kohler, Michael and Langer, Sophie},
  journal={The Annals of Statistics},
  volume={49},
  number={4},
  pages={2231--2249},
  year={2021},
  publisher={Institute of Mathematical Statistics}
}

@article{huangjian,
doi = {10.48550/ARXIV.2104.06708},
author = {Yuling Jiao and Guohao Shen and Yuanyuan Lin and Jian Huang},
title = {{Deep nonparametric regression on approximate manifolds: Nonasymptotic error bounds with polynomial prefactors}},
volume = {51},
journal = {The Annals of Statistics},
number = {2},
publisher = {Institute of Mathematical Statistics},
pages = {691 -- 716},
year = {2023},
doi = {10.1214/23-AOS2266},
URL = {https://doi.org/10.1214/23-AOS2266}
}

@article{SHEN201974,
title = {Nonlinear approximation via compositions},
journal = {Neural Networks},
volume = {119},
pages = {74-84},
year = {2019},
issn = {0893-6080},
author = {Zuowei Shen and Haizhao Yang and Shijun Zhang},
//DOI = {https://doi.org/10.1016/j.neunet.2019.07.011},
url = {https://www.sciencedirect.com/science/article/pii/S0893608019301996},
keywords = {Deep neural networks, ReLU activation function, Nonlinear approximation, Function composition, Hölder continuity, Parallel computing},
abstract = {Given a function dictionary D and an approximation budget N∈N, nonlinear approximation seeks the linear combination of the best N terms {Tn}1≤n≤N⊆D to approximate a given function f with the minimum approximation error εL,f≔min{gn}⊆R,{Tn}⊆D‖f(x)−∑n=1NgnTn(x)‖.Motivated by recent success of deep learning, we propose dictionaries with functions in a form of compositions, i.e., T(x)=T(L)∘T(L−1)∘⋯∘T(1)(x)for all T∈D, and implement T using ReLU feed-forward neural networks (FNNs) with L hidden layers. We further quantify the improvement of the best N-term approximation rate in terms of N when L is increased from 1 to 2 or 3 to show the power of compositions. In the case when L>3, our analysis shows that increasing L cannot improve the approximation rate in terms of N. In particular, for any function f on [0,1], regardless of its smoothness and even the continuity, if f can be approximated using a dictionary when L=1 with the best N-term approximation rate εL,f=O(N−η), we show that dictionaries with L=2 can improve the best N-term approximation rate to εL,f=O(N−2η). We also show that for Hölder continuous functions of order α on [0,1]d, the application of a dictionary with L=3 in nonlinear approximation can achieve an essentially tight best N-term approximation rate εL,f=O(N−2α∕d). Finally, we show that dictionaries consisting of wide FNNs with a few hidden layers are more attractive in terms of computational efficiency than dictionaries with narrow and very deep FNNs for approximating Hölder continuous functions if the number of computer cores is larger than N in parallel computing.}
}

@article{suzuki2018adaptivity,
  title={Adaptivity of deep ReLU network for learning in Besov and mixed smooth Besov spaces: optimal rate and curse of dimensionality},
  author={Suzuki, Taiji},
  journal={ArXiv preprint. arXiv:1810.08033},
  year={2018}
}

@article{LeisnerPhD,
 ISSN = {00222518, 19435258},
 URL = {http://www.jstor.org/stable/24902859},
 author = {Christopher Leisner},
 journal = {Indiana University Mathematics Journal},
 number = {2},
 pages = {437--455},
 publisher = {Indiana University Mathematics Department},
 title = {Nonlinear Wavelet Approximation in Anisotropic Besov Spaces},
 urldate = {2024-11-11},
 volume = {52},
 year = {2003}
}

@article{lin2000nonparametric,
  title={Nonparametric function estimation for clustered data when the predictor is measured without/with error},
  author={Lin, Xihong and Carroll, Raymond J},
  journal={Journal of the American statistical Association},
  volume={95},
  number={450},
  pages={520--534},
  year={2000},
  publisher={Taylor \& Francis}
}

@article{cai2011optimal,
  title={Optimal estimation of the mean function based on discretely sampled functional data: Phase transition},
  author={Cai, T Tony and Yuan, Ming},
  journal={The Annals of Statistics},
  volume={39},
  number={5},
  pages={2330--2355},
  year={2011},
  //mrnumber={MR2906870},
  publisher={Institute of Mathematical Statistics},
  //DOI={10.1214/11-AOS898}
}

@article{zhang2016sparse,
  title={From sparse to dense functional data and beyond},
  author={Zhang, Xiaoke and Wang, Jane-Ling},
  journal={The Annals of Statistics},
  volume={44},
  number={5},
  pages={2281--2321},
  year={2016},
  publisher={Institute of Mathematical Statistics}
}

@book{ramsay2005, 
place={Springer}, 
series={Springer Series in Statistics}, 
title={Functional Data Analysis}, 
//DOI={10.1017/CBO9781107337862}, 
publisher={Springer New York}, 
author={Ramsay, Jim and Silverman, Bernard}, 
year={2005}, 
edition={1}
}

@book{diggle2002analysis,
  title={Analysis of Longitudinal Data},
  author={Diggle, Peter and Diggle, Peter J and Heagerty, Patrick and Liang, Kung-Yee and Zeger, Scott and others},
  year={2002},
  publisher={Oxford university press},
  edition = 2,
  address={Oxford}
}

@book{hsiao2022analysis,
  title={Analysis of Panel Data},
  author={Hsiao, Cheng},
  year={2022},
  publisher={Cambridge university press},
  address={Cambridge},
  edition = 3,
}

@inproceedings{telgarsky2016benefits,
  title={Benefits of depth in neural networks},
  author={Telgarsky, Matus},
  booktitle={Conference on learning theory},
  pages={1517--1539},
  year={2016},
  organization={PMLR}
}

@article{yarotsky2017error,
  title={Error bounds for approximations with deep {ReLU} networks},
  author={Yarotsky, Dmitry},
  journal={Neural Networks},
Year = {2017},
Volume = {94},
Pages = {103-114},
Month = {OCT},
DOI = {10.1016/j.neunet.2017.07.002},
ISSN = {0893-6080},
EISSN = {1879-2782},
ResearcherID-Numbers = {Yarotsky, Dmitry/A-8811-2016},
ORCID-Numbers = {Yarotsky, Dmitry/0000-0002-5432-7143},
Unique-ID = {WOS:000410973300011},
}

@inproceedings{yarotsky2018optimal,
  title={Optimal approximation of continuous functions by very deep ReLU networks},
  author={Yarotsky, Dmitry},
  booktitle={Conference on learning theory},
  pages={639--649},
  year={2018},
  organization={PMLR}
}

@article{Shen2020,
	//DOI = {10.4208/cicp.oa-2020-0149},
	url = {https://doi.org/10.4208\%2Fcicp.oa-2020-0149},
	year = 2020,
	month = {jun},
	publisher = {Global Science Press},
	volume = {28},
	number = {5},
	pages = {1768--1811},
	author = {Zuowei Shen},
	title = {Deep Network Approximation Characterized by Number of Neurons},
	journal = {Communications in Computational Physics}
}

@article{shen2022optimal,
  title={Optimal approximation rate of ReLU networks in terms of width and depth},
  author={Shen, Zuowei and Yang, Haizhao and Zhang, Shijun},
  journal={Journal de Math{\'e}matiques Pures et Appliqu{\'e}es},
  volume={157},
  pages={101--135},
  year={2022},
  publisher={Elsevier}
}

@book{gyorfi2002distribution,
  title={A Distribution-Free Theory of Nonparametric Regression},
  author={Gy{\"o}rfi, L{\'a}szl{\'o} and Kohler, Michael and Krzyzak, Adam and Walk, Harro and others},
  edition={1},
  year={2002},
  publisher={Springer},
  address={New York}
}

@article{bauer2019deep,
  title={On deep learning as a remedy for the curse of dimensionality in nonparametric regression},
  author={Bauer, Benedikt and Kohler, Michael},
  journal={The Annals of Statistics},
  volume={47},
  number={4},
  pages={2261--2285},
  year={2019},
  //mrnumber={MR3953451}, 
  publisher={Institute of Mathematical Statistics},
  //DOI={10.1214/18-AOS1747},
}

@article{farrell2021deep,
  title={Deep neural networks for estimation and inference},
  author={Farrell, Max H and Liang, Tengyuan and Misra, Sanjog},
  journal={Econometrica},
  volume={89},
  number={1},
  pages={181--213},
  year={2021},
  publisher={Wiley Online Library}
}

@inproceedings{suzuki2021deep,
Author = {Suzuki, Taiji and Nitanda, Atsushi},
Title = {Deep learning is adaptive to intrinsic dimensionality of model
   smoothness in anisotropic Besov space},
Booktitle = {Advances in Neural Information Processing System 34 (NeurIPS 2021)},
\\Series = {Advances in Neural Information Processing Systems},
Year = {2021},
Volume = {34},
ISSN = {1049-5258},
ResearcherID-Numbers = {Nitanda, Atsushi/HSG-6879-2023},
Unique-ID = {WOS:000901616406065},
}

@article{barron1993universal,
  title={Universal approximation bounds for superpositions of a sigmoidal function},
  author={Barron, Andrew R},
  journal={IEEE Transactions on Information theory},
  volume={39},
  number={3},
  pages={930--945},
  year={1993},
  publisher={IEEE},
  //mrnumber={MR1237720},
  //DOI = {10.1109/18.256500},
}

@article{schmidt2019deep,
  title={Deep ReLU network approximation of functions on a manifold},
  author={Schmidt-Hieber, Johannes},
  journal={ArXiv preprint. arXiv:1908.00695},
  year={2019}
}

@article{nakada2020adaptive,
  title={Adaptive approximation and generalization of deep neural network with intrinsic dimensionality},
  author={Nakada, Ryumei and Imaizumi, Masaaki},
  journal={Journal of Machine Learning Research},
  volume={21},
  number={174},
  pages={1--38},
  year={2020}
}

@article{koltchinskii2006local,
  title={Local Rademacher complexities and oracle inequalities in risk minimization},
  author={Koltchinskii, Vladimir},
  journal={The Annals of Statistics},
  volume={34},
  number={6},
  pages={2593--2656},
  year={2006},
  //mrnumber={MR2329442},
  publisher={Institute of Mathematical Statistics}
}

@article{wang2021estimation,
  title={Estimation of the mean function of functional data via deep neural networks},
  author={Wang, Shuoyang and Cao, Guanqun and Shang, Zuofeng and Alzheimer's Disease Neuroimaging Initiative},
  journal={Stat},
  volume={10},
  number={1},
  pages={e393},
  year={2021},
  publisher={Wiley Online Library}
}

@article{shao2022intrinsic,
  title={Intrinsic Riemannian functional data analysis for sparse longitudinal observations},
  author={Shao, Lingxuan and Lin, Zhenhua and Yao, Fang},
  journal={The Annals of Statistics},
  volume={50},
  number={3},
  pages={1696--1721},
  year={2022},
  publisher={Institute of Mathematical Statistics}
}

@article{yang2022online,
Author = {Yang, Ying and Yao, Fang},
Title = {Online Estimation for Functional Data},
journal = {Journal of the American Statistical Association},
Year = {2023},
Volume = {118},
Number = {543},
Pages = {1630-1644},
Month = {JUL 3},
DOI = {10.1080/01621459.2021.2002158},
EarlyAccessDate = {DEC 2021},
ISSN = {0162-1459},
EISSN = {1537-274X},
ResearcherID-Numbers = {yan, jie/HNJ-0097-2023},
Unique-ID = {WOS:000739192900001},
}

@book{hedeker2006longitudinal,
  title={Longitudinal Data Analysis},
  author={Hedeker, Donald and Gibbons, Robert D},
  year={2006},
  publisher={Wiley-Interscience},
  address={Hoboken, NJ},
  edition={1}
}

@book{baltagi2008econometric,
  title={Econometric Analysis of Panel Data},
  author={Badi H. Baltagi},
  edition = {6th},
  year={2021},
  publisher={Springer Cham},

  isbn={978-3-030-53953-5},
  //DOI = {10.1007/978-3-030-53953-5},
}

@book{hsing2015theoretical,
  title={Theoretical foundations of functional data analysis, with an introduction to linear operators},
  author={Hsing, Tailen and Eubank, Randall},
  year={2015},
  publisher={John Wiley \& Sons},
  address={Chichester, West Sussex},
  edition=1
}

@article{deepcox,
author = {Qixian Zhong and Jonas M\"uller and Jane-Ling Wang},
title = {{Deep learning for the partially linear Cox model}},
volume = {50},
journal = {The Annals of Statistics},
number = {3},
publisher = {Institute of Mathematical Statistics},
pages = {1348 -- 1375},
keywords = {Censored data, minimax estimation, neural network, partial likelihood, Semiparametric efficiency, Survival analysis},
year = {2022},
//DOI = {10.1214/21-AOS2153},
URL = {https://doi.org/10.1214/21-AOS2153}
}

@inproceedings{deephazard,
 author = {Zhong, Qixian and M\"uller, Jonas W and Wang, Jane-Ling},
 booktitle = {Advances in Neural Information Processing Systems},
 pages = {15111--15124},
 publisher = {Curran Associates, Inc.},
 title = {Deep Extended Hazard Models for Survival Analysis},
 url = {https://proceedings.neurips.cc/paper/2021/file/7f6caf1f0ba788cd7953d817724c2b6e-Paper.pdf},
 volume = {34},
 year = {2021}
}

@article{zeger1988models,
  title={Models for longitudinal data: a generalized estimating equation approach},
  author={Zeger, Scott L and Liang, Kung-Yee and Albert, Paul S},
  journal={Biometrics},
  pages={1049--1060},
  year={1988},
  publisher={JSTOR}
}

@article{stone1980,
author = {Charles J. Stone},
title = {{Optimal Rates of Convergence for Nonparametric Estimators}},
volume = {8},
journal = {The Annals of Statistics},
number = {6},
publisher = {Institute of Mathematical Statistics},
pages = {1348 -- 1360},
keywords = {density function, nonparametric estimator, Optimal rate of convergence, regression function},
year = {1980},
//DOI = {10.1214/aos/1176345206},
//URL = {https://doi.org/10.1214/aos/1176345206}
}

@book{weiss2005modeling,
  title={Modeling Longitudinal Data},
  author={Weiss, Robert E},
  year={2005},
  publisher={Springer-Verlag},
  address={New York},
  edition={1}
}

@incollection{chamberlain1984panel,
title = { Panel data},
booktitle = {Handbook of Econometrics},
publisher = {Elsevier},
volume = {2},
year = {1984},
issn = {1573-4412},
//pages = {1247-1318},
//DOI = {https://doi.org/10.1016/S1573-4412(84)02014-6},
url = {https://www.sciencedirect.com/science/article/pii/S1573441284020146},
author = {Gary Chamberlain},
abstract = {Publisher Summary
This chapter discusses the models that are static conditional on a latent variable. The panel aspect of the data has been primarily used to control for the latent variable. Much work needs to be done on models that incorporate uncertainty and interesting dynamics. Exploiting the martingale implications of time additive utility seems fruitful. There is, however, a potentially important distinction between time averages and cross-section averages. A time average of forecast errors over T periods should converge to zero as T→ ∞. But an average of forecast errors across N individuals surely need not converge to zero as N→ ∞; there is a common component in those errors, due to economy-wide innovations. The same point applies when considering covariances of forecast errors with variables that are in the agent's information sets. If those conditioning variables are discrete, one can think of averaging over subsets of the forecast errors; as T→ ∞, these averages should converge to zero but not necessarily as N → ∞.}
}

@book{ferraty2006nonparametric,
  title={Nonparametric Functional Data Analysis: Theory and Practice},
  author={Ferraty, Fr{\'e}d{\'e}ric and Vieu, Philippe},
  volume={76},
  year={2006},
  publisher={Springer-Verlag},
  address={New York},
  edition=1
}

@article{fan2022noise,
author = {Jianqing Fan and Yihong Gu and Wen-Xin Zhou},
title = {{How do noise tails impact on deep ReLU networks?}},
volume = {52},
journal = {The Annals of Statistics},
number = {4},
publisher = {Institute of Mathematical Statistics},
pages = {1845 -- 1871},
year = {2024},
doi = {10.1214/24-AOS2428},
URL = {https://doi.org/10.1214/24-AOS2428}
}

@article{fan2022factor,
author = {Jianqing Fan and Yihong Gu},
title = {Factor Augmented Sparse Throughput Deep ReLU Neural Networks for High Dimensional Regression},
journal = {Journal of the American Statistical Association},
volume = {0},
number = {0},
//pages = {1--15},
year = {2023},
publisher = {ASA Website},
doi = {10.1080/01621459.2023.2271605},
URL = {https://doi.org/10.1080/01621459.2023.2271605},
eprint = { https://doi.org/10.1080/01621459.2023.2271605}
}

@article{Bartlett2019nearly,
  author  = {Peter L. Bartlett and Nick Harvey and Christopher Liaw and Abbas Mehrabian},
  title   = {Nearly-tight VC-dimension and Pseudodimension Bounds for Piecewise Linear Neural Networks},
  journal = {Journal of Machine Learning Research},
  year    = {2019},
  volume  = {20},
  number  = {63},
  pages   = {1--17},
  //mrnumber={MR3960917},
  url     = {http://jmlr.org/papers/v20/17-612.html}
}

@article{carlsson2009topology,
  title={Topology and data},
  author={Carlsson, Gunnar},
  journal={Bulletin of the American Mathematical Society},
  volume={46},
  number={2},
  pages={255--308},
  year={2009},
  //mrnumber={MR2476414},
  //DOI={10.1090/S0273-0979-09-01249-X}
}

@article{Anirban2014anisot,
author = {Anirban Bhattacharya and Debdeep Pati and David Dunson},
title = {{Anisotropic function estimation using multi-bandwidth Gaussian processes}},
volume = {42},
journal = {The Annals of Statistics},
number = {1},
publisher = {Institute of Mathematical Statistics},
pages = {352 -- 381},
keywords = {adaptive, anisotropic, Bayesian nonparametrics, Function estimation, Gaussian process, rate of convergence},
year = {2014},
//doi = {10.1214/13-AOS1192},
//URL = {https://doi.org/10.1214/13-AOS1192}
}

@article{covnet,
author = {Sarkar, Soham and Panaretos, Victor M.},
title = {CovNet: Covariance networks for functional data on multidimensional domains},
journal = {Journal of the Royal Statistical Society: Series B (Statistical Methodology)},
volume = {84},
number = {5},
pages = {1785-1820},
keywords = {deep learning, FDA, neural network, non-parametric model, random field, universal approximation},
//doi = {https://doi.org/10.1111/rssb.12551},
//url = {https://rss.onlinelibrary.wiley.com/doi/abs/10.1111/rssb.12551},
//eprint = {https://rss.onlinelibrary.wiley.com/doi/pdf/10.1111/rssb.12551},
abstract = {Abstract Covariance estimation is ubiquitous in functional data analysis. Yet, the case of functional observations over multidimensional domains introduces computational and statistical challenges, rendering the standard methods effectively inapplicable. To address this problem, we introduce Covariance Networks (CovNet) as a modelling and estimation tool. The CovNet model is universal—it can be used to approximate any covariance up to desired precision. Moreover, the model can be fitted efficiently to the data and its neural network architecture allows us to employ modern computational tools in the implementation. The CovNet model also admits a closed-form eigendecomposition, which can be computed efficiently, without constructing the covariance itself. This facilitates easy storage and subsequent manipulation of a covariance in the context of the CovNet. We establish consistency of the proposed estimator and derive its rate of convergence. The usefulness of the proposed method is demonstrated via an extensive simulation study and an application to resting state functional magnetic resonance imaging data.},
year = {2022}
}

@article{bhattacharya2023deep,
      title={Deep Neural Networks for Nonparametric Interaction Models with Diverging Dimension}, 
      author={Sohom Bhattacharya and Jianqing Fan and Debarghya Mukherjee},
      journal={ArXiv preprint. arXiv:2302.05851},
      year={2023}
}

@article{rice1983smoothing,
  title={Smoothing splines: regression, derivatives and deconvolution},
  author={Rice, John and Rosenblatt, Murray},
  journal={The Annals of Statistics},
Year = {1983},
Volume = {11},
Number = {1},
Pages = {141-156},
DOI = {10.1214/aos/1176346065},
ISSN = {0090-5364},
Unique-ID = {WOS:A1983QE68600015},
}

@article{fan1992variable,
  title={Variable bandwidth and local linear regression smoothers},
  author={Fan, Jianqing and Gijbels, Irene},
  journal={The Annals of Statistics},
  pages={2008--2036},
  year={1992},
Volume = {20},
Number = {4},
Month = {DEC},
  publisher={JSTOR}
}

@article{brumback1998smoothing,
  title={Smoothing spline models for the analysis of nested and crossed samples of curves},
  author={Brumback, Babette A and Rice, John A},
  journal={Journal of the American Statistical Association},
  volume={93},
  number={443},
  pages={961--976},
  year={1998},
  publisher={Taylor \& Francis}
}

@article{fan2021selective,
  title={A selective overview of deep learning},
  author={Fan, Jianqing and Ma, Cong and Zhong, Yiqiao},
  journal={Statistical science: a review journal of the Institute of Mathematical Statistics},
  volume={36},
  number={2},
  pages={264},
  year={2021},
  publisher={NIH Public Access}
}

@article{shaham2018provable,
  title={Provable approximation properties for deep neural networks},
  author={Shaham, Uri and Cloninger, Alexander and Coifman, Ronald R},
  journal={Applied and Computational Harmonic Analysis},
  volume={44},
  number={3},
  pages={537--557},
  year={2018},
  publisher={Elsevier}
}

@article{cloninger2021deeparxiv,
Author = {Cloninger, Alexander and Klock, Timo},
Title = {A deep network construction that adapts to intrinsic dimensionality
   beyond the domain},
Journal = {Neural Networks},
Year = {2021},
Volume = {141},
Pages = {404-419},
Month = {SEP},
DOI = {10.1016/j.neunet.2021.06.004},
EarlyAccessDate = {JUN 2021},
ISSN = {0893-6080},
EISSN = {1879-2782},
ORCID-Numbers = {Klock, Timo/0000-0002-0122-3017
   Cloninger, Alexander/0000-0002-1423-9624},
Unique-ID = {WOS:000681162400017},
}

@article{chen2019efficient,
  title={Nonparametric regression on low-dimensional manifolds using deep ReLU networks: Function approximation and statistical recovery},
  author={Chen, Minshuo and Jiang, Haoming and Liao, Wenjing and Zhao, Tuo},
  journal={Information and Inference: A Journal of the IMA},
  volume={11},
  number={4},
  pages={1203--1253},
  year={2022},
  publisher={Oxford University Press}
}

@article{HANG2021337,
title = {Optimal learning with anisotropic Gaussian SVMs},
journal = {Applied and Computational Harmonic Analysis},
volume = {55},
pages = {337-367},
year = {2021},
issn = {1063-5203},
author = {Hanyuan Hang and Ingo Steinwart}
}

@article{bickel2007local,
  title={Local polynomial regression on unknown manifolds},
  author={Bickel, Peter J and Li, Bo},
  journal={Lecture Notes-Monograph Series},
  pages={177--186},
  volume = {54},
  year={2007},
  publisher={JSTOR}
}

@article{aswani2011regression,
    Author = {Aswani, Anil and Bickel, Peter and Tomlin, Claire},
    Title = {REGRESSION ON MANIFOLDS: ESTIMATION OF THE EXTERIOR DERIVATIVE},
    Journal = {Annals of Statistics},
    Year = {2011},
    Volume = {39},
    Number = {1},
    Pages = {48-81},
    Month = {FEB},
    DOI = {10.1214/10-AOS823},
    ISSN = {0090-5364},
    ORCID-Numbers = {Aswani, Anil/0000-0001-5777-7185},
    Unique-ID = {WOS:000288183800002},
}

@article{cheng2013local,
  title={Local linear regression on manifolds and its geometric interpretation},
  author={Cheng, Ming-Yen and Wu, Hau-tieng},
  journal={Journal of the American Statistical Association},
  volume={108},
  number={504},
  pages={1421--1434},
  year={2013},
  publisher={Taylor \& Francis}
}

@article{zhang2002wavelet,
  title={Wavelet threshold estimation of a regression function with random design},
  author={Zhang, Shuanglin and Wong, Man-Yu and Zheng, Zhongguo},
  journal={Journal of multivariate analysis},
  volume={80},
  number={2},
  pages={256--284},
  year={2002},
  publisher={Elsevier}
}

@article{horowitz2007,
author = {Joel L. Horowitz and Enno Mammen},
title = {{Rate-optimal estimation for a general class of nonparametric regression models with unknown link functions}},
volume = {35},
journal = {The Annals of Statistics},
number = {6},
publisher = {Institute of Mathematical Statistics},
pages = {2589 -- 2619},
year = {2007},
doi = {10.1214/009053607000000415},
URL = {https://doi.org/10.1214/009053607000000415}
}

@article{zhang2024nonparametric,
      title={Nonparametric Classification on Low Dimensional Manifolds using Overparameterized Convolutional Residual Networks}, 
      author={Kaiqi Zhang and Zixuan Zhang and Minshuo Chen and Yuma Takeda and Mengdi Wang and Tuo Zhao and Yu-Xiang Wang},
      year={2024},
      journal={ArXiv preprint. arXiv:2307.01649}
}

@inproceedings{zhang2023effective,
  title={Effective minkowski dimension of deep nonparametric regression: function approximation and statistical theories},
  author={Zhang, Zixuan and Chen, Minshuo and Wang, Mengdi and Liao, Wenjing and Zhao, Tuo},
  booktitle={International Conference on Machine Learning},
  pages={40911--40931},
  year={2023},
  organization={PMLR}
}

@inproceedings{liu2021besov,
  title={Besov function approximation and binary classification on low-dimensional manifolds using convolutional residual networks},
  author={Liu, Hao and Chen, Minshuo and Zhao, Tuo and Liao, Wenjing},
  booktitle={International Conference on Machine Learning},
  pages={6770--6780},
  year={2021},
  organization={PMLR}
}

@article{xu2023sample,
Author = {Xu, Zhenghao and Ji, Xiang and Chen, Minshuo and Wang, Mengdi and Zhao,
   Tuo},
Title = {Sample Complexity of Neural Policy Mirror Descent for Policy
   Optimization on Low-Dimensional Manifolds},
Journal = {Journal of Machine Learning Research},
Year = {2024},
Volume = {25},
Article-Number = {226},
ISSN = {1532-4435},
ResearcherID-Numbers = {Zhao, Tuo/AAG-7648-2021},
Unique-ID = {WOS:001293872000001},
}

@article{dahal2022deep,
  title={On deep generative models for approximation and estimation of distributions on manifolds},
  author={Dahal, Biraj and Havrilla, Alexander and Chen, Minshuo and Zhao, Tuo and Liao, Wenjing},
  journal={Advances in Neural Information Processing Systems},
  volume={35},
  pages={10615--10628},
  year={2022}
}

@article{zheng2024dynamic,
  title={Dynamic synthetic control method for evaluating treatment effects in auto-regressive processes},
  author={Zheng, Xiangyu and Chen, Song Xi},
  journal={Journal of the Royal Statistical Society Series B: Statistical Methodology},
  volume={86},
  number={1},
  pages={155--176},
  year={2024},
  publisher={Oxford University Press US}
}

@article{zou2022estimation,
  title={Estimation of low rank high-dimensional multivariate linear models for multi-response data},
  author={Zou, Changliang and Ke, Yuan and Zhang, Wenyang},
  journal={Journal of the American Statistical Association},
  volume={117},
  number={538},
  pages={693--703},
  year={2022},
  publisher={Taylor \& Francis}
}

@article{chen2017detecting,
  title={Detecting the causality influence of individual meteorological factors on local PM2. 5 concentration in the Jing-Jin-Ji region},
  author={Chen, Ziyue and Cai, Jun and Gao, Bingbo and Xu, Bing and Dai, Shuang and He, Bin and Xie, Xiaoming},
  journal={Scientific Reports},
  volume={7},
  number={1},
  pages={40735},
  year={2017},
  publisher={Nature Publishing Group UK London}
}

@Manual{gampackage,
  title = {gam: Generalized Additive Models},
  author = {Trevor Hastie},
  year = {2023},
  note = {R package version 1.22-3},
  url = {https://CRAN.R-project.org/package=gam},
}

@Manual{PLSiMCpp,
  title = {PLSiMCpp: Methods for Partial Linear Single Index Model},
  author = {Wu, Shunyao and Zhang, Qi and Li, Zhiruo and Liang, Hua},
  year = {2022},
  note = {R package version 1.0.4},
  url = {https://CRAN.R-project.org/package=PLSiMCpp},
}

\end{document}